\documentclass[journal]{IEEEtran}

\ifCLASSINFOpdf
\else
\fi

\usepackage{xcolor}
\usepackage{cite}
\usepackage{amsmath}
\usepackage{amssymb}
\usepackage{amsfonts}
\usepackage{algorithmic}
\usepackage{enumerate}
\usepackage{graphicx}
\usepackage{subfigure}
\usepackage{url}
\usepackage{float}
\usepackage{textcomp}
\usepackage{etoolbox}
\usepackage{epstopdf}
\usepackage{amsthm}
\newtheorem{theorem}{Theorem}
\newtheorem{remark}{Remark}
\usepackage{quoting}

\usepackage{color}

\newtheorem{assum}{Assumption}

\newcommand{\sgn}{\operatorname{sgn}}

\begin{document}
\title{Extremum and Nash Equilibrium Seeking with Delays and PDE{s}: Designs $\&$ Applications} 

	
	\author{Tiago Roux Oliveira, Miroslav Krsti\'{c}, and Tamer Ba{\c s}ar. 
		\thanks{T. R. Oliveira is with the Department of Electronics and Telecommunication Engineering,
			 State University of Rio de Janeiro (UERJ), Brazil
			(e-mail: tiagoroux@uerj.br).}
        \thanks{M. Krsti\'{c} is with the  Department of Mechanical and Aerospace Engineering, University of California - San Diego, USA (e-mail: krstic@ucsd.edu).}
        \thanks{T. Ba{\c s}ar is with the Department of Electrical and Computer Engineering, University of Illinois Urbana-Champaign, Champaign, IL, USA (e-mail: basar1@illinois.edu).}
		\thanks{Corresponding author: T. R. Oliveira. 
            }
		\thanks{This work was supported by the Brazilian agencies CNPq, CAPES and FAPERJ.}}

\markboth{Submitted to IEEE Control Systems Magazine (Special Issue: Into the Second Century of
Extremum Seeking Control)}%
{Oliveira \MakeLowercase{\textit{et al.}}: Extremum and Nash Equilibrium Seeking with Delays and PDE{s}}

	\maketitle

	\begin{abstract}
           The development of extremum seeking (ES) has progressed, over the past hundred years,
from static maps, to finite-dimensional dynamic systems, to networks of static and dynamic agents.
Extensions from ODE dynamics to maps and agents that incorporate delays or even partial differential
equations (PDEs) is the next natural step in that progression through ascending research challenges.
This paper reviews results on algorithm design and theory of ES for such infinite-dimensional systems.
Both 
hyperbolic and parabolic dynamics are presented:
delays or transport equations, heat-dominated equation, wave equations, and reaction-advection-diffusion equations.
Nash equilibrium seeking (NES) methods are introduced for noncooperative game scenarios of the model-free kind and then specialized to single-agent optimization. 
Even heterogeneous PDE games, such as a duopoly with one parabolic and one hyperbolic agent, are
considered. Several engineering applications are touched upon for illustration, including flow-traffic control for urban mobility, oil-drilling systems, deep-sea
cable-actuated source seeking, additive manufacturing modeled by the Stefan PDE, biological reactors, light-source seeking with flexible-beam structures, and 
neuromuscular electrical stimulation.
	\end{abstract}

	\begin{IEEEkeywords}
		 extremum seeking; Nash equilibrium seeking; partial differential equations; time delays; infinite-dimensional systems, game theory.
	\end{IEEEkeywords}

	\IEEEpeerreviewmaketitle

\section{Introduction}


Classical extremum seeking (ES) dealts with model-free online optimization problems in which, regardless of whether there is a single input or multiple inputs, only a single payoff is being maximized (or cost minimized). Optimization problems become  more interesting when there are multiple payoffs, being maximized by the respective multiple inputs/actors/players. When the different payoffs are distinct functions of the inputs of different players, it is clear that conflicting objectives may arise among the players. This scenario leads to what is called a noncooperative game. In a noncooperative game, no player can improve its payoff over the payoff that results from applying a certain optimal action. The collection of such optimal actions by the players is called a Nash equilibrium \cite{N:1951,Basar:1999}. 

The realization that extremum seeking is not only applicable to single-player optimization (possibly with a vector of inputs maximizing a single payoff), in an online model-free setting, but also to multiple competing players in the noncooperative game setting, emerged in \cite{FKB:2012}. Each player employing an ES algorithm maximizes its own payoff, irrespective of what the other players' actions are. It is proven in \cite{FKB:2012} that if all the players employ ES algorithms they collectively converge to a Nash equilibrium. In other words, each of the players finds its optimal strategy, in an online fashion, even though they do not know the analytical forms of the payoff functions (neither the other players' nor their own) and neither have access to the actions applied by the other players nor to the payoffs achieved by the other players. Recent versions of the Nash equilibrium seeking algorithm originally proposed in \cite{FKB:2012} can be found in \cite{PKB:2023}, where fixed-time convergence is incorporated in time-varying networks. 

Delays and PDEs \cite{K:2009,K:2010} in a game context arise in applications like network virtualization, software defined networks, cloud computing, the Internet of Things, context-aware networks, green communications, and security \cite{Basar:2019,Sastry:2013,AB:2011}. In particular, differential games with delays are dealt with (in a partly or fully model-based fashion) in  \cite{C:1972,MS:1974,KT:1984,EH:1989,GS:2001,P:2015,CFMS:2018}. PDE dynamics arise in the Black-Scholes model of behavior in financial markets \cite{R:1975,LM:2014}. Hence, strong motivation exists for designs that ensure convergence to Nash equilibria in the presence of delays and PDEs.

The development of ES has progressed from static maps to finite-dimensional dynamic systems,
to networks of static and dynamic players. Extensions from ODE dynamics to maps and players that incorporate delays or even PDEs was the next natural step in that progression through ascending research challenges \cite{OK:2022}. In \cite{Tiago} and \cite{FKKO:2018,TAC:2020}, we dealt with classical ES with single payoff functions, where the actuation dynamics are governed by hyperbolic and parabolic PDEs, respectively. In \cite{JOTA:2021}, we have considered noncooperative games where players act through arbitrarily long delays or transport hyperbolic PDEs. The delays may be distinct and, in general, each player knows only the length of its own delay. In order to compensate for distinct delays at the inputs of the players, we have employed predictor feedback.


In the first half of this paper, we  consider noncooperative games where players act through a more complex and complete scenario of PDE dynamics. 

First, following a natural transition from the Nash equilibrium seeking (NES) design through delays to design through parabolic PDEs, we progress to developing a non-model based strategy for locally stable convergence to Nash equilibria in quadratic noncooperative games with player actions subject to diffusion (heat) PDE dynamics with distinct diffusion coefficients and each player having access only to its own payoff value. The proposed approach employs extremum seeking, with sinusoidal perturbation signals employed to estimate the Gradient (first derivative) and Hessian (second derivative) of unknown quadratic functions. 
In order to compensate for distinct diffusion processes in the inputs of the $N$-players, we employ boundary control with averaging-based estimates. 
%

Second, when players are competing in a noncooperative game, there is no reason to assume that players are pursuing even remotely related objectives. One could be maximizing its profit while another one maximizes some social good. Likewise, there is no reason to assume that different players are subject to the same type of physics. One might be propagating an input through some social opinion dynamics while another one may be propagating it through some epidemiological dynamics. This means that, while first we provide useful designs and guarantees of Nash equilibrium attainment in the presence of delays or heat dynamics, \textit{i.e.}, in the presence of hyperbolic and parabolic PDE dynamics, it is of interest to see whether, when different players, unaware of the competing optimization pursuits in different physical domains, are employing ES algorithms with compensation of their own specific PDEs, will be able to achieve the Nash-optimality under the interference of other players who are compensating their own PDE dynamics, which are from different classes.


In this sense, one could formulate a problem of a noncooperative game with dynamics from different PDE classes (hyperbolic, parabolic, Korteweg-de-Vriess, Schroedinger, Kuramoto Sivashinsky, etc.) 
but we do not pursue all such possibilities here. We only demonstrate the achievement of convergence to a Nash equilibrium in a game with two ``heterogeneous'' players (a duopoly) where one player is compensating a transport PDE and the other one a heat PDE. While not coupled directly, the coupling of the players in the payoff functions, and the coupling that results from the ES algorithms acting at the boundary conditions of the PDEs, results in a coupled heterogeneous pair of PDEs in the closed-loop system, which happens also to be time-varying and nonlinear.
\


\textcolor{black}{
In both scenarios (homogeneous or heterogeneous games), we apply an adequate small-gain analysis for the resulting Input-to-State Stable (ISS) parabolic PDE-ODE loop, as well as averaging theory in infinite dimensions, due to the infinite-dimensional state of the heat PDEs and delays, in order to obtain local convergence  to a small neighborhood of the Nash equilibrium. We quantify the size of these residual sets for this distinct class of parabolic-hyperbolic PDEs.}

The material in the present paper and the related  results \cite{ORKB:2020, ACC:2021,TempoAward:2021}, presented only in conferences but not published in journals, is the first instance of noncooperative games being tackled in a model-free fashion with extremum seeking in the presence of PDE dynamics, represented by heat PDEs and delays. The employment of PDE motion planning (for obtaining the probing signals) combined with the PDE-ISS formulation as well as the distributed boundary control feedback design and the closed-loop analysis via small-gain theorems for PDEs are the main bits of new knowledge that a reader already acquainted with extremum seeking will gain from this paper. 

We start the second part of the paper by highlighting the relation between NES and classical ES for infinite-dimensional systems. In this discussion, classes of PDEs more general than transport and heat equations are addressed. 

The second half of the paper is dedicated to a number of engineering applications: is presented:  traffic control, oil drilling control, deep-sea cable-actuated control, additive manufacturing, bioreactors, light-source seeking with cantilever beams, and neuromuscular electrical stimulation under delays. These applications may draw interest based either on their physical nature or the mathematical nature.  

Physically, the applications selected, to introduce the reader to some of the technological possibilities of ES, are grouped as follows. First, we present an application representative of flow: traffic. Next, we show applications entailing mechanical vibrations: drill string instabilities in oil drilling and oscillation of cables in undersea manipulation and construction. Third, we turn to thermal diffusion in additive manufacturing and the diffusion of chemical or biological matter in tubular reactors. Fourth, in a light-seeking application, we consider cantilever beams (vibrating systems qualitatively different from the torsional vibrations in drill strings and translational vibrations in cables). Fifth, we consider neuromuscular electrical stimulation (NMES) in patient rehabilitation. 

The applications featured in the paper are not selected only for their physical variety and representativeness, but also for their mathematical representativeness. The traffic and NMES applications feature a pure transport PDE or delay. The drilling and cable applications feature second-order hyperbolic PDEs. The additive manufacturing and bioreactor applications feature parabolic PDEs. Through these applications, the reader is exposed to the full spectrum of basic PDE types in the presence of which ES can be successfully pursued. 

The paper is organized as follows. Section~``$N$-Player Game with Quadratic Payoffs for Nash Equilibrium Seeking'' introduces the standard assumptions (such as, diagonal dominance for the Hessian matrix) for noncooperative games with $N$ players and quadratic payoff functions. While the setup taking into account delays is shown in Section~``Noncooperative Scenario with Delays'', we consider the inclusion of heat PDEs in noncooperative games in the Section~``Noncooperative Games with Heat PDEs''. In these two sections, we consider delays and heat PDEs separately, each player handling a distinct heat PDE or transport PDE. The extension from homogeneous noncooperative games to heterogeneous ones is carried out in the Section~``Noncooperative Duopoly with Heterogeneous Transport-Heat PDE Dynamics''. In that section, we seek Nash equilibria in a quadratic noncooperative duopoly game ($N=2$) with player actions subject to concomitantly transport-heat PDE dynamics---one player compensating for a delay (transport PDE) and the other one for a heat (diffusion) PDE---and each player having access only to its own payoff value. The stability analysis for every specific case is conducted in the corresponding Sections~``Heat PDEs Case'', ``Delay Case'' and ``Stability Analysis'', with the proofs of the main theorems made available in the Appendix. We illustrate the theoretical results numerically on an example combining hyperbolic and parabolic dynamics in a 2-player setting, given in the Section~``Simulations''. The Section ``From Nash Equilibrium to Extremum Seeking'' brings further results on ES for infinite-dimensional systems governed by a wider class of PDE systems. A series of application sections follows, covering traffic control, oil drilling control, deep-sea cable-actuated control, additive manufacturing, bioreactors, light-source seeking with cantilever beams, and neuromuscular electrical stimulation under delays; showing the vast space of possibilities of designing ES algorithms for real-world applications. Concluding remarks are included in the Section~``Conclusion.'' 
Notation and terminology used in the paper are listed in ``Notation, Norms and Terminology''.

\section{$N$-Player Game with Quadratic Payoffs for Nash Equilibrium Seeking}
\label{doritos_section2}

Game theory provides an important framework for mathematical modeling and analysis of scenarios involving different players where there is coupling in their actions, in the sense that their respective outcomes (outputs) $y_{i}(t) \in \mathbb{R}$ do not depend exclusively on their own actions/strategies (input signals) $\theta_{i}(t)\in \mathbb{R}$, with $i=1 ,\ldots, N$, but at least on a subset of others'. Moreover, defining $\theta := [\theta_1, \ldots, \theta_N]^T$, each player's payoff function $J_{i}(\theta) : \mathbb{R}^{N} \to \mathbb{R}$ depends on the action $\theta_j$ of at least one other Player $j$, $j\not= i$. An $N$-tuple of actions, $\theta^*$, is said to be in Nash equilibrium, if no Player $i$ can improve its payoff  by unilaterally deviating from $\theta_i^*$, this being so for all $i$ \cite{Basar:1999}. 

Consider games where the payoff function of each player is quadratic, 
expressed as a strictly concave\footnote{By strict concavity, we mean $J_i(\theta)$ is strictly concave in $\theta_i$ for all $\theta_{-i}$, this being so for each $i=1,\ldots, N$.} combination of their actions 
\begin{align}
J_{i}(\theta(t))=&\frac{1}{2}\sum_{j=1}^{N}\sum_{k=1}^{N}\epsilon_{jk}^{i}H_{jk}^{i}\theta_{j}(t)\theta_{k}(t) 
+\sum_{j=1}^{N}h_{j}^{i}\theta_{j}(t)+c_{i}\,, \label{eq:Ji}
\end{align}  
where $\theta_{j}(t) \!\in\! \mathbb{R}$ is the decision variable (action) of Player $j$, 
$H_{jk}^{i}$, $h_{j}^{i}$, $c_{i} \!\in\! \mathbb{R}$ are constants, $H_{ii}^{i}\!<\!0$, $H_{jk}^{i}\!=\!H_{kj}^{i}$ and $\epsilon_{jk}^{i}\!=\!\epsilon_{kj}^{i}\!>\!0$, $\forall i,j,k$. 

Quadratic payoff functions, of the type above, are of particular interest in game theory, first because they constitute second-order approximations to other types of non-quadratic payoff functions, and second because they are analytically tractable, leading in general to closed-form equilibrium solutions which provide insight into the properties and features of the equilibrium solution concept under consideration \cite{Basar:1999}. 

For the sake of completeness, we provide here in mathematical terms, the definition of a Nash equilibrium $\theta^*=[\theta^*_1\,, \ldots\,,\theta_N^*]^T$ in an $N$-player game:
\begin{align} \label{Nashcu}
J_i(\theta_i^*\,,\theta_{-i}^*) \!\geq\! J_i(\theta_i\,,\theta_{-i}^*)\,, \quad 
\forall \theta_i \in \mathcal{U}_i\,, \quad i \!\in\! \{1\,, \ldots\,,N\}\,,
\end{align}
where $J_i$ is the payoff function of player $i$, the variable $\theta_i$ corresponds to its action, while $\mathcal{U}_i$ is its action set and $\theta_{-i}$ denotes the actions of the other players. Hence, no player has an incentive to unilaterally deviate its action from $\theta^*$. In a duopoly example ($N=2$), to be considered later, $\mathcal{U}_1=\mathcal{U}_2=\mathbb{R}$, where $\mathbb{R}$ denotes the set of real numbers.

In order to determine the Nash equilibrium solution in strictly concave quadratic games with $N$ players, where each action set is the entire real line, one should differentiate $J_{i}$ with respect to $\theta_{i}(t) \,, \forall i=1 ,\ldots, N$, setting the resulting expressions equal to zero, and solving the set of equations thus obtained.

This set of equations, which also provides a sufficient condition due to the strict concavity, is
\begin{align}
\sum_{j =1}^{N}\epsilon_{ij}^{i}H_{ij}^{i}\theta_{j}^{*}+h_{i}^{i}=0\,,\quad i=1 ,\ldots, N\,, \label{eq:NE_v0}
\end{align}
which can be written in compact form as 
\begin{equation}
\begin{aligned}
\begin{bmatrix}
\epsilon_{11}^{1} H_{11}^{1} & \epsilon_{12}^{1} H_{12}^{1} & \hdots & \epsilon_{1N}^{1} H_{1N}^{1} \\
\epsilon_{21}^{2} H_{21}^{2} & \epsilon_{22}^{2} H_{22}^{2} & \hdots & \epsilon_{2N}^{2} H_{2N}^{2} \\
\vdots                       & \vdots                       &        & \vdots     \\
\epsilon_{N1}^{N} H_{N1}^{N} & \epsilon_{N2}^{N} H_{N2}^{N} & \hdots & \epsilon_{NN}^{N} H_{NN}^{N}   
\end{bmatrix}
\begin{bmatrix}
\theta_{1}^{*} \\
\theta_{2}^{*} \\
\vdots \\
\theta_{N}^{*} 
\end{bmatrix}
=-
\begin{bmatrix}
h_{1}^{1} \\
h_{2}^{2} \\
\vdots    \\
h_{N}^{N}   
\end{bmatrix}
\,.
\end{aligned}
\end{equation}

Defining the Hessian matrix $H$ and vectors $\theta^*$ and $h$ by 
\begin{align}
H&:=
\begin{bmatrix}
\epsilon_{11}^{1} H_{11}^{1} & \epsilon_{12}^{1} H_{12}^{1} & \hdots & \epsilon_{1N}^{1} H_{1N}^{1} \\
\epsilon_{21}^{2} H_{21}^{2} & \epsilon_{22}^{2} H_{22}^{2} & \hdots & \epsilon_{2N}^{2} H_{2N}^{2} \\
\vdots                       & \vdots                       &        & \vdots     \\
\epsilon_{N1}^{N} H_{N1}^{N} & \epsilon_{N2}^{N} H_{N2}^{N} & \hdots & \epsilon_{NN}^{N} H_{NN}^{N}   
\end{bmatrix}
\,, \nonumber \\
\theta^{*}&:=
\begin{bmatrix}
\theta_{1}^{*} \\
\theta_{2}^{*} \\
\vdots \\
\theta_{N}^{*} 
\end{bmatrix}
\,, \quad
h:=
\begin{bmatrix}
h_{1}^{1} \\
h_{2}^{2} \\
\vdots    \\
h_{N}^{N}   
\end{bmatrix}
\,, \label{eq:Htheta*h}
\end{align}
there exists a unique Nash Equilibrium at  $\theta^{*}=-H^{-1}h$, if $H$ is invertible:


\begin{align}
\begin{bmatrix}
\theta_{1}^{*} \\
\theta_{2}^{*} \\
\vdots \\
\theta_{N}^{*} 
\end{bmatrix}
&= -
\begin{bmatrix}
\epsilon_{11}^{1} H_{11}^{1} & \epsilon_{12}^{1} H_{12}^{1} & \hdots & \epsilon_{1N}^{1} H_{1N}^{1} \\
\epsilon_{21}^{2} H_{21}^{2} & \epsilon_{22}^{2} H_{22}^{2} & \hdots & \epsilon_{2N}^{2} H_{2N}^{2} \\
\vdots                       & \vdots                       &        & \vdots     \\
\epsilon_{N1}^{N} H_{N1}^{N} & \epsilon_{N2}^{N} H_{N2}^{N} & \hdots & \epsilon_{NN}^{N} H_{NN}^{N}   
\end{bmatrix}^{-1} \nonumber \\
&\quad \times
\begin{bmatrix}
h_{1}^{1} \\
h_{2}^{2} \\
\vdots    \\
h_{N}^{N}   
\end{bmatrix}. \label{eq:NE_v2}
\end{align}

%
For more details, see \cite[Chapter 4]{Basar:1999}.

In addition to Assumption~\ref{ch.13.Assumption 1.} formulated in \cite{FKB:2012}, we further assume/formalize the following Assumption~\ref{ch13.Assumption 2.} for noncooperative games.

\begin{assum} \label{ch.13.Assumption 1.}
The Hessian matrix $H$ given by \textcolor{black}{\textnormal{(\ref{eq:Htheta*h})}} is strictly diagonal dominant, \textit{i.e.},
\begin{equation}
\sum_{j\neq i}^{N}|\epsilon_{ij}^{i}H_{ij}^{i}| < |\epsilon_{ii}^{i}H_{ii}^{i}|\,, \quad i \in \{1\,, \ldots N\}\,.    
\end{equation}
\end{assum}

\begin{assum} \label{ch13.Assumption 2.}
The parameters $\epsilon_{jk}^{i}$ and $\epsilon_{kj}^{i}$ which appear in the Hessian matrix $H$ given by (\ref{eq:Htheta*h}) satisfy the conditions below: 
\begin{eqnarray}
\epsilon_{ii}^{i}=1\,, \quad \epsilon_{jk}^{i}&=&\epsilon_{kj}^{i}=\epsilon\,, \quad \forall j\neq k\,,
\end{eqnarray}
with $0<\epsilon<1$ in the payoff functions  \textcolor{black}{\textnormal{(\ref{eq:Ji})}}.
\end{assum}

By Assumption~\ref{ch.13.Assumption 1.}, the Nash Equilibrium $\theta^*$ exists and is unique since strictly diagonally dominant matrices are nonsingular by the Levy-Desplanques Theorem \cite{HJ:1985}. 
Assumption~\ref{ch13.Assumption 2.} could be relaxed to consider different values of the coupling parameters $\epsilon_i$ for each Player $i$. However, without loss of generality, we have assumed the same weights for the interconnection channels among the players in order to facilitate the proofs of our theorems, but also to guarantee that the considered noncooperative game is not favoring any specific player.

\section{Noncooperative Scenario with Delays} \label{noncooperative_games}

\begin{figure}[htb!]
\begin{center}
\includegraphics[width=3.0 in]{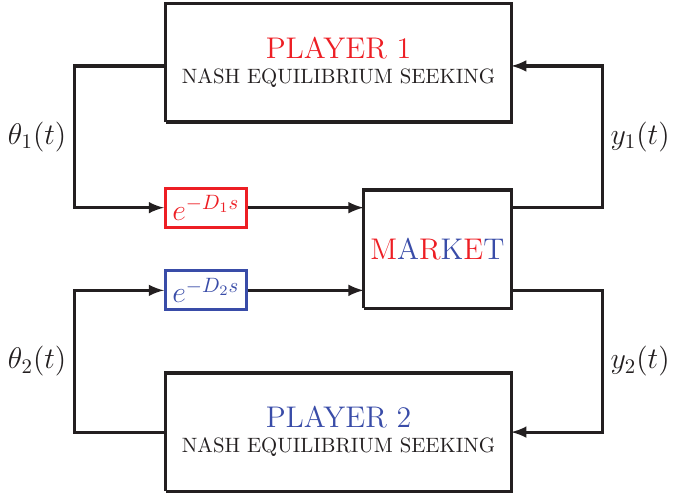}
\end{center}
\caption{Nash equilibrium seeking schemes applied by two players ($N=2$) in a duopoly market structure with delayed players' actions.}
\label{fig:blockDiagram_v2}
\end{figure}

For the sake of completeness and to keep the material sel-contained, we briefly review the case of noncooperative games subject to multiple and distinct delays originally addressed in \cite{JOTA:2021}. In this scenario, the purpose of the extremum seeking is still to estimate the Nash equilibrium vector $\theta^*$, but without allowing any sharing of information among the players. As mentioned earlier, each player only needs to measure the value of its own payoff function described by  

\begin{equation} \label{ajagambiarra_noncooperative}
\begin{split}
y_i(t) &= J_i(\theta(t - D)) \\
&= \frac{1}{2} \sum_{j=1}^{N} \sum_{k=1}^{N} \epsilon_{jk}^{i} H_{jk}^{i} \theta_{j}(t - D_{j}) \theta_{k}(t - D_{k}) \\
& \quad + \sum_{j=1}^{N} h_{j}^{i} \theta_{j}(t - D_{j}) + c_{i}.
\end{split}
\end{equation}
with $J_i$ given by (\ref{eq:Ji}). In this sense, we are able to formulate the closed-loop system in a \textit{decentralized} fashion, where no knowledge about the payoffs $y_{-i}$ or actions $\theta_{-i}$ of the other players is required, as illustrated in Fig.~\ref{fig:blockDiagram_v2}. 
%
%

Without loss of generality, we assume that the inputs have distinct known (constant) delays which are ordered so that
\begin{equation}\label{equ:4_delay_order}
	D = \text{diag}\{D_1,D_2,\cdots, D_N\}, \quad 0\leq D_1 \leq \cdots \leq D_N\,.
\end{equation}
Further, given any $\mathbb{R}^N$-valued signal $f$, we introduce 
\begin{small}
\begin{equation} \label{4_fD}
	f^D(t) \!:=\! f(t-D)\!=\!  
	\begin{bmatrix}
		f_1(t\!-\!D_1) \!&\! f_2(t\!-\!D_2) \!&\! ... \!&\! f_N(t\!-\!D_N)
	\end{bmatrix}^T.
\end{equation} 
\end{small}

\begin{figure}[ht]
\begin{center}
\includegraphics[width=3.25 in]{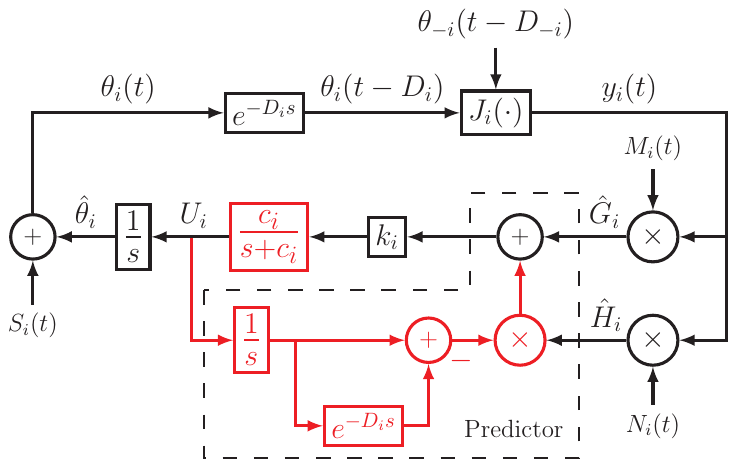}
\end{center}
\caption{Block diagram illustrating the Nash Equilibrium seeking strategy performed for each player. \textcolor{black}{In red color, the predictor feedback used to compensate the individual delay $D_i$ for the noncooperative case. With some abuse of notation, the constants $c_i$ were chosen to denote the parameters of the filters $c_i/(s+c_i)$, but they are not necessarily the same constants which appear in the payoff functions given in (\ref{eq:Ji}).}}
\label{fig:blockDiagram_v3}
\end{figure}

Fig,~\ref{fig:blockDiagram_v3} contains a schematic diagram that summarizes the proposed Nash Equilibrium policy for each $i$-th player where its output is given by (\ref{ajagambiarra_noncooperative}), 
where the vector $\theta_{-i}(t-D_{-i})$ in Fig.~\ref{fig:blockDiagram_v3} represents the delayed actions of all other players.

The additive-multiplicative dither signals $S_i(t)$, $M_i(t)$ are
\begin{align}
S_{i}(t)&=a_{i}\sin(\omega_{i}(t+D_{i})) \,, \label{eq:Si} \\
M_{i}(t)&=\frac{2}{a_{i}}\sin(\omega_{i}t)\,, \label{eq:Mi}
\end{align}
with nonzero constant amplitudes $a_i>0$ at frequencies $\omega_i \neq \omega_j$. 
%
Such probing frequencies $\omega_i$ can be selected as
\begin{equation}\label{omegadefinition}
\omega_i=\omega_{i}'\omega=\mathcal{O}(\omega)\,, \quad
i\in{1,2,\ldots, N}\,,
\end{equation}
where $\omega$ is a positive constant and $\omega_{i}'$ is a
rational number. One possible choice is given in \cite{GKN:2012} as
\begin{equation}\label{possible_implementation}
\omega_{i}'\not\in \left\{\omega_{j}'\,, \ \frac{1}{2}(\omega_{j}'+\omega_{k}')\,, \ \omega_{j}'+2\omega_{k}'\,, \ \omega_{j}'+\omega_{k}'\pm \omega_{l}'\right\}\,,
\end{equation}
for all distinct $i,j,k$ and $l$. 

By considering $\hat{\theta}_{i}(t)$ as an estimate of $\theta^{*}_{i}$, one can define the \textit{estimation error}:
\begin{align}
\tilde{\theta}_{i}(t)&=\hat{\theta}_{i}(t)-\theta_{i}^{*}\,.\label{eq:tildeThetai}
\end{align}
The estimate $\hat{G}_{i}$ of the unknown gradient of each payoff $J_i$ is given by
\begin{align}
\hat{G}_{i}(t)&=M_{i}(t)y_{i}(t)\,. \label{eq:hatGi_nonaverage}
\end{align}
%
Computing the average of the resulting signal leads us to
\begin{align}
\hat{G}_{i}^{\rm{av}}(t)
&=\sum_{j=1}^{N}\epsilon_{ij}^{i}H_{ij}^{i}\tilde{\theta}_{j}^{\rm{av}}(t-D_{j})\,,\label{eq:hatGi}
\end{align}

with $\Pi$ defined as
\begin{equation} \label{period}
\Pi:= 2 \pi \times \text{LCM}\left\{\frac{1}{\omega_i} \right\}\,, 
\end{equation}
and LCM standing  for the least common multiple.

At this point, if we ignore the prediction loop and the low-pass filter (both indicated in red color) in Fig.~\ref{fig:blockDiagram_v3}, the control law $U_i(t)=k_i\hat{G}_i(t)$ could be obtained as in the classical ES approach. In this case, from equations (\ref{eq:tildeThetai}) and (\ref{eq:hatGi}), we could write the average version of
\begin{equation} \label{cudoerro}
\dot{\tilde{\theta}}_i(t)=U_i(t)
\end{equation}
as 
\begin{align}
\dot{\tilde{\theta}}_{i}^{\rm{av}}(t)&=k_{i}\hat{G}_{i}^{\rm{av}}(t) \nonumber \\
&=k_{i}\sum_{j=1}^{N}\epsilon_{ij}^{i}H_{ij}^{i}\tilde{\theta}_{j}^{\rm{av}}(t-D_{j})\,.\label{dotTildeAvi_v0}
\end{align}
Therefore, by defining $$\tilde{\theta}^{\rm{av}}(t): = [\tilde{\theta}_{1}^{\rm{av}}(t)\,,\tilde{\theta}_{2}^{\rm{av}}(t)\,,\ldots\,,\tilde{\theta}_{N}^{\rm{av}}(t)]^T \in \mathbb{R}^{N}$$ and taking into account all players, one has
\begin{align}
\dot{\tilde{\theta}}^{\rm{av}}(t)&=KH\tilde{\theta}^{\rm{av}}(t-D)\,,\label{dotTildeAv_v0}
\end{align}
with $K:=\text{diag}\{k_{1}\,,\ldots\,,k_{N}\}$ and $H$ given by (\ref{eq:Htheta*h}). Equation (\ref{dotTildeAv_v0}) means that, even if $KH$ was a Hurwitz matrix, the equilibrium $\tilde{\theta}_{\rm{e}}^{\rm{av}}=0$ of the closed-loop average system would not necessarily be stable for arbitrary values of the time-delays $D_i$. This reinforces the demand of employing the prediction feedback $U_i(t)=k_i\hat{G}_i(t+D_i)$---or even its filtered version---for each player to stabilize collectively the closed-loop system, as illustrated with red color in Fig.~\ref{fig:blockDiagram_v3}.

In such a decentralized scenario, the dither frequencies $\omega_{-i}$, the excitation amplitudes $a_{-i}$, and consequently, the individual control laws $U_{-i}(t)$ are not available to Player $i$. Recalling that the model of the payoffs (\ref{eq:Ji}) and (\ref{ajagambiarra_noncooperative}) are also assumed to be unknown, it becomes impossible to reconstruct individually or estimate completely the Hessian matrix $H$ given in (\ref{eq:Htheta*h}) by using demodulating signals such as in \cite{Tiago}. 

%
%

Following the non-sharing information paradigm, the $i$th-player is only able to estimate the element $H_{ii}^{i}$ of the $H$ matrix (\ref{eq:Htheta*h}) by itself, and this being so for all players. Therefore, only the diagonal of $H$ can be properly recovered in the average sense. In this way, the signal $N_{i}(t)$ is now simply defined by:
\begin{align}
N_{i}(t):= N_{ii}(t)=	\dfrac{16}{a^2_i} \bigg( \sin^2 (\omega_i t) - \dfrac{1}{2} \bigg)\,, \label{eq:Ni}
\end{align} 
according to \cite{Tiago}. Then, the average version of 
\begin{align}
\hat{H}_{i}(t)= N_i(t) y_i(t) \label{eq:hatHi_nonaverage}
\end{align}
is given by
\begin{align}
\hat{H}_{i}^{\rm{av}}(t)=\left[N_{i}(t)y_{i}(t)\right]_{\rm{av}}=H_{ii}^{i}\,. \label{eq:hatHi}
\end{align}

In order to compensate for the time delays, the following predictor-based update law was proposed in \cite{JOTA:2021}:
\begin{equation} \label{adaptation_NC}
	\dot{\hat{\theta}}_{i}(t) = U_{i}(t)\,,
\end{equation}
\begin{equation} \label{4_dU_NC}
	\dot{U}_{i}(t) = -c_{i}U_{i}(t) + c_{i}k_{i} \bigg(\hat{G}_{i}(t) +\hat{H}_{i}(t) \int_{t-D_i}^{t} U_i(\tau)d\tau \bigg)\,,
\end{equation}
\noindent for positive constants $k_{i}$ and $c_{i}$.

\section{Noncooperative Games with Heat PDEs} \label{noncooperative_games_CP14}

Recapitulating, the optimality of the respective outcomes (outputs) of Players $i$ and $j$, respectively $y_{i}(t) \in \mathbb{R}$ and $y_{j}(t) \in \mathbb{R}$, do not depend exclusively on their own strategies (input signals) $\Theta_{i}(t)\in \mathbb{R}$ and $\Theta_{j}(t)\in \mathbb{R}$. Moreover, defining $\Theta := [\Theta_1\,, \ldots\,, \Theta_N]^T$, each player's payoff function $J_i$ depends also on 
$\Theta_j$ of the other player $j$, $j\not=i$. An $N$-tuple of $\Theta^*=[\Theta^{*}_{1}\,, \ldots\,,\Theta^{*}_{N}]^T $ is said to be in Nash equilibrium, if no player $i$ can improve its payoff by unilaterally deviating from $\Theta_i^*$, this being so for
all $i\in\{1\,, \ldots\,, N\}$ \cite{Basar:1999}. 

\begin{figure}[ht!]
\begin{center}
\includegraphics[width=3.0 in]{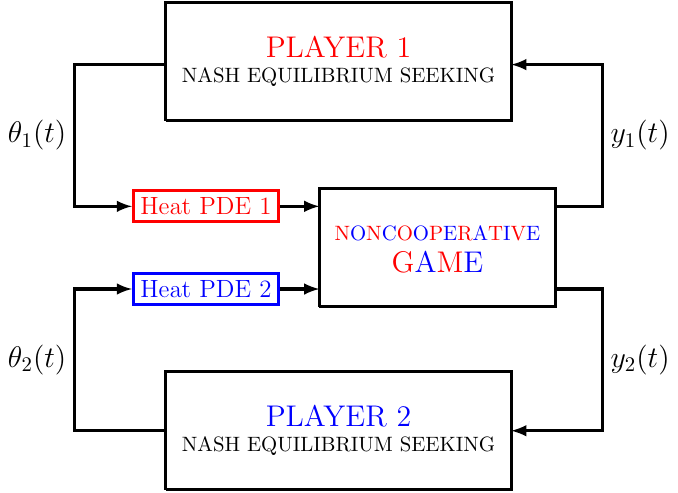}
\end{center}
\caption{Nash equilibrium seeking in a duopoly ($N=2$) noncooperative game with players acting through heat PDE dynamics.}
\label{fig:blockDiagram_v2_CP14}
\end{figure}

\begin{figure}[ht]
\centering
\includegraphics[width=3.25 in]{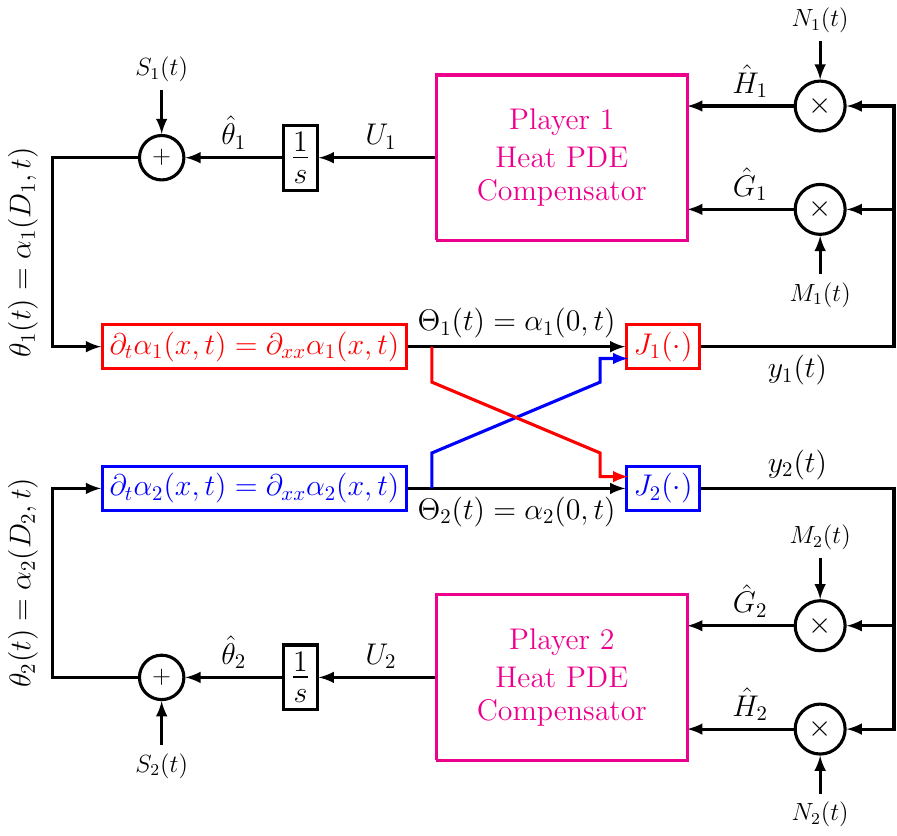}
\caption{Block diagram illustrating the Nash Equilibrium seeking strategy performed for each player ($N=2$). \textcolor{black}{In magenta color are the boundary controllers used to compensate the individual heat PDEs for the noncooperative game.}}
\label{fig:blockDiagram_v3_CP14}
\end{figure}

As shown in Fig. \ref{fig:blockDiagram_v2_CP14} and Fig.  \ref{fig:blockDiagram_v3_CP14}, distinct heat equations (with Dirichlet actuation) are assumed in the vector of player actions
$\theta(t) \in \mathbb{R}^2$, in the particular duopoly game with $N=2$. 

In the $N$-player game, the propagated actuator vector $\Theta(t) \in \mathbb{R}^N$ is given by 
\begin{align}
\label{eq:heat_eqn_start}
&\Theta_i(t)=\alpha_i(0,t)\,, \quad \forall i \in \{ 1, \ldots, N \}\,,  \\
\label{eq:heat_eqn_1}
&\partial_t\alpha_i(x,t)=\partial_{xx}\alpha_i(x,t),\quad x\in(0,D_i)\,, \\
&\partial_x\alpha_i(0,t) = 0\,, \\
\label{eq:heat_eqn_end}
&\alpha_i(D_i,t)= \theta_i(t)\,,
\end{align}
where $\alpha_i: [0\,, D_i] \times \mathbb{R}_+ \to \mathbb{R}$ and each domain length $D_i$ is known. The solution of (\ref{eq:heat_eqn_start})--(\ref{eq:heat_eqn_end}) is given by
\begin{equation} 
\alpha_i(x,t) = \mathcal{L}^{-1}\left[\frac{\cosh(x\sqrt{s})}{\cosh(D_i\sqrt{s})}\right]*\,\theta_i(t),
\label{eq:sol_heat_eqn}
\end{equation}
where $\mathcal{L}^{-1}(\,\cdot\,)$ denotes the inverse Laplace transformation and $*$ is the convolution operator. Given this relation, we define the \textit{diffusion operator} for the PDE (\ref{eq:heat_eqn_start})--(\ref{eq:heat_eqn_end}) with boundary input and measurement given by
$\mathcal{D}\!=\!\text{diag}\{\mathcal{D}_1, \ldots, \mathcal{D}_N\}$ with  
\begin{equation}  \label{bucetasso}
\mathcal{D}_i[\varphi(t)] \!\!=\!\! \mathcal{L}^{-1}\!\left[\frac{1}{\cosh(D_i\sqrt{s})}\right]*\,\varphi(t),\quad 
		\text{s.t.}\quad \Theta(t)\!\!=\!\! \mathcal{D}\left[\theta(t)\right]. 
\end{equation}

As in (\ref{eq:Ji}), we consider games where the payoff function $y_i(t) \!=\! J_i(\mathcal{D}[\theta(t)]) \!=\! J_i(\Theta(t))$ of each player is quadratic 
\cite{Basar:1999}, expressed as 
a strictly concave combination of their actions propagated through distinct heat PDE dynamics 
\begin{align}
J_{i}(\Theta(t))=&\frac{1}{2}\sum_{j=1}^{N}\sum_{k=1}^{N}\epsilon_{jk}^{i}H_{jk}^{i}\Theta_{j}(t)\Theta_{k}(t) 
+\sum_{j=1}^{N}h_{j}^{i}\Theta_{j}(t)+c_{i}\,, \label{eq:Ji_heat}
\end{align} 
where $J_{i}(\Theta) : \mathbb{R}^{N} \!\to\! \mathbb{R}$, obtained by replacing $\theta$ for $\Theta$ in (\ref{eq:Ji}). Equations (\ref{Nashcu})--(\ref{eq:NE_v2}) can be repeated here simply replacing $\theta^*\in \mathbb{R}^{N}$ for $\Theta^* \in \mathbb{R}^{N}$, where $\Theta^*$ represents the Nash equilibrium defined in (\ref{Nashcu}). 
The \textit{objective} is to design an extremum seeking-based strategy to reach the Nash Equilibrium in noncooperative games subjected to heat PDEs in the decision variables of the players (input signals). 


Since our goal is to find the unknown optimal inputs $\Theta^*$ (and $\theta^*$), we define the estimation errors 
\begin{align}
\tilde{\theta}(t) = \hat{\theta}(t)-\theta^*\,, \qquad \vartheta(t)\!=\!\hat{\Theta}(t)-\Theta^* \,,
\label{eq:def_theta_tilde}
\end{align}
where the vectors $\hat{\theta}(t)$ and $\hat{\Theta}(t)$ are the estimates of $\theta^*$ and $\Theta^*$. In order to make \eqref{eq:def_theta_tilde} coherent with the optimizer of the static map $\Theta^*$, we apply the diffusion operator (\ref{bucetasso}) to $\tilde{\theta}_i$ in (\ref{eq:def_theta_tilde}) and get 
%
%
\begin{align}
\label{eq:def_vartheta_start}
&\vartheta_i(t):=\bar{\alpha}_i(0,t), \quad \forall i \in \{ 1, \ldots, N \}\,,  \\ 
&\partial_t\bar{\alpha}_i(x,t)=\partial_{xx}\bar{\alpha}_i(x,t), \quad x\in(0,D_i)\,, \label{eq:sys_diff_theta_2}  \\
&\partial_x\bar{\alpha}_i(0,t)=0\,, \label{eq:sys_diff_theta_3}\\
&\bar{\alpha}_i(D_i,t)=\tilde{\theta}_i(t)\,,
\label{eq:def_vartheta_end}
\end{align}
where $\bar{\alpha}_i: [0\,, D_i] \times \mathbb{R}_+ \to \mathbb{R}$ and $\vartheta(t)\!:=\!\mathcal{D}[\tilde{\theta}(t)]\!=\!\hat{\Theta}(t)-\Theta^*$ is the propagated estimation error $\tilde{\theta}(t)$ through the diffusion domain. For
$\lim\limits_{t\rightarrow \infty} \theta(t)=\theta_{c}$, we have $\lim\limits_{t\rightarrow \infty} \Theta(t) = \Theta_c = \theta_c$, where the index $c$ indicates a constant signal. 
Indeed, from (\ref{eq:sol_heat_eqn}), for a constant input $\theta=\theta_c$, one has $\mathcal{L}\lbrace\theta_c\rbrace = \theta_c/s$, and applying the Laplace limit theorem we get $\lim\limits_{t\rightarrow \infty} \mathcal{D}_i[\theta_{ci}] = \lim\limits_{s\rightarrow 0} \left\lbrace\frac{\theta_{ci}}{\cosh(\sqrt{s}D_i)} \right\rbrace = \theta_{ci}$. In the particular case $\theta=\theta_c=\theta^*$, one has 
%
\begin{equation} \label{meucucheirosa}
\Theta^*=\theta^*\,.
\end{equation}

Fig.~\ref{fig:blockDiagram_v3_CP14} depicts a schematic diagram that summarizes the proposed Nash Equilibrium policy for each player, where their outputs are given by
\begin{align}
y_{i}(t)&=J_{i}(\Theta(t))\,.
\label{eq:yi_v0_CP14}
\end{align}
%
The additive dither signals in the presence of heat PDE dynamics \cite{FKKO:2018} are defined according to 
\begin{align}
S_{i}(t)&=\frac{1}{2}a_{i}e^{\sqrt{\frac{\omega_{i}}{2}}D_{i}}\sin\left(\omega_{i}t+\sqrt{\frac{\omega_{i}}{2}}D_{i}\right) \nonumber \\
&\quad+\frac{1}{2}a_{i}e^{-\sqrt{\frac{\omega_{i}}{2}}D_{i}}\sin\left(\omega_{i}t-\sqrt{\frac{\omega_{i}}{2}}D_{i}\right)
\,,\label{eq:Si_CP14}
\end{align}
and the multiplicative demodulation signals are given by
\begin{align}
M_{i}(t)&=\frac{2}{a_{i}}\sin(\omega_{i}t)
\,, \label{eq:Mi_CP14}
\end{align}
%
%
with nonzero constant amplitudes $a_{i}>0$ at frequencies $\omega_i \neq \omega_j$. 
%
Such probing frequencies $\omega_i$ can be selected as

\begin{equation}\label{omegadefinition_CP14}
\omega_i=\omega_{i}'\omega=\mathcal{O}(\omega)\,, \quad
\forall i \in \{ 1, \ldots, N \}\,,
\end{equation}
where $\omega$ is a positive constant and $\omega_{i}'$ is a
rational number---one possible choice is given in \cite{GKN:2012}. 

Following the non-sharing information paradigm, only the diagonal elements of $H$ can be properly recovered in the average sense by players. 
In this sense, the signals $N_{i}(t)$ 
are simply defined by \cite{GKN:2012}:
\begin{align}
N_{i}(t)=	\frac{16}{a^2_i} ( \sin^2 (\omega_i t) - \frac{1}{2}) \,.
\label{eq:Ni_CP14}
\end{align} 
%
Then, the average version of 
\begin{align}
\hat{H}_{i}(t)= N_i(t) y_i(t) 
\label{eq:hatHi_nonaverage_CP14}
\end{align}
%
is given by
\begin{align}
\hat{H}_{i}^{\rm{av}}(t)=\left[N_{i}(t)y_{i}(t)\right]_{\rm{av}}=H_{ii}^{i}
\,. \label{eq:hatHi_CP14}
\end{align} 
%

Considering $\hat{\theta}_{i}(t)$ 
as the estimates of $\theta^{*}_{i}$, 
one can define
from (\ref{eq:def_theta_tilde}) the individual estimation errors: 
\begin{align}
\tilde{\theta}_{i}(t)=\hat{\theta}_{i}(t)-\theta_{i}^{*}\,, \quad \vartheta_{i}(t)=\hat{\Theta}_{i}(t)-\Theta_{i}^{*}
\,.\label{eq:tildeThetai_CP14}
\end{align}
The estimate of the unknown gradients of the payoff functions are given by
\begin{align}
\hat{G}_{i}(t)&=M_{i}(t)y_{i}(t)
\,, \label{eq:hatGi_nonaverage_CP14}
\end{align}
%
and computing the average of the resulting signal, leads us to
\begin{align}
\hat{G}^{\rm{av}}(t)=H\vartheta^{\rm{av}}(t)
\,,\label{eq:hatGi_CP14}
\end{align}
where the Hessian $H$ is given in (\ref{eq:Htheta*h}). 
%

Additionally, from the block diagram in Fig. \ref{fig:blockDiagram_v3_CP14} for arbitrary $N$, one has
\begin{align}
\dot{\hat{\theta}}_i(t) = U_i(t)\,, \quad \forall i \in \{ 1, \ldots, N \}  \,,
\label{eq:theta_hat_U}
\end{align}
and, consequently, 
\begin{align}
\dot{\tilde{\theta}}_i(t)=U_i(t)\,, \quad \forall i \in \{ 1,\ldots, N \} \,,
\label{eq:theta_hat_Ui}
\end{align}
since $\dot{\tilde{\theta}}(t)=\dot{\hat{\theta}}(t)$, once $\theta^*$ is constant. 
%
%
Taking the time derivative of \eqref{eq:def_vartheta_start}--\eqref{eq:def_vartheta_end} and with the help of \eqref{eq:def_theta_tilde} and \eqref{eq:theta_hat_Ui}, the \textit{propagated error dynamics} are written as
\begin{align}
\label{eq:error_dyn_start}
&\dot{\vartheta}_i(t)=u_i(0,t)\,,  \quad \forall i \in \{ 1, \ldots, N \}\,,  \\ 
&\partial_t u_i(x,t)=\partial_{xx}u_i(x,t)\,, \quad x\in(0,D_i)\,, \\
&\partial_x u_i(0,t)=0\,, \\
&u_i(D_i,t)=U_i(t)\,,
\label{eq:error_dyn_end}
\end{align}
where $u_i: [0\,, D_i] \times \mathbb{R}_+ \to \mathbb{R}$ and  $u_i(x,t):=\partial_t\bar{\alpha}_i(x,t)$. 



Hence, from (\ref{eq:hatGi_CP14}) and (\ref{eq:error_dyn_start})--(\ref{eq:error_dyn_end}),  
it is possible to find a compact form for the overall average estimated gradient  according to 

\begin{align}
\dot{\hat{G}}^{\rm{av}}(t) &= H \dot{\vartheta}^{\rm{av}}(t) = H \mathcal{D}[U^{\rm{av}}(t)]\,, \label{eq:dotHatGAv_CP14} \\ \nonumber 
\end{align}
\parbox{\linewidth}{where $\vartheta^{\rm{av}}(t)\!:=\! [\vartheta_{1}^{\rm{av}}(t), \dots, \vartheta_{N}^{\rm{av}}(t)]^T \in \mathbb{R}^N$, 
$\hat{G}^{\rm{av}}(t)\!:=\! [\hat{G}_{1}^{\rm{av}}(t), \dots, \hat{G}_{N}^{\rm{av}}(t)]^T \in \mathbb{R}^N$ and 
$U^{\rm{av}}(t)\!:=\! [U_{1}^{\rm{av}}(t), \dots, U_{N}^{\rm{av}}(t)]^T \in \mathbb{R}^N$ are the average versions of 
$U(t)\!:=\! [U_{1}(t), \dots, U_{N}(t)]^T, \vartheta(t)\!:=\! [\vartheta_{1}(t), \dots, \vartheta_{N}(t)]^T$ \and  $\hat{G}(t):= [\hat{G}_{1}(t), \dots, \hat{G}_{N}(t)]^T$.}

Now, we use an extremum seeking strategy based on boundary control to compensate for the diffusion operator
$\mathcal{D}[\cdot]$ in (\ref{eq:dotHatGAv_CP14}) due to the multiple and distinct heat PDEs in the players' actions. Basically, the control laws must be able to ensure exponential stabilization of $\hat{G}^{\rm{av}}(t)$ and, consequently, of $\vartheta^{\rm{av}}(t)\!=\!\hat{\Theta}^{\rm{av}}(t)\!-\!\Theta^*$. From (\ref{eq:hatGi_CP14}), it is clear that, if $H$ is invertible, $\vartheta^{\rm{av}}(t) \!\to\! 0$ as $\hat{G}^{\rm{av}}(t) \!\to\! 0$. Hence, the convergence of $\vartheta^{\rm{av}}(t)$ to the origin results in the convergence of $\Theta(t)$
to a small neighborhood of the Nash Equilibrium point $\Theta^*=\theta^*$---see (\ref{meucucheirosa})---via averaging theory \cite{HL:1990}, without allowing sharing of any information among the players.


\subsection{Decentralized PDE Boundary Control using only the Known Diagonal Terms of the Hessian}

In this sense, we are able to formulate the closed-loop system in a \textit{decentralized} fashion, where no knowledge about the payoff or action of any other player is required. 

%

Inspired by \cite{ORKB:2020}, where ES was considered for PDE compensation but not in the context of games, we propose the following boundary-based update laws
$\dot{\hat{\theta}}_i(t)=U_i(t)$, $i\in\{1\,,\ldots\,, N\}$:
%
%
\begin{align} \label{4_dU_NC_CP14}
	\dot{U}_{i}(t) &= -c_{i}U_{i}(t) \nonumber \\
 &\quad+ c_{i}k_{i} \left(\hat{G}_{i}(t) +\hat{H}_{i}(t) \int_{0}^{D_i}(D_i-\tau) u_i(\tau,t)d\tau \right)
\,,
\end{align}
%
for positive constants $k_{i}>0$ and $c_{i}>0$, in order to compensate for the heat PDEs in (\ref{eq:error_dyn_start})--(\ref{eq:error_dyn_end}).  Again, with some abuse of notation, constants $c_i$ were chosen to denote the parameters of the control laws, but they are not related to those which appear in the payoff functions given in (\ref{eq:Ji_heat}). 
%
%

As discussed in \cite[Remark~2]{FKKO:2018}, 
the boundary control law (\ref{4_dU_NC_CP14}) could be rewritten as
\begin{align}\label{predictor_new}
\dot{U}_i(t)&=-c_i U_i(t) \nonumber \\
&\quad+c_ik_i \left[\hat{G}_i(t)\!+\!\hat{H}_i(t)\left(\hat{\theta}_i(t)\!-\!\Theta_i(t)\!+\!a_i \sin (\omega_i t)\right)\right], 
\end{align}
using the diffusion equations $\partial_t\alpha_i(x,t)\!=\!\partial_{xx}\alpha_i(x,t)$, $\forall i\!\in\!\{1\,, \ldots\,, N\}$, and the integration by parts, associated with
(\ref{eq:heat_eqn_start})--(\ref{eq:heat_eqn_end}), (\ref{eq:def_theta_tilde}) and recalling that $\vartheta_i+a_i\sin(\omega_i t)=
\Theta_i(t)-\Theta_i^*$, analogously to \cite[Eq. (25)]{FKKO:2018}.

\subsection{ISS-Like Properties for Parabolic PDE Representation}

For the sake of simplicity, let us assume that $c_{i} \!\to\! +\infty$ in (\ref{4_dU_NC_CP14}), resulting in the following general expression:
\begin{equation} \label{rabao_NC_CP14}
U_{i}(t)=k_{i}\left(\hat{G}_{i}(t) +\hat{H}_{i}(t) \int_{0}^{D_i}(D_i-\tau) u_i(\tau,t)d\tau\right)\,.
\end{equation}
Recalling (\ref{eq:error_dyn_start})--(\ref{eq:error_dyn_end}), the infinite-dimensional closed-loop system (\ref{eq:dotHatGAv_CP14}) and (\ref{rabao_NC_CP14}) in its average version can be written in the corresponding PDE
representation form, as 
\begin{eqnarray}
\dot{\hat{G}}^{\rm{av}}(t) &=& H u^{\rm{av}}(0,t)\,, \label{saco1_NC_CP14}\\
\partial_t u^{\rm{av}}(x,t)&=&D^{-2}\partial_{xx} u^{\rm{av}}(x,t)\,, \quad x\in (0\,,1)\,, \label{saco2_NC_CP14}\\
\partial_x u^{\rm{av}}(0,t)&=&0\,, \\
u^{\rm{av}}(1\,,t)&=&U^{\rm{av}}(t) \label{saco3_NC_CP14}\,,
\end{eqnarray}

with $D\!=\!\text{diag}\{D_1, \ldots, D_N\}$. 

In the \textit{reduction-like approach} \cite{a14} (or finite-spectrum assignment), we use the transformation (for $i \in \{1\,, \ldots\,, N\}$):
\begin{eqnarray} 
\bar{G}^{\rm{av}}_i(t)\!\!\!\!\!&=&\!\!\!\!\!\hat{G}^{\rm{av}}_i(t)\!+\!\sum_{j=1}^{N} \epsilon_{ij}^{i} H_{ij}^{i} \!\int_{0}^{D_j}\!\!\!(D_j\!-\!\tau) u^{\rm{av}}_j(\tau,t) d\tau
\nonumber \\ 
\!\!\!\!\!&=&\!\!\!\!\! \hat{G}^{\rm{av}}_i(t) \!+\! \sum_{j=1}^{N} \epsilon_{ij}^{i} H_{ij}^{i} \!\int_{0}^{1}\!\!\!D_j^2(1\!-\!\xi) u^{\rm{av}}_j(\xi,t) d\xi\,,~~~~
\label{transformacao_cu_NC_CP14}
\end{eqnarray}
where 
$\int_{0}^{D_j}(D_j-\tau) u^{\rm{av}}_j(\tau,t) d\tau=\int_{0}^{1} D_j^2(1-\xi) u^{\rm{av}}_j(\xi,t) d\xi$. 
%
%
With some mathematical manipulations, it is not difficult to see that $\bar{G}^{\rm{av}}$ satisfies
\begin{eqnarray}  \label{Zreduction_NC_CP14}
\dot{\bar{G}}^{\rm{av}}(t)=H U^{\rm{av}}(t)\,.
\end{eqnarray}
Now, after adding and subtracting
the next terms in blue and red in (\ref{rabao_NC_CP14}), it is rewritten as:
%
\begin{align} \label{rabao_NC_1_CP14}
&U_{i}(t)=k_{i}\left(\hat{G}_{i}(t) +\hat{H}_{i}(t) \int_{0}^{D_i} (D_i-\tau) u_i(\tau,t)d\tau\right. \nonumber \\
&\quad\left.+\textcolor{blue}{\sum_{j\neq i}\epsilon_{ij}^{i}H_{ij}^{i}\int_{0}^{D_j} (D_j-\tau) u_j(\tau,t) d\tau}\right) \nonumber\\
&~~~-k_{i}\textcolor{red}{\sum_{j\neq i}\epsilon_{ij}^{i}H_{ij}^{i} \int_{0}^{D_j} (D_j-\tau) u_j(\tau,t) d\tau}\,,
\end{align}
%
%
whose average compact form is
\begin{align} \label{rabao_NC_1_cu}
U^{\rm{av}}(t)&=K\bar{G}^{\rm{av}}(t)+\epsilon K\phi^{\rm{av}}(D,t), 
\end{align}
where the matrix $K:=\text{diag}\{k_{1},\ldots\,,k_{N}\}$ with entries $k_i>0$
and the auxiliary variable $\phi(D,t):=[\phi_{1}(D,t),\ldots\,,\phi_{N}(D,t)]^{T}\in\mathbb{R}^{N}$ is defined as

\begin{align}
\phi_{i}(D,t)&:=-\sum_{j\neq i}H_{ij}^{i}\int_{0}^{D_j} (D_j-\tau) u_j(\tau) d\tau\,, \nonumber \\
\phi_{i}(1,t)&:=-\sum_{j\neq i}H_{ij}^{i}\int_{0}^{1} D_j^2 (1-\xi) u_j(\xi,t) d\xi. \label{eq:phii_CP14} 
\end{align}

%
%
%
%
Then, it is possible to find a compact form for the overall average game from equations (\ref{Zreduction_NC_CP14}) and
(\ref{rabao_NC_1_cu}), such as
\begin{eqnarray}
\dot{\bar{G}}^{\rm{av}}(t) &=& HK \bar{G}^{\rm{av}}(t)+\epsilon HK \phi^{\rm{av}}(1,t)\,, \label{saco1M_NC_CP14}\\
\partial_t u^{\rm{av}}(x,t)&=&D^{-2}\partial_{xx} u^{\rm{av}}(x,t)\,, \quad x\in (0\,,1)\,, \label{saco2M_NC_CP14}\\
\partial_x u^{\rm{av}}(0,t)&=&0\,, \\
u^{\rm{av}}(1\,,t)&=&K \bar{G}^{\rm{av}}(t)+\epsilon K \phi^{\rm{av}}(1,t) \label{saco3M_NC_CP14}\,.
\end{eqnarray}
%

From (\ref{saco1M_NC_CP14}), if $HK$ is Hurwitz, it is clear that the dynamics of the ODE state variable $\bar{G}^{\rm{av}}(t)$ is exponentially Input-to-State Stable (ISS) \cite{KK:2018} with respect to the PDE state $u(x,t)$ by means of the function $\phi^{\rm{av}}(1,t)$. Moreover, the PDE subsystem  (\ref{saco2M_NC_CP14}) is ISS \cite{KK:2018} with respect to $\bar{G}^{\rm{av}}(t)$ in the boundary condition $u^{\rm{av}}(1\,,t)$.

\section{Stability Analysis}
\label{doritos_stability}

We will now show that the ODE-PDE loops (e.g.,  
(\ref{saco1M_NC_CP14})--(\ref{saco3M_NC_CP14})), for the corresponding delay and heat case, contain small-parameters $\epsilon$ which can lead to closed-loop stability if they are chosen sufficiently small. To attain $\theta^*$ stably in real time, without any model information (except for the delays or difusion domains $D_i$), each Player $i$ employs the noncooperative extremum seeking strategy (\ref{4_dU_NC}) or (\ref{4_dU_NC_CP14}) via boundary control feedback.


\subsection{Heat PDEs Case}
\label{doritos_stability_case1}

The next theorem provides the stability/convergence properties of the closed-loop extremum seeking feedback for the $N$-player noncooperative game with heat PDEs. 


\begin{theorem} \label{ch16.theorem.16.1.vacina}
Consider the closed-loop system (\ref{saco1_NC_CP14})--(\ref{saco3_NC_CP14}) under \textcolor{black}{Assumptions~\ref{ch.13.Assumption 1.} and \ref{ch13.Assumption 2.}}, and multiple heat PDEs  
(\ref{eq:heat_eqn_start})--(\ref{eq:heat_eqn_end}) with distinct diffusion coefficients $D_i$ for the $N$-player quadratic noncooperative game with payoff functions given in (\ref{eq:Ji_heat}) and (\ref{eq:yi_v0_CP14}) and control laws $U_i(t)$ defined in (\ref{4_dU_NC_CP14}). There exist $c_i>0$ and $\omega>0$ sufficiently large as well as $\epsilon>0$ sufficiently small such that the closed-loop system  
with state $\vartheta_i(t)$, $u_i(x,t)$, $\forall i\in \{1\,,\ldots\,,N\}$, has a unique locally exponentially stable periodic solution in $t$ of period $\Pi$ in (\ref{period}), with $\omega_i$ in (\ref{omegadefinition_CP14}) of order $\mathcal{O}(\omega)$ according to (\ref{omegadefinition}), denoted by $\vartheta_i^{\Pi}(t)$, $u_i^{\Pi}(x,t)$, which satisfies, $\forall t\geq 0$:
%
\begin{eqnarray}\label{periodic_solution_NCG_CP14}
\left(\sum_{i=1}^{2}\left[\vartheta_{i}^{\Pi}(t)\right]^2  \!+\! \int_{0}^{D_i} \left[u_{i}^{\Pi}(x,t)\right]^2dx \right)^{1/2} \!\leq\! \mathcal{O}(1/\omega)\,.
\end{eqnarray}
%
In particular,

\begin{eqnarray}\label{limsup1_NCG_CP14}
\limsup_{t\to+\infty}|\Theta(t)\!-\!\Theta^*|\!\!\!\!\!\!&=&\!\!\!\!\!\!\mathcal{O}(|a|\!+\!1/\omega), 
\\
\limsup_{t\to+\infty}|\theta(t)\!-\!\theta^*|\!\!\!\!\!\!&=&\!\!\!\!\!\!\mathcal{O}\!\left(|a|e^{\max(D_i)\sqrt{\omega/2}}\!+\!1/\omega\right), \label{limsup2_CP14}
\end{eqnarray}
where $a=[a_1 \ a_2]^T$ and $\theta^*\!=\!\Theta^*$ is the unique (unknown) \textcolor{black}{Nash Equilibrium given by (\ref{eq:NE_v2}).} 
\end{theorem}

\begin{proof}
See the Appendix.
\end{proof}

\subsection{Delay Case}
\label{doritos_stability_case2}

On the other hand, the derivative of (\ref{eq:hatGi}) is
\begin{align}
\dot{\hat{G}}_{i}^{\rm{av}}
(t)&=\sum_{j=1}^{N}\epsilon_{ij}^{i}H_{ij}^{i}\dot{\tilde{\theta}}_{j}^{\rm{av}}(t-D_{j})\,. \label{eq:dotHatGi_v1}
\end{align} 
By delaying by $D_i$ units the time-argument of both sides of the average version of 
(\ref{cudoerro}), we obtain 
\begin{align}
\dot{\tilde{\theta}}_{i}^{\rm{av}}(t-D_{i})=U_{i}^{\rm{av}}(t-D_{i})\,. \label{dotTildeAvi}
\end{align}
Thus, equation (\ref{eq:dotHatGi_v1}) can be rewritten as 
\begin{align}
\dot{\hat{G}}_{i}^{\rm{av}}(t)&=\sum_{j=1}^{N}\epsilon_{ij}^{i}H_{ij}^{i}U_{j}^{\rm{av}}(t-D_{j})\,. \label{eq:dotHatGi_v2}
\end{align} 
Taking into account all players, from  (\ref{eq:hatGi}) and (\ref{eq:dotHatGi_v2}), it is possible to find a compact form for the overall average estimated gradient $\hat{G}^{\rm{av}}
(t)\!:=\! [\hat{G}_{1}^{\rm{av}}(t)\,,\ldots\,,\hat{G}_{N}^{\rm{av}}
(t)]^T\!\in\!\mathbb{R}^{N}$ according to 
\begin{align}
\hat{G}^{\rm{av}}(t)&=H \tilde{\theta}^{\rm{av}}(t-D)\,, \label{eq:hatGAv}\\
\dot{\hat{G}}^{\rm{av}}(t)&=H U^{\rm{av}}(t-D)\,, \label{eq:dotHatGAv}
\end{align}
where $H$ is given in (\ref{eq:Htheta*h}) and $$U^{\rm{av}}(t)\!:=\! [U_{1}^{\rm{av}}
(t)\,,U_{2}^{\rm{av}}(t)\,,\ldots\,,U_{N}^{\rm{av}}(t)]^T \!\in\!\mathbb{R}^{N}.$$

The next theorem provides the stability/convergence properties of the closed-loop extremum seeking feedback for the $N$-player noncooperative game with delays and non-sharing of information. 



\begin{theorem} \label{scheme2} Consider the closed-loop system 
(\ref{dotTildeAvi}), or (\ref{eq:dotHatGi_v2}), 
under Assumptions~\ref{ch.13.Assumption 1.} and
\ref{ch13.Assumption 2.}, and multiple and distinct input delays $D_i$ for the $N$-player quadratic noncooperative game with no information sharing, with payoff functions given in \textcolor{black}{\textnormal{(\ref{ajagambiarra_noncooperative})}} and control laws $U_i(t)$ defined in \textcolor{black}{\textnormal{(\ref{4_dU_NC})}}. There exist $c_i>0$ and $\omega>0$ sufficiently large as well as $\epsilon>0$ sufficiently small such that closed-loop system  
with state $\tilde\theta_i(t-D_i)$, $U_i(\tau)$, $\forall \tau \in [t-D_i,t]$ and $\forall i\in{1,2,\ldots, N}$, has a unique locally exponentially stable periodic solution in $t$ of period $\Pi$ in (\ref{period}), denoted by $\tilde{\theta}_i^{\Pi}(t-D_i)$, $U_i^{\Pi}(\tau)$, $\forall \tau \in [t-D_i,t]$ satisfying, $\forall t\geq 0$:
\begin{eqnarray}\label{periodic_solution_NCG}
\left(\sum_{i=1}^{N}\left[\tilde{\theta}_{i}^{\Pi}(t\!-\!D_i)\right]^2  \!+\! \int_{t-D_i}^{t} \!\!\!\left[U_{i}^{\Pi}(\tau)\right]^2d\tau \right)^{1/2} \!\!\!\!\leq\! \mathcal{O}(1/\omega)\,.
\end{eqnarray}
In particular,
\begin{eqnarray}\label{limsup1_NCG}
\limsup_{t\to+\infty}|\theta(t)-\theta^*|&=&\mathcal{O}(|a|+1/\omega)\,, 
\end{eqnarray}
where $a=[a_1 \ a_2 \ \cdots \ a_N]^T$ and $\theta^*$ is the unique Nash Equilibrium given by \textcolor{black}{\textnormal{(\ref{eq:NE_v2})}}. 
\end{theorem}

\begin{proof}
See the Appendix.
\end{proof}

\section{Noncooperative Duopoly with Heterogeneous Transport-Heat PDE Dynamics} \label{ch15.noncooperative_games}

Here, we propose a non-model based strategy for locally stable convergence to Nash equilibrium in a quadratic noncooperative game with player actions $\Theta_i(t)\in \mathbb{R}$ subject to heterogeneous PDE dynamics. For the sake of simplicity, we consider the duopoly scenario ($N=2$), where different players use different types of PDEs, one player compensating for a delay (transport PDE) and the other one a heat (diffusion) PDE, with each player having access only to its own payoff value, $y_{i}(t) \in \mathbb{R}$, for $i\in\{1\,, 2\}$. Heretofore, we solved Nash equilibrium seeking problems with homogeneous games, where the PDE dynamics of distinct nature were not allowed.

\subsection{Quadratic Payoffs and Heterogeneous PDEs}

\begin{figure}[ht]
\begin{center}
\includegraphics[width=3.2 in]{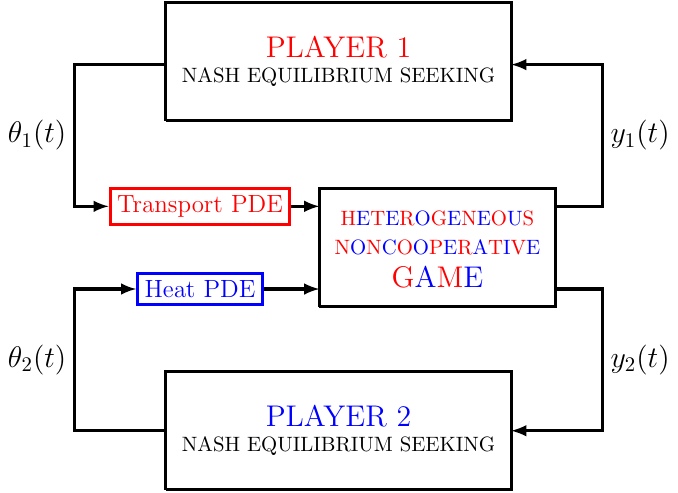}
\end{center}
\caption{Nash equilibrium seeking in a heterogeneous noncooperative game with players acting through transport-heat PDE dynamics.}
\label{ch15.fig:blockDiagram_v2}
\end{figure}

\begin{figure}[ht]
\begin{center}
\includegraphics[width=3.25 in]{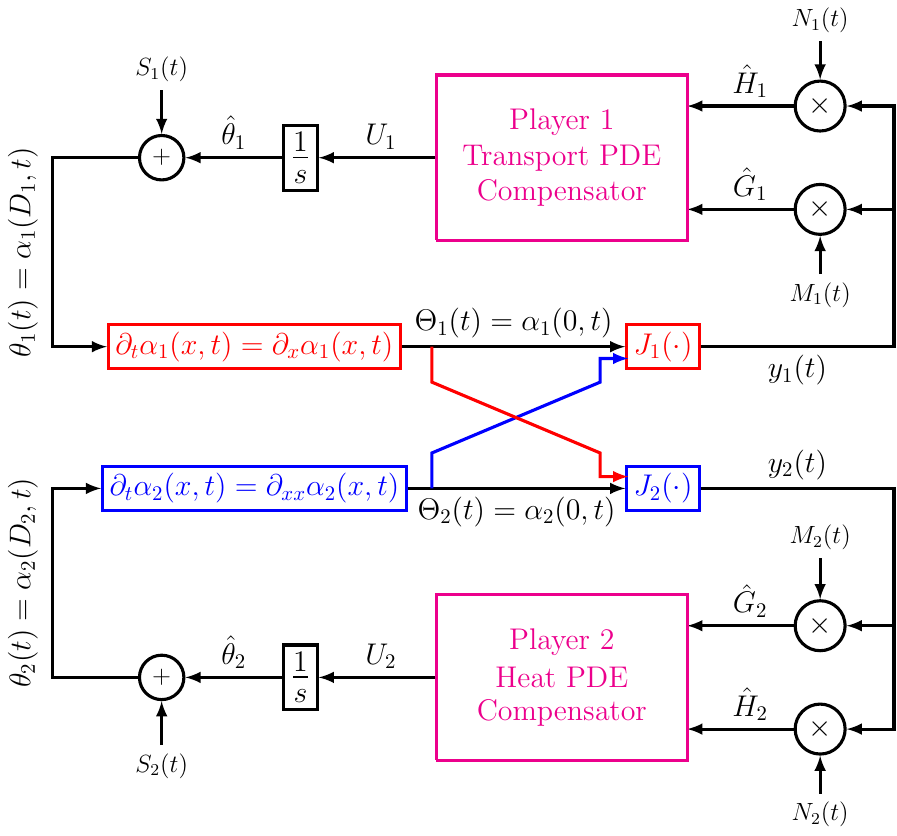}
\end{center}
\vspace{-0.5cm}
\caption{Block diagram illustrating the Nash Equilibrium seeking strategy performed for each player. \textcolor{black}{In magenta color are the boundary controllers used to compensate the individual transport-heat PDEs for the heterogeneous noncooperative game.}}
\label{ch15.fig:blockDiagram_v3}
\end{figure}

As shown in Fig.~\ref{ch15.fig:blockDiagram_v2} and Fig.~\ref{ch15.fig:blockDiagram_v3}, now distinct (transport and heat) PDE equations (with Dirichlet actuation) are assumed in the vector of player actions
$\theta(t) \in \mathbb{R}^2$. Thus, the propagated actuator vector $\Theta(t):= [\Theta_1(t)\,, \Theta_2(t)]^T \in \mathbb{R}^2$ is given by the following transport PDE for player $P_1$ 
\begin{align}
\label{ch15.eq:heat_eqn_start_delay}
&\Theta_1(t)= \theta_1(t-D_1) = \alpha_1(0,t)\,, \\
\label{ch15.eq:heat_eqn_1_delay}
&\partial_t\alpha_1(x,t)=\partial_{x}\alpha_1(x,t),\quad x\in(0,D_1)\,, \\
&\partial_x\alpha_1(0,t) = 0\,, \\
\label{ch15.eq:heat_eqn_end_delay}
&\alpha_1(D_1,t)= \theta_1(t)\,,
\end{align}
and the next heat PDE for player $P_2$ 
\begin{align}
\label{ch15.eq:heat_eqn_start}
&\Theta_2(t)=\alpha_2(0,t)\,,  \\
\label{ch15.eq:heat_eqn_1}
&\partial_t\alpha_2(x,t)=\partial_{xx}\alpha_2(x,t),\quad x\in(0,D_2)\,, \\
&\partial_x\alpha_2(0,t) = 0\,, \\
\label{ch15.eq:heat_eqn_end}
&\alpha_2(D_2,t)= \theta_2(t)\,,
\end{align}
where $\alpha_i: [0\,, D_i] \times \mathbb{R}_+ \to \mathbb{R}$, $\forall i \in \{ 1, 2 \}$, and each domain length $D_i$ is known. 

The solution of (\ref{ch15.eq:heat_eqn_start_delay})--(\ref{ch15.eq:heat_eqn_end_delay}) is 
\begin{equation} 
\alpha_1(x,t) = \theta_1(t+x-D_1)\,,
\label{ch15.eq:sol_heat_eqn_delay}
\end{equation}
which represents a delayed action for player $P_1$ at $x=0$.
 
On the other hand, the solution of (\ref{ch15.eq:heat_eqn_start})--(\ref{ch15.eq:heat_eqn_end}) is given by
\begin{equation} 
\alpha_2(x,t) = \mathcal{L}^{-1}\left[\frac{\cosh(x\sqrt{s})}{\cosh(D_2\sqrt{s})}\right]*\,\theta_2(t),
\label{ch15.eq:sol_heat_eqn}
\end{equation}
where $\mathcal{L}^{-1}(\,\cdot\,)$ denotes the inverse Laplace transformation and $*$ is the convolution operator. Given these relations, we define the \textit{heterogeneous} \textit{transport}-\textit{diffusion operator} $\mathcal{H}=\text{diag}\{\mathcal{H}_1\,,\mathcal{H}_2\}$ for the PDEs (\ref{ch15.eq:heat_eqn_start_delay})--(\ref{ch15.eq:heat_eqn_end_delay}) and (\ref{ch15.eq:heat_eqn_start})--(\ref{ch15.eq:heat_eqn_end}) with boundary inputs and measurements given by
\begin{small}
\begin{align}  \label{ch15.bucetasso}
\mathcal{H}_1[\varphi(t)] \!=\! \varphi(t+x-D_1),  
		~~\text{s.t.}~~ \Theta_1(t)\!\!=\!\! \mathcal{H}_1\left[\theta_1(t)\right], \nonumber\\
\mathcal{H}_2[\varphi(t)] \!\!=\!\! \mathcal{L}^{-1}\!\left[\frac{1}{\cosh(D_2\sqrt{s})}\right]*\,\varphi(t),  
		~~\text{s.t.}~~ \Theta_2(t)\!\!=\!\! \mathcal{H}_2\left[\theta_2(t)\right]. 
\end{align}
\end{small}

We consider games where the payoff function $y_i(t) \!=\! J_i(\Theta(t))$, $\forall i \in \{ 1, 2 \}$, of each player is quadratic 
\cite{Basar:1999}, expressed as 
a strictly concave combination of their actions propagated through distinct transport-heat PDE dynamics 
\begin{align}
J_{1}(\Theta(t))&=\frac{H_{11}^{1}}{2}\Theta_{1}^{2}(t)+\frac{H_{22}^{1}}{2}\Theta_{2}^{2}(t)+ \epsilon H_{12}^{1}\Theta_{1}(t)\Theta_{2}(t) + \nonumber \\
&\quad+h_{1}^{1}\Theta_{1}(t)+h_{2}^{1}\Theta_{2}(t)+c_{1}\,, \label{ch15.eq:J1}  \\
J_{2}(\Theta(t))&=\frac{H_{11}^{2}}{2}\Theta_{1}^{2}(t)+\frac{H_{22}^{2}}{2}\Theta_{2}^{2}(t)+ \epsilon H_{21}^{2}\Theta_{1}(t)\Theta_{2}(t) + \nonumber \\
&\quad+h_{1}^{2}\Theta_{1}(t)+h_{2}^{2}\Theta_{2}(t)+c_{2}\,, \label{ch15.eq:J2} 
\end{align}  
where $J_{1}(\Theta)\,, J_{2}(\Theta) : \mathbb{R}^{2} \!\to\! \mathbb{R}$,
$H_{jk}^{i}$, $h_{j}^{i}$, $c_{i} \!\in\! \mathbb{R}$ are constants, $H_{ii}^{i}\!<\!0$,  $\forall i,j,k \!\in\! \{1\,,2\}$, \textcolor{black}{and $\epsilon\!>\!0$ without loss of generality.}



For the sake of completeness, we provide here in mathematical terms, the definition of a Nash equilibrium $\Theta^*=[\Theta^*_1\,, \Theta_2^*]^T$ in a $2$-player game:
%
\begin{equation} \label{ch15.Nashcu}
J_1(\Theta_1^*\,,\Theta_{2}^*) \!\geq\! J_1(\Theta_1\,,\Theta_{2}^*) ~~  \mbox{and} ~~ J_2(\Theta_1^*\,,\Theta_{2}^*) \!\geq\! J_2(\Theta_1^{*}\,,\Theta_{2}). \!\!
\end{equation}
%
Hence, no player has any incentive to unilaterally deviate its action from $\Theta^*$.

As done before to determine (\ref{eq:NE_v2}) in the general case, the Nash equilibrium solution $\Theta^{*}=-H^{-1}h$ for the duopoly is simply:
%
%
\begin{align}
\begin{bmatrix}
\Theta_{1}^{*} \\
\Theta_{2}^{*} 
\end{bmatrix}
=-
\begin{bmatrix}
 H_{11}^{1} & \epsilon H_{12}^{1} \\
\epsilon H_{21}^{2} &  H_{22}^{2}
\end{bmatrix}^{-1}
\begin{bmatrix}
h_{1}^{1} \\
h_{2}^{2}  
\end{bmatrix} \label{ch15.eq:NE_v2} \,.
\end{align}
%

The \textit{objective} is to design an extremum seeking-based strategy to reach the Nash Equilibrium in heterogeneous noncooperative games subjected to transport-heat PDEs in the decision variables of the players (input signals). 

Since our goal is to find the unknown optimal inputs $\Theta^*$ (and $\theta^*$), we define the estimation errors 
\begin{align}
\tilde{\theta}(t) = \hat{\theta}(t)-\theta^*\,, \qquad \vartheta(t)\!=\!\hat{\Theta}(t)-\Theta^* \,,
\label{ch15.eq:def_theta_tilde}
\end{align}
where the vectors $\hat{\theta}(t)$ and $\hat{\Theta}(t)$ are the estimates of $\theta^*$ and $\Theta^*$. In order to make \eqref{ch15.eq:def_theta_tilde} coherent with the optimizer of the static map $\Theta^*$, we apply the heterogeneous transport-diffusion operator (\ref{ch15.bucetasso}) to $\tilde{\theta}_i$ in (\ref{ch15.eq:def_theta_tilde}), arriving at   
%
%
%
\begin{align}
\label{ch15.eq:def_vartheta_start_delay}
&\vartheta_1(t)= \tilde{\theta}_1(t-D_1) :=\bar{\alpha}_1(0,t)\,,   \\ 
&\partial_t\bar{\alpha}_1(x,t)=\partial_{x}\bar{\alpha}_1(x,t)\,, \quad x\in(0,D_1)\,, \label{ch15.eq:sys_diff_theta_2_delay}  \\
&\partial_x\bar{\alpha}_1(0,t)=0\,, \label{ch15.eq:sys_diff_theta_3_delay}\\
&\bar{\alpha}_1(D_1,t)=\tilde{\theta}_1(t)\,,
\label{ch15.eq:def_vartheta_end_delay}
\end{align}
and
\begin{align}
\label{ch15.eq:def_vartheta_start}
&\vartheta_2(t):=\bar{\alpha}_2(0,t)\,,   \\ 
&\partial_t\bar{\alpha}_2(x,t)=\partial_{xx}\bar{\alpha}_2(x,t)\,, \quad x\in(0,D_2)\,, \label{ch15.eq:sys_diff_theta_2}  \\
&\partial_x\bar{\alpha}_2(0,t)=0\,, \label{ch15.eq:sys_diff_theta_3}\\
&\bar{\alpha}_2(D_2,t)=\tilde{\theta}_2(t)\,,
\label{ch15.eq:def_vartheta_end}
\end{align}
where $\bar{\alpha}_i: [0\,, D_i] \times \mathbb{R}_+ \to \mathbb{R}$, $\forall i \in \{ 1, 2 \}$, and $\vartheta(t)\!:=\!\mathcal{H}[\tilde{\theta}(t)]\!=\!\hat{\Theta}(t)-\Theta^*$ is the propagated estimation error $\tilde{\theta}(t)$ through the transport-diffusion domain. For
$\lim\limits_{t\rightarrow \infty} \theta(t)\!=\!\theta_{c}$, we have $\lim\limits_{t\rightarrow \infty} \Theta(t) \!=\! \Theta_c \!=\! \theta_c$, where the index $c$ indicates a constant signal. 
Indeed, from (\ref{ch15.bucetasso}), for a constant input $\theta=\theta_c$, one has $\lim\limits_{t\rightarrow \infty} \mathcal{H}_1[\theta_{c1}] = \theta_{c1}$ for Player $P_1$. For Player $P_2$, we have $\mathcal{L}\lbrace\theta_{c2}\rbrace = \theta_{c2}/s$, and applying the Laplace limit theorem we get $\lim\limits_{t\rightarrow \infty} \mathcal{H}_2[\theta_{c2}] \!=\! \lim\limits_{s\rightarrow 0} \left\lbrace\frac{\theta_{c2}}{\cosh(\sqrt{s}D_2)} \right\rbrace \!=\! \theta_{c2}$. Thus, in the particular case $\theta\!=\!\theta_c\!=\!\theta^*$, one has 
%
\begin{equation} \label{ch15.meucucheirosa}
\Theta^*=\theta^*\,.
\end{equation}

Fig.~\ref{ch15.fig:blockDiagram_v3} depicts a schematic diagram that summarizes the proposed Nash Equilibrium policy for each player. While the signals $M_i(t)$, $N_i(t)$, $\hat{H}_i(t)$ and $\hat{G}_i(t)$, $\forall i\in\{1\,, 2\}$, are the same as in \textcolor{black}{Sections~``Noncooperative Scenario with Delays'' and ``Noncooperative Games with Heat PDEs'',}  
%
%
the additive dither signals in the presence of transport-heat PDE dynamics \cite{Tiago,FKKO:2018} are re-defined according to 
\begin{align} 
\begin{cases}
S_{1}(t)\!\!\!\!&=a_{1}\sin\left(\omega_{1}(t+D_{1})\right) \\
S_{2}(t)\!\!\!\!&=\frac{1}{2}a_{2}e^{\sqrt{\frac{\omega_{2}}{2}}D_{2}}\sin\left(\omega_{2}t+\sqrt{\frac{\omega_{2}}{2}}D_{2}\right)  \\
&\quad+\frac{1}{2}a_{2}e^{-\sqrt{\frac{\omega_{2}}{2}}D_{2}}\sin\left(\omega_{2}t-\sqrt{\frac{\omega_{2}}{2}}D_{2}\right)
\end{cases}  \,. \label{ch15.eq:Si}
\end{align}
Hence, computing the average of $\hat{G}_{1}(t)$ and $\hat{G}_{2}(t)$, leads us to
\begin{align}
\begin{cases}
\hat{G}_{1}^{\rm{av}}(t)=H_{11}^{1}\vartheta_{1}^{\rm{av}}(t)+\epsilon H_{12}^{1}\vartheta_{2}^{\rm{av}}(t) \\
\hat{G}_{2}^{\rm{av}}(t)=\epsilon H_{21}^{2}\vartheta_{1}^{\rm{av}}(t)+H_{22}^{2}\vartheta_{2}^{\rm{av}}(t)
\end{cases}\,.\label{ch15.eq:hatGi}
\end{align}
%

Additionally, from the block diagram in Fig.~\ref{ch15.fig:blockDiagram_v3}, one has
\begin{align}
\dot{\hat{\theta}}_i(t) = U_i(t)\,, \quad \dot{\tilde{\theta}}_i(t)=U_i(t)\,, \quad \forall i \in \{ 1, 2 \}  \,,
\label{ch15.eq:theta_hat_U}
\end{align}
since $\dot{\tilde{\theta}}(t)=\dot{\hat{\theta}}(t)$, once $\theta^*$ is constant. 
%
%
Taking the time derivative of \eqref{ch15.eq:def_vartheta_start_delay}--\eqref{ch15.eq:def_vartheta_end_delay} and \eqref{ch15.eq:def_vartheta_start}--\eqref{ch15.eq:def_vartheta_end}, with the help of \eqref{ch15.eq:def_theta_tilde} and \eqref{ch15.eq:theta_hat_U}, the \textit{propagated error dynamics} can be written as
\begin{align}
\label{ch15.eq:error_dyn_start_delay}
&\dot{\vartheta}_1(t)= u_1(0,t) = U(t-D_1) \,,   \\ 
&\partial_t u_1(x,t)=\partial_{x}u_1(x,t)\,, \quad x\in(0,D_1)\,, \\
&\partial_x u_1(0,t)=0\,, \\
&u_1(D_1,t)=U_1(t)\,,
\label{ch15.eq:error_dyn_end_delay}
\end{align}
and
\begin{align}
\label{ch15.eq:error_dyn_start}
&\dot{\vartheta}_2(t)=u_2(0,t) \,,  \\ 
&\partial_t u_2(x,t)=\partial_{xx}u_2(x,t)\,, \quad x\in(0,D_2)\,, \\
&\partial_x u_2(0,t)=0\,, \\
&u_2(D_2,t)=U_2(t)\,,
\label{ch15.eq:error_dyn_end}
\end{align}
where $u_i: [0\,, D_i] \times \mathbb{R}_+ \to \mathbb{R}$, $u_i(x,t):=\partial_t\bar{\alpha}_i(x,t)$, $\forall i \in \{ 1, 2 \}$,
and $\bar{\alpha}_i(x,t)=\alpha_i(x,t)-\beta_i(x,t)-\Theta_i^*$. The term $\beta_i(x,t)$ is the PDE state of the \textit{trajectory
generation problem} \cite[Chap.~12]{krstic2008boundary} solved to obtain $S_1(t)=\beta_1(D_1,t)$ and $S_2(t)=\beta_2(D_2,t)$ in (\ref{ch15.eq:Si})---for more details, see \cite[Eqs. (19) to (22)]{TAC:2020}.


Hence, from (\ref{ch15.eq:hatGi}), (\ref{ch15.eq:error_dyn_start_delay})--(\ref{ch15.eq:error_dyn_end_delay}) and (\ref{ch15.eq:error_dyn_start})--(\ref{ch15.eq:error_dyn_end}),  
it is possible to find a compact form for the overall average estimated gradient  according to 
\begin{align}
\hat{G}^{\rm{av}}(t)&=H \vartheta^{\rm{av}}(t)\,, \label{ch15.eq:hatGAv}\\
\dot{\hat{G}}^{\rm{av}}(t)&=H \dot{\vartheta}^{\rm{av}}(t)= H\mathcal{H}[U^{\rm{av}}(t)]\,, \label{ch15.eq:dotHatGAv}
\end{align}
where the Hessian $H$ is given in equation (\ref{ch15.eq:NE_v2}), $\vartheta^{\rm{av}}(t)\!:=\! [\vartheta_{1}^{\rm{av}}(t)\,,\vartheta_{2}^{\rm{av}}(t)]^T \!\in\!\mathbb{R}^{2}$, 
$\hat{G}^{\rm{av}}(t)\!:=\! [\hat{G}_{1}^{\rm{av}}(t)\,,\hat{G}_{2}^{\rm{av}}(t)]^T\!\in\!\mathbb{R}^{2}$ and 
$U^{\rm{av}}(t)\!:=\! [U_{1}^{\rm{av}}(t)\,,U_{2}^{\rm{av}}(t)]^T \!\in\!\mathbb{R}^{2}$ are the average versions of 
$U(t)\!:=\! [U_{1}(t)\,,U_{2}(t)]^T$, $\vartheta(t)\!:=\! [\vartheta_{1}(t)\,,\vartheta_{2}(t)]^T$ and
$\hat{G}(t)\!:=\! [\hat{G}_{1}(t)\,,\hat{G}_{2}(t)]^T$, respectively.



\subsection{Decentralized PDE Boundary Control using only the Known Diagonal Terms of the Hessian}


%

The control laws must be able to ensure exponential stabilization of $\hat{G}^{\rm{av}}(t)$ by compensating the heterogeneous transport-diffusion operator
$\mathcal{H}[\cdot]$ in (\ref{ch15.eq:dotHatGAv}). 
%
Following \cite{Tiago,FKKO:2018}, we propose the following boundary-based update laws $\dot{\hat{\theta}}_i(t)=U_i(t)$, $i\in\{1\,,2\}$:
%
%
\begin{equation} \label{ch15.4_dU_NC}
\begin{cases}
	\dot{U}_{1}(t) \!\!\!\!&=\!\! -c_{1}U_{1}(t)  \\
 &\quad+ c_{1}k_{1} \left(\hat{G}_{1}(t) +\hat{H}_{1}(t) \int_{0}^{D_1} u_1(\tau,t)d\tau \right)\\
	\dot{U}_{2}(t) \!\!\!\!&=\!\! -c_{2}U_{2}(t)  \\
 &\quad+ c_{2}k_{2} \left(\hat{G}_{2}(t) +\hat{H}_{2}(t) \int_{0}^{D_2}(D_2-\tau) u_2(\tau,t)d\tau \right)
\end{cases},
\end{equation}
%
for positive constants $k_{1}>0$, $k_{2}>0$, $c_{1}>0$ and $c_{2}>0$, in order to compensate for the transport-heat PDEs in (\ref{ch15.eq:error_dyn_start_delay})--(\ref{ch15.eq:error_dyn_end_delay}) and (\ref{ch15.eq:error_dyn_start})--(\ref{ch15.eq:error_dyn_end}). As before two sections and with some abuse of notation, constants $c_1$ and $c_2$ were chosen to denote the parameters of the control laws, but they have no relation to those which appear in the payoffs given by (\ref{ch15.eq:J1}) and (\ref{ch15.eq:J2}). 
%
%

The boundary control law (\ref{ch15.4_dU_NC}) could be rewritten as
\begin{align}\label{ch15.predictor_new}
\dot{U}_{1}(t) \!&=\!\! -c_{1}U_{1}(t) \nonumber \\
&\quad+ c_{1}k_{1} \left(\hat{G}_{1}(t) \!+\!\hat{H}_{1}(t) \int_{t-D_1}^{t} U_1(\tau,t)d\tau \right), \nonumber \\
\dot{U}_2(t)\!&=\!\!-c_2 U_2(t) \nonumber \\
&\quad+ c_2k_2 \!\left[\hat{G}_2(t)\!+\!\hat{H}_2(t)\!\left(\hat{\theta}_2(t)\!-\!\Theta_2(t)\!+\!a_2 \sin (\omega_2 t)\right)\right], 
\end{align}
using the relation $u_1(x,t)=U_1(t+x-D_1)$ for the transport PDE and the diffusion equation $\partial_t\alpha_2(x,t)\!=\!\partial_{xx}\alpha_2(x,t)$ as well as the integration by parts, associated with
(\ref{ch15.eq:heat_eqn_start})--(\ref{ch15.eq:heat_eqn_end}), (\ref{ch15.eq:def_theta_tilde}) and recalling that $\vartheta_2+a_2\sin(\omega_2 t)=
\Theta_2(t)-\Theta_2^*$, analogously to \cite[Eq. (25)]{FKKO:2018}.

\subsection{ISS-Like Properties for Hyperbolic-Parabolic PDE\\ Representation}

For the sake of simplicity, we assume $c_{1}\,,c_{2} \!\to\! +\infty$ in (\ref{ch15.4_dU_NC}), resulting in the following general expressions:
\begin{equation}
\begin{aligned} \label{ch15.rabao_NC}
U_{1}(t)&=k_{1}\Bigg(\hat{G}_{1}(t) +\hat{H}_{1}(t) \int_{0}^{D_1} u_1(\tau,t)d\tau\Bigg), \nonumber \\
U_{2}(t)&=k_{2}\Bigg(\hat{G}_{2}(t) +\hat{H}_{2}(t) \int_{0}^{D_2}(D_2-\tau) u_2(\tau,t)d\tau\Bigg).
\end{aligned}
\end{equation}
Recalling (\ref{ch15.eq:error_dyn_start_delay})--(\ref{ch15.eq:error_dyn_end_delay}) and (\ref{ch15.eq:error_dyn_start})--(\ref{ch15.eq:error_dyn_end}), the infinite-dimensional closed-loop system (\ref{ch15.eq:dotHatGAv}) and (\ref{ch15.rabao_NC}) in its average version can be written in the corresponding PDE representation form, given by
\begin{eqnarray}
\dot{\hat{G}}_1^{\rm{av}}(t) &=& H_{11}^{1} u_{1}^{\rm{av}}(0,t) +\epsilon H_{12}^{1}u_{2}^{\rm{av}}(0,t) \,, \label{ch15.saco1_NC_delay}\\
\partial_t u_1^{\rm{av}}(x,t)&=&D_1^{-1}\partial_{x} u_1^{\rm{av}}(x,t)\,, \quad x\in (0\,,1)\,, \label{ch15.saco2_NC_delay}\\
\partial_x u_1^{\rm{av}}(0,t)&=&0\,, \\
u_1^{\rm{av}}(1\,,t)&=&U_1^{\rm{av}}(t) \label{ch15.saco3_NC_delay}\,,
\end{eqnarray}
and
\begin{eqnarray}
\dot{\hat{G}}_2^{\rm{av}}(t) &=& \epsilon H_{21}^{2}u_{1}^{\rm{av}}(0,t)+H_{22}^{2}u_{2}^{\rm{av}}(0,t) \,, \label{ch15.saco1_NC}\\
\partial_t u_2^{\rm{av}}(x,t)&=&D_2^{-2}\partial_{xx} u_2^{\rm{av}}(x,t)\,, \quad x\in (0\,,1)\,, \label{ch15.saco2_NC}\\
\partial_x u_2^{\rm{av}}(0,t)&=&0\,, \\
u_2^{\rm{av}}(1\,,t)&=&U_2^{\rm{av}}(t) \label{ch15.saco3_NC}\,.
\end{eqnarray}
%

In the \textit{reduction-like approach} \cite{a14} (or finite-spectrum assignment), we use the following transformations to write:
\begin{align} 
\bar{G}^{\rm{av}}_1(t)\!&=\!\hat{G}^{\rm{av}}_1(t)\!+\! \epsilon_{11}^{1} H_{11}^{1} \!\int_{0}^{D_1}\!\!\!u^{\rm{av}}_1(\tau,t) d\tau \nonumber \\ 
&\quad\!+\! \epsilon_{12}^{1} H_{12}^{1} \!\int_{0}^{D_2}\!\!\!(D_2\!-\!\tau) u^{\rm{av}}_2(\tau,t) d\tau~~~~
\label{ch15.transformacao_cu_NC_delay}
\end{align}
and
\begin{align} 
\bar{G}^{\rm{av}}_2(t)\!&=\!\hat{G}^{\rm{av}}_2(t)\!+\! \epsilon_{22}^{2} H_{22}^{2} \!\int_{0}^{D_2}\!\!\!(D_2\!-\!\tau) u^{\rm{av}}_2(\tau,t) d\tau \nonumber \\
&\quad\!+\! \epsilon_{21}^{2} H_{21}^{2} \!\int_{0}^{D_1}\!\!\! u^{\rm{av}}_1(\tau,t) d\tau\,,~~~~
\label{ch15.transformacao_cu_NC}
\end{align}
%
where $\epsilon_{11}^{1}=\epsilon_{22}^{2}=1$ and $\epsilon_{12}^{1}=\epsilon_{21}^{2}=\epsilon$. 

%
With some mathematical manipulations, it is not difficult to see that $\bar{G}^{\rm{av}}$ satisfies
\begin{eqnarray}  \label{ch15.Zreduction_NC}
\dot{\bar{G}}^{\rm{av}}(t)=H U^{\rm{av}}(t)\,.
\end{eqnarray}
After adding and subtracting
the next terms in blue and red in (\ref{ch15.rabao_NC}), it can be rewritten as:
%
\begin{align} \label{ch15.rabao_NC_1}
&U_{1}(t)=k_{1}\left(\hat{G}_{1}(t) +\hat{H}_{1}(t) \int_{0}^{D_1}  u_1(\tau,t)d\tau\right. \nonumber \\
&\quad\left.+\textcolor{blue}{\epsilon H_{12}^{1}\int_{0}^{D_2} (D_2-\tau) u_2(\tau,t) d\tau}\right)  \nonumber\\
&~~~-k_{1}\textcolor{red}{\epsilon H_{12}^{1}\int_{0}^{D_2} (D_2-\tau) u_2(\tau,t) d\tau}\,, \\
&U_{2}(t)=k_{2}\left(\hat{G}_{2}(t) +\hat{H}_{2}(t) \int_{0}^{D_2} (D_2-\tau) u_2(\tau,t)d\tau\right. \nonumber \\
&\quad\left.+\textcolor{blue}{\epsilon H_{21}^{2}\int_{0}^{D_1}  u_1(\tau,t) d\tau}-\textcolor{red}{\epsilon H_{21}^{2}\int_{0}^{D_1}  u_1(\tau,t) d\tau}\right) \,, 
\end{align}
whose average compact form is
\begin{align} 
U^{\rm{av}}(t)&=K\bar{G}^{\rm{av}}(t)+\epsilon K\phi^{\rm{av}}(D,t) 
\,, \label{ch15.rabao_NC_1_cu}
\end{align}
where the matrix $K:=\text{diag}\{k_{1}\,,k_{2}\}$ with entries $k_1>0$, $k_2>0$ and the auxiliary variable $\phi(D,t)$ is defined as
\begin{align}
\phi(D,t)&:=-\begin{bmatrix}
	H_{12}^{1}\int_{0}^{D_2} (D_2-\tau) u_2(\tau,t) d\tau \\
  H_{21}^{2}\int_{0}^{D_1} u_1(\tau,t) d\tau
\end{bmatrix}\,, \nonumber \\
\phi(1,t)&:=-\begin{bmatrix}
	H_{12}^{1}\int_{0}^{1}D_2^2 (1-\xi) u_2(\xi,t) d\xi \\
  H_{21}^{2}\int_{0}^{1} D_1 u_1(\xi,t) d\xi
\end{bmatrix}\,,
\label{ch15.eq:phii} 
\end{align}
since $\int_{0}^{D_j}(D_j-\tau) u_j(\tau,t) d\tau=\int_{0}^{1} D_j^2(1-\xi) u_j(\xi,t) d\xi$, for $j \in \{1\,,2\}$. 
Then, it is possible to find a compact form for the overall average game from (\ref{ch15.Zreduction_NC}) and
(\ref{ch15.rabao_NC_1_cu}), such as
\begin{eqnarray}
\dot{\bar{G}}^{\rm{av}}(t) &=& HK \bar{G}^{\rm{av}}(t)+\epsilon HK \phi^{\rm{av}}(1,t)\,, \label{ch15.saco1M_NC}\\
\partial_t u_1^{\rm{av}}(x,t)&=&D_1^{-1}\partial_{x} u_1^{\rm{av}}(x,t)\,, \quad ~x\in (0\,,1)\,, \label{ch15.saco2M_NC_delay}\\
\partial_t u_2^{\rm{av}}(x,t)&=&D_2^{-2}\partial_{xx} u_2^{\rm{av}}(x,t)\,, \!\quad x\in (0\,,1)\,, \label{ch15.saco2M_NC}\\
\partial_x u^{\rm{av}}(0,t)&=&0\,, \\
u^{\rm{av}}(1\,,t)&=&K \bar{G}^{\rm{av}}(t)+\epsilon K \phi^{\rm{av}}(1,t) \label{ch15.saco3M_NC}\,.
\end{eqnarray}
%

From (\ref{ch15.saco1M_NC}), if $HK$ is Hurwitz, it is clear that the dynamics of the ODE state variable $\bar{G}^{\rm{av}}(t)$ is exponentially Input-to-State Stable (ISS) \cite{KK:2018} with respect to the PDE state $u^{\rm{av}}(x,t)=[u_1^{\rm{av}}(x,t)\,,u_2^{\rm{av}}(x,t)]^T$ by means of the function $\phi^{\rm{av}}(1,t)$. Moreover, the PDE subsystem (\ref{ch15.saco2M_NC_delay}) and (\ref{ch15.saco2M_NC}) is ISS \cite{KK:2018} with respect to $\bar{G}^{\rm{av}}(t)$ in the boundary condition $u^{\rm{av}}(1\,,t)$.

\subsection{Stability Analysis}
\label{doritos_stability_hetero}

Next, we will show here that this hyperbolic-parabolic PDE-ODE loop (\ref{ch15.saco1M_NC})--(\ref{ch15.saco3M_NC}) contains a small-parameter $\epsilon$ which can lead to closed-loop stability if it is chosen sufficiently small. To this end, we assume the following particular condition for duopoly games \cite{FKB:2012}, which is equivalent to  Assumptions \ref{ch.13.Assumption 1.} and \ref{ch13.Assumption 2.} when the payoff functions of the form of (\ref{ch15.eq:J1}) and (\ref{ch15.eq:J2}) are considered.


\begin{assum} \label{ch17.Assumption 1.nene}
The Hessian matrix $H$ given by (\ref{ch15.eq:NE_v2}) is strictly diagonal dominant, \textit{i.e.}, 
\begin{equation}
\begin{cases}
|\epsilon H_{12}^{1}| < |H_{11}^{1}|\\
|\epsilon H_{21}^{2}| < |H_{22}^{2}|\\    
\end{cases}\,.
\end{equation}
\end{assum}
%

\smallskip



%


The next theorem provides the stability/convergence properties of the closed-loop error system of the proposed ES feedback for the $2$-player noncooperative game with transport-heat PDEs. 


\begin{theorem} \label{duoply_heterogeneous_queirozBday} Consider the closed-loop system (\ref{ch15.eq:error_dyn_start_delay})--(\ref{ch15.eq:error_dyn_end}) under transport-heat PDEs (\ref{ch15.eq:heat_eqn_start_delay})--(\ref{ch15.eq:heat_eqn_end}) of distinct transport-diffusion coefficients $D_1$ and $D_{2}$ for the heterogeneous duopoly quadratic game with payoff functions $y_i(t) \!=\! J_i(\Theta(t))$, $\forall i \in \{ 1, 2 \}$, given in (\ref{ch15.eq:J1}) and (\ref{ch15.eq:J2}), satisfying Assumption~\ref{ch17.Assumption 1.nene} and control laws $U_i(t)$ defined in (\ref{ch15.4_dU_NC}) or (\ref{ch15.predictor_new}).
Thus, there exist $c_1>0$, $c_2>0$ and $\omega>0$ sufficiently large as well as $\epsilon>0$ sufficiently small such that the closed-loop error system  
with state $\vartheta_i(t)$, $u_i(x,t)$, $\forall i\in \{1\,,2\}$, has a unique locally exponentially stable periodic solution in $t$ of period $\Pi$ in (\ref{period}),
denoted by $\vartheta_i^{\Pi}(t)$, $u_i^{\Pi}(x,t)$ and satisfying, $\forall t\geq 0$:
%
\begin{eqnarray}\label{ch15.periodic_solution_NCG}
\left(\sum_{i=1}^{2}\left[\vartheta_{i}^{\Pi}(t)\right]^2  \!+\! \int_{0}^{D_i} \left[u_{i}^{\Pi}(x,t)\right]^2dx \right)^{1/2} \!\leq\! \mathcal{O}(1/\omega)\,.
\end{eqnarray}
%
In particular,
\begin{eqnarray}\label{ch15.limsup1_NCG}
\limsup_{t\to+\infty}|\Theta(t)\!-\!\Theta^*|\!&=&\!\mathcal{O}(|a|\!+\!1/\omega),
\\
\limsup_{t\to+\infty}|\theta_1(t)\!-\!\theta_1^*|\!&=&\!\mathcal{O}\left(a_1\!+\!1/\omega\right), \label{ch15.limsup2}
\\
\limsup_{t\to+\infty}|\theta_2(t)\!-\!\theta_2^*|\!&=&\!\mathcal{O}\left(a_2e^{D_2\sqrt{\omega/2}}\!+\!1/\omega\right), \label{ch15.limsup3}
\end{eqnarray}
where $a=[a_1 \ a_2]^T$ and $\theta^*\!=\!\Theta^*$ is the unique (unknown) Nash Equilibrium given by (\ref{ch15.eq:NE_v2}).
\end{theorem}

\begin{proof}
See the Appendix.
\end{proof}

\section{Simulations}
\label{doritos_simulation}

Due to space limitations, we will restrict ourselves to a unique numerical example considering a heterogeneous noncooperative game with $2$ players that employ the proposed ES strategy for PDE compensation of \textcolor{black}{Section~``Noncooperative Duopoly with Heterogeneous Transport-Heat PDE Dynamics''}. We
revisit the example in \cite{ACC:2021}, and consider the following payoff functions (\ref{ch15.eq:J1}) and (\ref{ch15.eq:J2}) subject to transport-heat PDEs  (\ref{ch15.eq:heat_eqn_start_delay})--(\ref{ch15.eq:heat_eqn_end_delay}) and  (\ref{ch15.eq:heat_eqn_start})--(\ref{ch15.eq:heat_eqn_end}) with distinct transport-diffusion coefficients $D_1=30$ and $D_2=3$ in the players' decisions, $i\in\{1\,,2\}$:
%
\begin{align} \label{ch15.J1_example} 
J_1(\Theta(t))&=-5~\Theta_1^2(t)+5~\epsilon \Theta_1(t)\Theta_2(t)+250~\Theta_1(t) \nonumber \\
&\quad-150~\Theta_2(t)-3000 \,,\\
J_2(\Theta(t))&=-5~\Theta_2^2(t)+5~\epsilon \Theta_1(t)\Theta_2(t)-150~\Theta_1(t) \nonumber \\
&\quad+150~\Theta_2(t)+2500\,, \label{ch15.J2_example}
\end{align}
%
which, according to (\ref{ch15.eq:NE_v2}), yield the unique Nash equilibrium 
\begin{align} 
\Theta_{1}^{*}=\theta_{1}^{*}=\frac{100+30\epsilon}{4-\epsilon^{2}}\,, \quad
\Theta_{2}^{*}=\theta_{2}^{*}=\frac{60+50\epsilon}{4-\epsilon^{2}}\,. \label{ch15.Nash_equilibrium_1_2}
\end{align}
In order to attain (\ref{ch15.Nash_equilibrium_1_2}), 
the players implement the non-model based real-time optimization strategy acting through the transport-heat PDE dynamics (see Fig.~\ref{ch15.fig:blockDiagram_v3}). 
%
%
%
%
%
%
%
%
For comparison purposes, except for the transport-heat PDEs in the players' input signals, the plant and
the controller parameters were chosen similarly to \cite{FKB:2012} in the simulation tests:
$\epsilon\!=\!1$, $a_1\!=\!0.075$, $a_2\!=\!0.05$, $k_1\!=\!2$, $k_2\!=\!5$, $\omega_1\!=\!26.75~$rad/s, $\omega_2\!=\!22~$rad/s and
$\theta_1(0)\!=\!\hat{\theta}_1(0)\!=\!50$, $\theta_2(0)\!=\!\hat{\theta}_2(0)\!=\!\theta_2^*\!=\!110/3$. In addition, the time-constants of the boundary control filters were set to $c_{1}\!=\!c_{2}\!=\!100$.  

\begin{figure}[ht]
	\centering
	\subfigure[Actions' time histories.\label{ch15.fig:orkb2020_u_delays}]{\includegraphics[width=6.5cm]{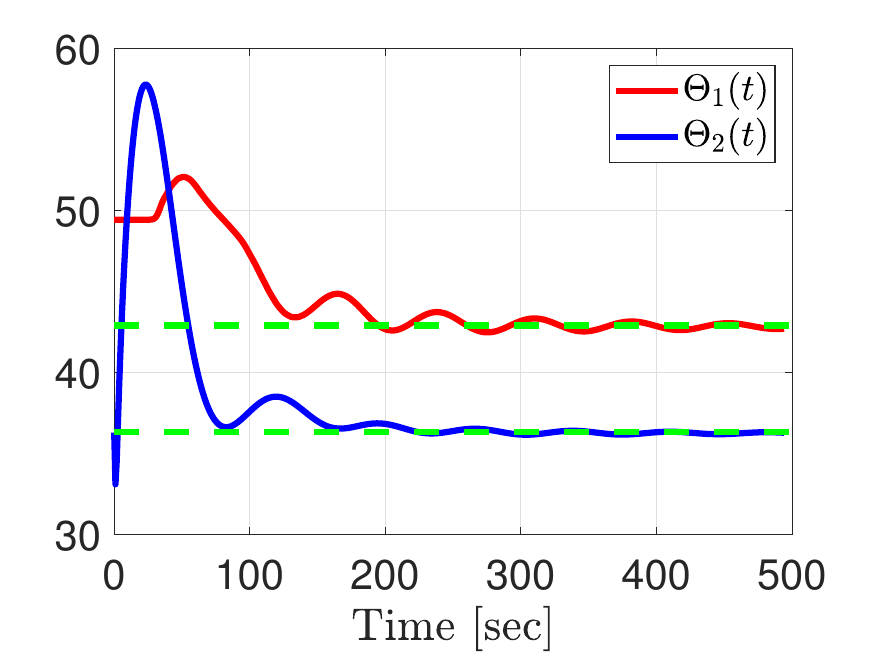}}
	~
	\subfigure[Payoffs' time histories. \label{ch15.fig:orkb2020_y_delays}]{\includegraphics[width=6.5cm]{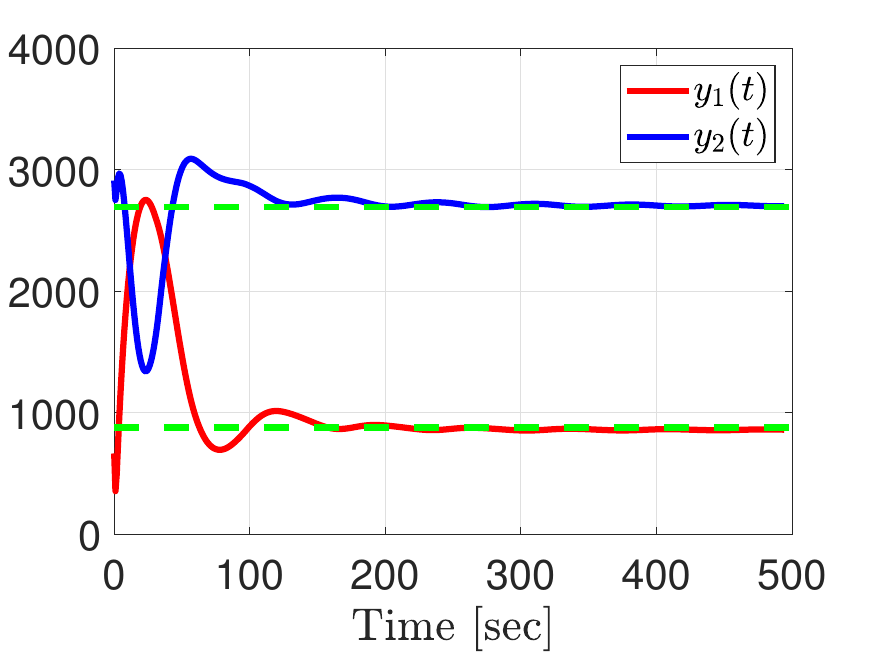}}
	\caption{(a) actions' time histories and (b) payoffs' time histories for $P_1$ and $P_2$, 
	for $\epsilon=1$. The dashed lines denote the values at the Nash Equilibrium, $\Theta^{*}_{1}=43.33$ and $\Theta^{*}_{2}=36.67$ (with $J_1(\Theta^*)=889$ and $J_2(\Theta^*)=2722$).}
\end{figure} 

We can check that the ES approach proposed in \cite{FKB:2012} is effective when transport-heat PDEs are not present in the decision variables. 
However, in the presence of the transport-heat PDEs in the input signals $\theta_1$ and $\theta_2$, but without considering
any kind of PDE compensation, 
the game collapses with the explosion of its variables (curves not shown). On the other hand, 
Fig.~\ref{ch15.fig:orkb2020_u_delays} and Fig.~\ref{ch15.fig:orkb2020_y_delays}
show that the proposed boundary control based scheme fixes this with
a remarkable evolution in searching the Nash Equilibrium and simultaneously compensating for the effect of the transport-heat PDEs in our heterogeneous noncooperative game.

The evolutions of the infinite-dimensional states $\alpha_1(x,t)$ and $\alpha_2(x,t)$ modeled by the transport-heat PDEs (\ref{ch15.eq:heat_eqn_start_delay})--(\ref{ch15.eq:heat_eqn_end_delay}) and (\ref{ch15.eq:heat_eqn_start})--(\ref{ch15.eq:heat_eqn_end}) are shown in
Fig.~\ref{ch15.fig:orkb2020_u_delays_3D1} to Fig.~\ref{ch15.fig:orkb2020_y_delays_3D2}. The values of the boundary inputs $\theta_1(t)$ and $\theta_2(t)$ as well as the boundary outputs $\Theta_1(t)$ and $\Theta_2(t)$ are highlighted in colors black and red, respectively. The initial condition is in blue color.   

\begin{figure}[ht!]
	\centering
	\subfigure[\textcolor{black}{Action plots for $\epsilon\!=\!0.75$. The dashed lines denote
	$\Theta^{*}_{1}\!=\!35.64$ and $\Theta^{*}_{2}\!=\!28.36$.}
	\label{ch15.fig:orkb2020_u_delays_epsilon0_1}]{\includegraphics[width=6.75cm]
{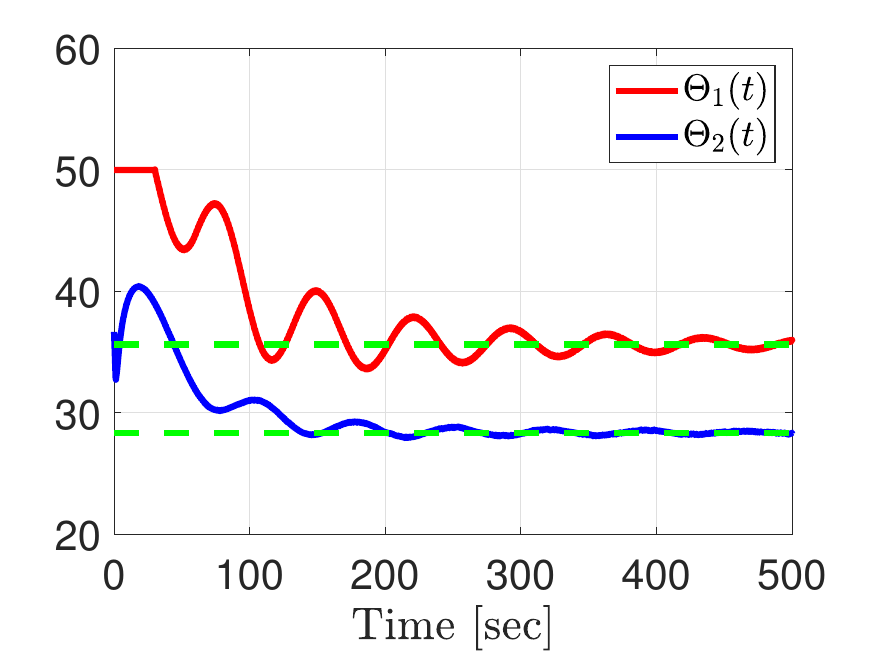}}
	~
	\subfigure[\textcolor{black}{Action plots for $\epsilon\!=\!0.5$. The dashed lines denote $\Theta^{*}_{1}\!=\!30.67$ and $\Theta^{*}_{2}\!=\!22.67$.}
	\label{ch15.fig:orkb2020_u_delays_epsilon0_5}]{\includegraphics[width=6.75cm]
{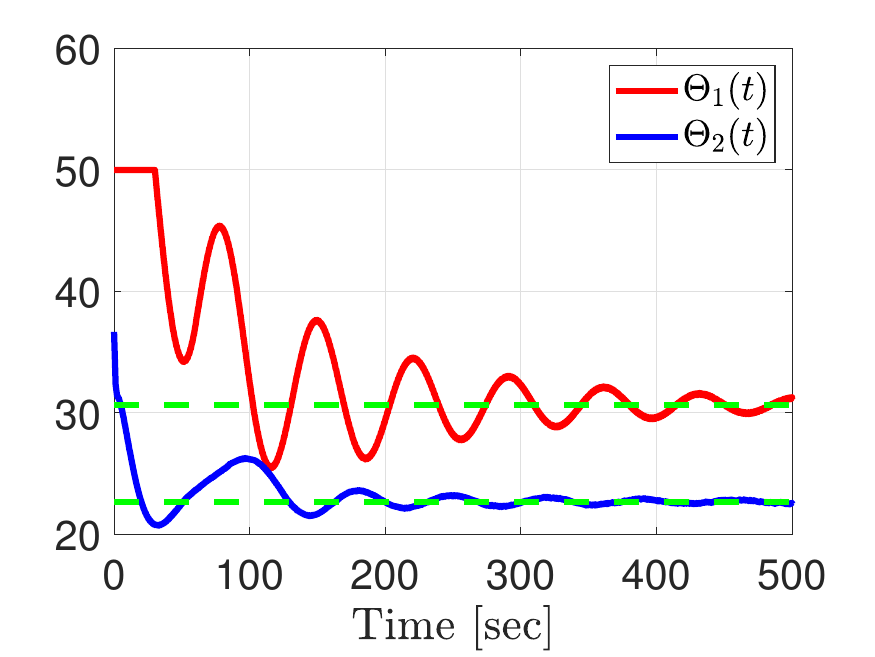}}
	~
	\subfigure[\textcolor{black}{Action plots for $\epsilon\!=\!0.25$. The dashed lines denote $\Theta^{*}_{1}\!=\!27.30$ and $\Theta^{*}_{2}\!=\!18.41$.}
	\label{ch15.fig:orkb2020_u_delays_epsilon0_01}]{\includegraphics[width=6.75cm]
{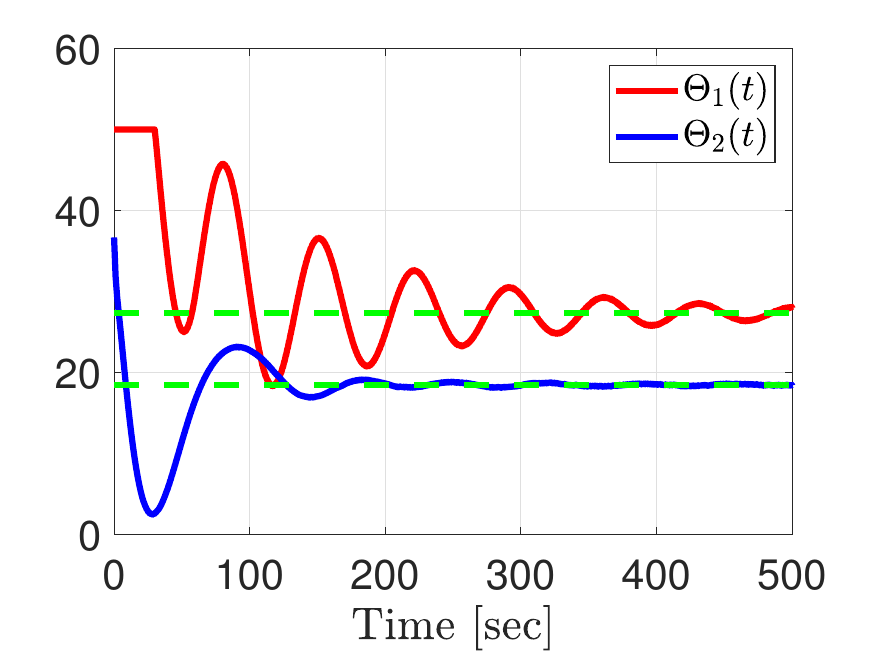}}
	\caption{Actions' time histories for $P_1$ and $P_2$
	for  distinct values of the coupling coefficient ($\epsilon=0.75$, $\epsilon=0.5$ and $\epsilon=0.25$).}
\end{figure}

\begin{figure*}[ht]
	\centering
	\subfigure[Parameter $\Theta_1(t)$ (red) converges to a $\mathcal{O}(|a|\!+\!1/\omega)$--neighborhood of $\Theta_1^*$ (dashed-green).\label{ch15.fig:orkb2020_u_delays_3D1}]{\includegraphics[width=7.0cm]{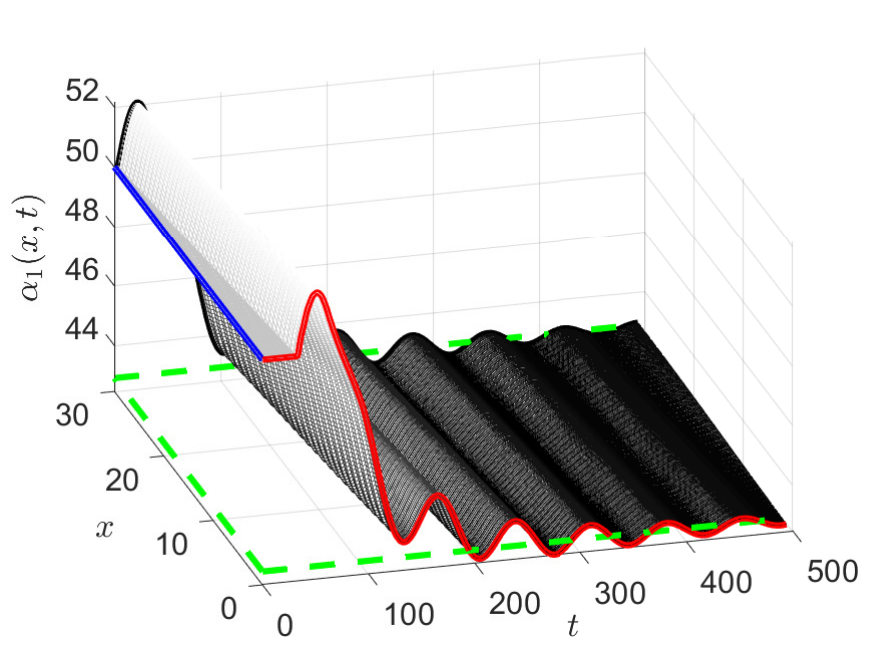}}
	~
	\subfigure[Parameter $\theta_1(t)$ (black) converges to a $\mathcal{O}(a_1\!+\!1/\omega)$--neighborhood of $\theta_1^*$ (dashed-green). \label{ch15.fig:orkb2020_y_delays_3D1}]{\includegraphics[width=7.0cm]{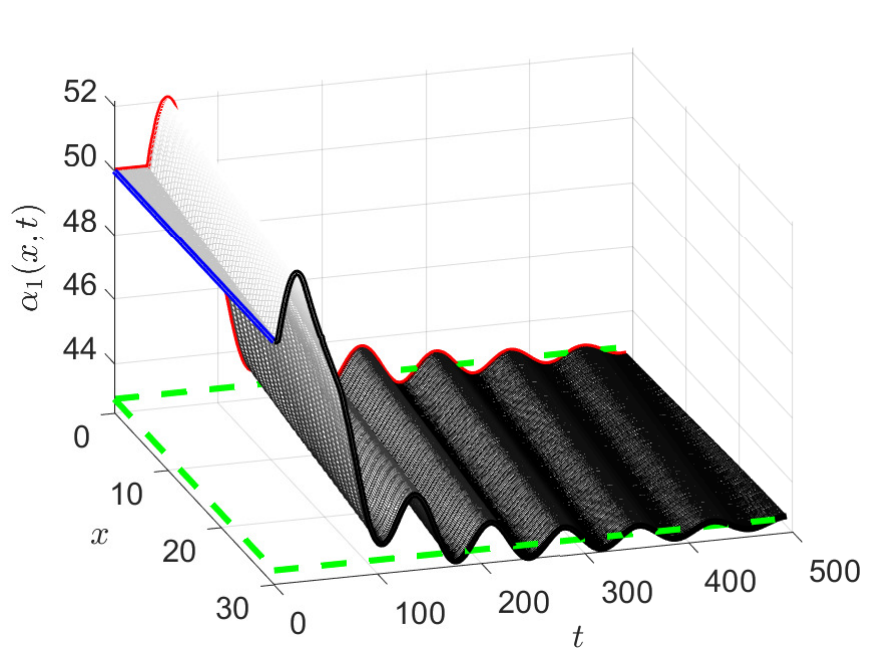}}
	~
	\subfigure[Parameter $\Theta_2(t)$ (red) converges to a $\mathcal{O}(|a|\!+\!1/\omega)$--neighborhood of $\Theta_2^*$ (dashed-green). \label{ch15.fig:orkb2020_u_delays_3D2}]{\includegraphics[width=7.0cm]{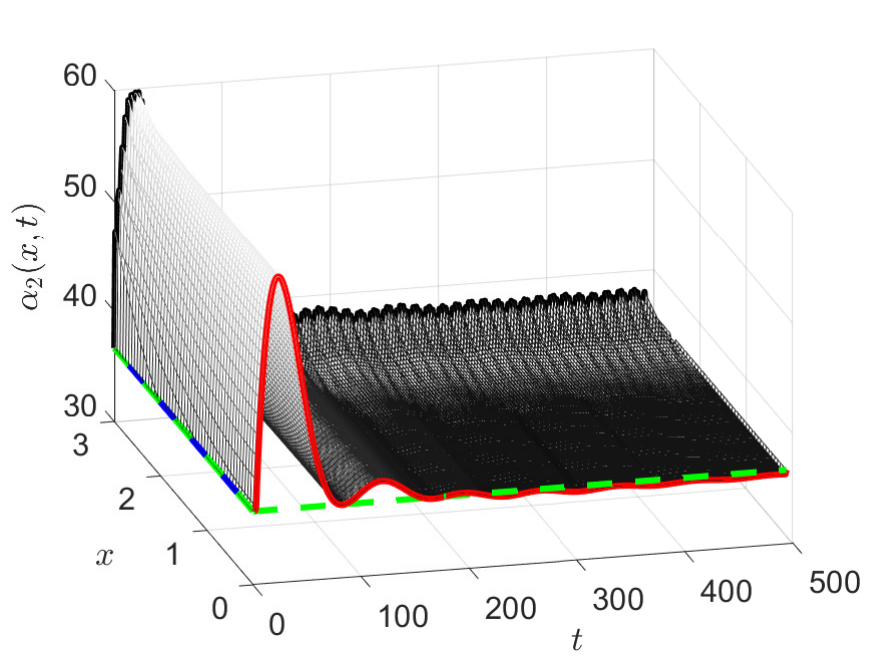}}
	~
	\subfigure[Parameter $\theta_2(t)$ (black) converges to a $\mathcal{O}(a_2e^{D_2\sqrt{\omega/2}}\!+\!1/\omega)$--neighborhood of $\theta_2^*$ (dashed-green). \label{ch15.fig:orkb2020_y_delays_3D2}]{\includegraphics[width=7.0cm]{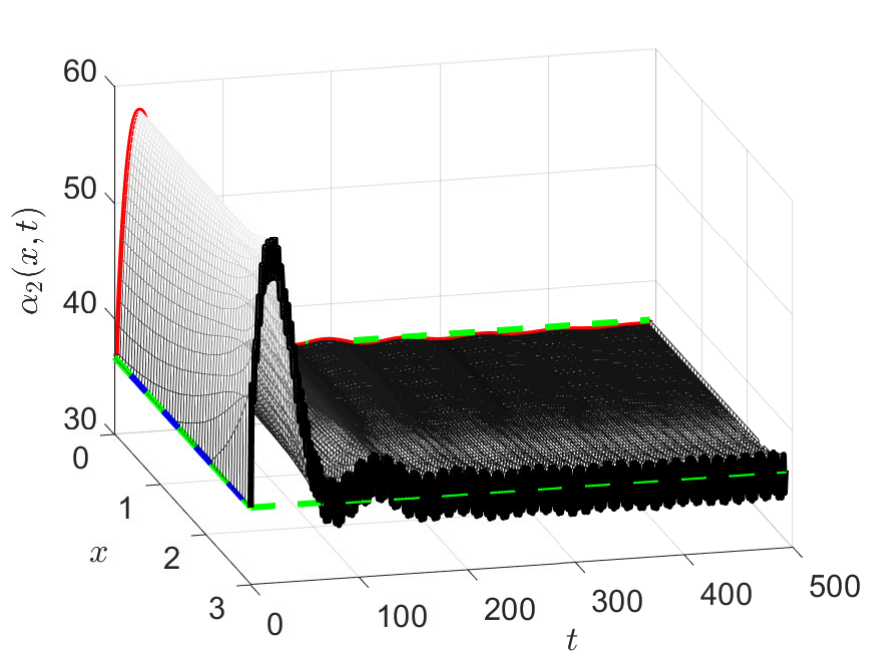}}
	\caption{Evolution of the infinite-dimensional states $\alpha_1(x,t)$ and $\alpha_2(x,t)$ of the transport-heat PDEs in a heterogeneous duopoly game with boundary Dirichlet actuation, according to (\ref{ch15.limsup1_NCG})--(\ref{ch15.limsup3}): from $\alpha_1(D_1,t)\!=\!\theta_1(t)$ to $\alpha_1(0,t)\!=\!\Theta_1(t)$, with $D_1\!=\!30$ for Player $P_1$ and from
	$\alpha_2(D_2,t)\!=\!\theta_2(t)$ to $\alpha_2(0,t)\!=\!\Theta_2(t)$, with $D_2\!=\!3$ for Player $P_2$.}
\end{figure*}

This first set of simulations indicates that even under an adversarial scenario of strong coupling between the players with
$\epsilon\!=\!1$, the proposed approach has behaved successfully. This suggests that our stability analysis may be conservative and the theoretical assumption
$0\!<\!\epsilon\!<\!1$ may be relaxed given the performance of the closed-loop 
control system. 
In Fig.~\ref{ch15.fig:orkb2020_u_delays_epsilon0_1}, Fig.~\ref{ch15.fig:orkb2020_u_delays_epsilon0_5} and Fig.~\ref{ch15.fig:orkb2020_u_delays_epsilon0_01}, different values of
$\epsilon\!=\!0.75$, $\epsilon\!=\!0.5$ and
$\epsilon\!=\!0.25$ are considered to evaluate the robustness of the proposed scheme under different levels of coupling between the two players and the corresponding impact on the transient responses. 

\section{From Nash Equilibrium to Extremum Seeking}
\label{label.CAP13.final}

We now present the scalar extremum seeking design for distinct families of PDEs. According to Remark~\ref{finalizandoremark1_1}, extremum seeking can be interpreted as a particular case of Nash Equilibrium Seeking.

\begin{remark}\label{finalizandoremark1_1}
\em Note that the material in the previous sections addressed the following special cases: 
\begin{itemize}
    \item Multivariable static extremum seeking \cite{Tiago,TAC:2020}, as corollaries of Theorems~\ref{ch16.theorem.16.1.vacina}, \ref{scheme2} and \ref{duoply_heterogeneous_queirozBday}, when $h_1=\dots=h_N$ and $\epsilon=1$, and
    \item Scalar static extremum seeking, when simply $N=1$.
\end{itemize}
\end{remark}


\subsection{Scalar Extremum Seeking as a Particular Case of Nash Equilibrium Seeking with One Single Agent} 

In this case, the ES goal is to optimize in real-time an unknown static map $Q(\cdot)$:  
\begin{equation}
y(t) = Q(\Theta(t)), 
\label{extrachapter.eq:initial_output_static_map_CSm}
\end{equation}
with maximum or minimum unknown output $y^{*}$ and optimizer $\Theta^{*}$, by measuring the output $y(t)$ and adjusting the input $\Theta(t)$. 
As illustrated in the block diagram of Fig.~\ref{ch12.fig:cascade}, the input $\Theta$ of the map is 
governed by a PDE, which must be properly compensated by means of an appropriate boundary control law.

For maximization problems, the unknown nonlinear map is assumed to be locally quadratic, such that
\begin{equation}
Q(\Theta) = y^{*} + \dfrac{H}{2}(\Theta - \Theta^{*})^{2},
\label{extrachapter.eq:static_map_CSm}
\end{equation}
where $\Theta^{*},y^{*} \in \mathbb{R}$ and the Hessian $H<0$ is also unknown (when $H>0$, we have a minimization problem). Hence, the output of the static map is given by
\begin{equation}
y(t) = y^{*} + \dfrac{H}{2}(\Theta(t) - \Theta^{*})^{2}.
\label{extrachapter.eq:final_output_static_map_CSm}
\end{equation}
For the sake of simplicity, we have considered scalar maps with one single-input and one single-output, but a more general setup of MISO maps ($\Theta \in \mathbb{R}^{N}$ and $y\in \mathbb{R}$) could also be addressed, as done in \cite{Tiago,TAC:2020}.

While the multiplicative perturbation signals $M(t)=\frac{2}{a}\sin(\omega t)$ and $N(t)=-\frac{8}{a^2}\cos(\omega t)$ follow classical ES designs \cite{KW:00,GKN:2012}, the additive dither $S(t)$ must be redesigned using the trajectory generation paradigm \cite[Chapter~12]{krstic2008boundary}. In the next, we show the ES design for the simplest cases of hyperbolic and parabolic PDEs: transport (delay) and
heat (diffusion) PDEs. While, in the case of a delay, it suffices to advance in time the sinusoidal perturbation, in the case of the heat dynamics one has to employ a solution to the motion planning problem where the output of the heat PDE system is a sinusoid at one of its boundaries and the input is the signal that must be applied on the other boundary in order to generate a sinusoid at the output. This input signal happens to consist of sinusoidal and exponential functions.

\subsection{Extremum Seeking for Transport Hyperbolic PDE}
In this case, the following infinite-dimensional and averaging-based predictor feedback is introduced in order to compensate the delay \cite{K:2009}

\begin{equation}\label{predictor}
U(t)= \frac{c}{s+c}\left\{k \left[G(t)+\hat{H}(t)\int_{t-D}^{t} U(\tau) d\tau \right] \right\}\,,
\end{equation}
where $k>0$ and $c>0$ is sufficiently large, \textit{i.e.}, the predictor feedback is of the form of a low-pass filtered of the non average version of 
\begin{align}\label{predictoraverage}
U_{\rm{av}}(t) &= k G_{\rm{av}}(t+D)\\ \nonumber
 &= k \left[G_{\rm{av}}(t)+H\int_{t-D}^{t}U_{\rm{av}}(\sigma) d\sigma
\right]\,.
\end{align}
The signal
\begin{equation}\label{hessian_estimate_H}
\hat{H}(t)=N(t)y(t) 
\end{equation}
in (\ref{predictor}) is used to obtain an estimate of the unknown Hessian $H$ and the additive dither is given by
\begin{eqnarray}
S(t) = a \sin(\omega (t+D))\,. \label{dither_signalS}
\end{eqnarray}

\subsection{Extremum Seeking for Diffusion-Dominated Parabolic PDE}

In the case of a diffusion process in the actuation dynamics, the trajectory generation problem for motion planning design \cite[Chapter~12]{krstic2008boundary} to be solved is:
\begin{align}
\label{eq:def_S_start}
&S(t):=\beta(D,t) \\ 
&\partial_t\beta(x,t)=\partial_{xx}\beta(x,t), \quad x\in(0,D) \\
&\partial_x\beta(0,t)=0 \\
&\beta(0,t)=a\sin(\omega t),
\label{eq:def_S_end}
\end{align}
where $\beta:[0,D]\times\mathbb{R}_+\rightarrow \mathbb{R}$. The explicit solution of \eqref{eq:def_S_start} is given by
\begin{align}
S(t) &= \frac{1}{2}ae^{\sqrt{\frac{\omega}{2}}D}\sin\left(\omega t \!+\! \sqrt{\frac{\omega}{2}D}\right) \nonumber \\
&+\frac{1}{2}ae^{-\sqrt{\frac{\omega}{2}}D}\sin\left(\omega t \!-\! \sqrt{\frac{\omega}{2}D}\right). 
\label{eq:sol_S}
\end{align} 

On the other hand, we write the average-based infinite-dimensional control law to compensate the diffusion process by 
\begin{align}
U(t) =\frac{c}{s+c}\left\lbrace K \left[G(t) + \hat{H}(t)\int_{0}^{D} (D-r)u(r,t)dr\right]\right\rbrace,
\label{eq:avg_pred}
\end{align}
where $c>0$ is sufficiently large.

\subsection{Extremum Seeking for Distinct Families of PDEs}

\begin{figure*}[ht]
	\centering
	\includegraphics[width=12cm]{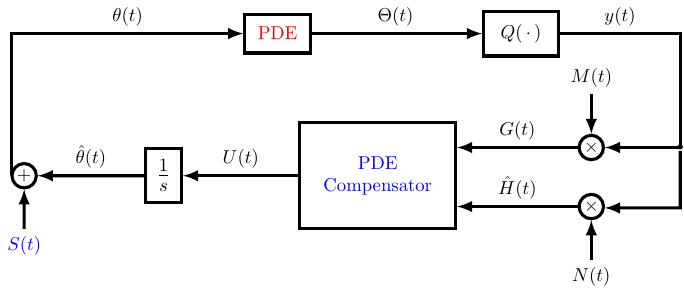}
	\caption{A general block diagram for implementation of ES control design for PDE connections (\textcolor{red}{in red}) at the input of nonlinear convex maps $Q(\cdot)$.
	Although the multiplicative perturbation signals $M(t)$ and $N(t)$ are the same as in the classical ES designs \cite{KW:00,GKN:2012}, the additive dither $S(t)$ (\textcolor{blue}{in blue}) must be redesigned using the trajectory generation paradigm \cite[Chapter~12]{krstic2008boundary} and the application of an adequate boundary control law (\textcolor{blue}{in blue}) for PDE compensation is necessary.}
	\label{ch12.fig:cascade}
\end{figure*}

In order to show that the proposed ES approach for infinite-dimensional systems is general and applicable to a wider class of PDEs,
we have formulated in Table~\ref{tabelademerda} the stabilizing boundary control (BC) law $U(t)$ and given the explicit solutions to the trajectory generation problem of $S(t)$  
for five other classes of distributed parameter systems \cite{krstic2008boundary}: (a) reaction-advection-diffusion (RAD)
PDEs \cite{Alcos:2019,IJACSP2020}, (b) wave equations \cite{OK:2019,ASME2021},
(c) hyperbolic transport PDEs ---for constant delays \cite{Tiago},
(d) time-varying delays \cite{ORDK:2018}, and (e) distributed delays \cite{TOK:2020,AUT2023}.


\begin{table*}[ht]
	\caption{ES for Distinct Classes of PDE Systems.}
	\centering
	\resizebox{15cm}{!}{
	\begin{tabular}{||c||l||}
		\hline
  \hline
		{\LARGE\textbf{RAD Equation}}    &\\                     
		              &
									{\LARGE$\textcolor{blue}{\text{PDE}:} ~\partial_t\alpha(x,t)\!=\!\epsilon\partial_{xx}\alpha(x,t)\!+\!
									b\partial_x\alpha(x,t)\!+\! \lambda\alpha(x,t)\,, \quad x\in[0,1]$} \\
		              &{\LARGE$\textcolor{blue}{\text{BC (Dirichlet):}} ~U(t)\!=\! \frac{c}{s+c}\left\lbrace k
									e^{-\frac{b}{2\epsilon}}\left[ \gamma(1)G(t)\!+\!\hat{H}(t)\int_{0}^{1}e^{\frac{b}{2\epsilon}\sigma} m(1\!-\!\sigma)
									u(\sigma,t)d\sigma \right]\right\rbrace,$}\\
									 &{\LARGE$\gamma(x) = \cosh\left(\sqrt{\frac{\xi}{\epsilon}}x\right)+\frac{b}{2\epsilon}\sqrt{\frac{\epsilon}{\xi}}									\sinh\left(\sqrt{\frac{\xi}{\epsilon}}x\right), \quad \xi:=b^2/(4\epsilon)-\lambda\ge 0\,, k>0$}\\
									&{\LARGE$m(x-\sigma) = \frac{1}{\epsilon}\sqrt{\frac{\epsilon}{\xi}}\sinh\left(\sqrt{\frac{\xi}{\epsilon}}
									(x-\sigma)\right), \quad \epsilon>0\,, b\geq 0\,, \lambda \geq 0$}\\
								  &{\LARGE$\textcolor{blue}{\text{Trajectory Generation}:} ~S(t) = e^{-\frac{b}{2\epsilon}}\sum_{k=0}^{\infty}
									\frac{a_{2k}(t)}{(2k)!}+\frac{b}{2\epsilon}\frac{a_{2k}(t)}{(2k+1)!},$}\\
									&{\LARGE$a_{2k} := \frac{a}{\epsilon^k}\sin(\omega t) \sum_{n=0}^{k}\binom{k}{2n}\xi^{k-2n}\omega^{2n}
									+ \frac{a}{\epsilon^k}\cos(\omega t)\sum_{n=0}^{k}\binom{k}{2n+1}\xi^{k-2n-1} \omega^{2n+1}$}\\
		&\\
		\hline
  \hline
		{\LARGE\textbf{Wave Dynamics}} &\\
									                                 &
									{\LARGE$\textcolor{blue}{\text{PDE}:} ~\partial_{tt}\alpha(x,t)=\partial_{xx}\alpha(x,t), \quad x\in[0,D]$} \\
		              &{\LARGE$\textcolor{blue}{\text{BC (Neumann):}} ~U(t)=\frac{c}{s+c}
									\left\{c\left[k \hat{H}(t)u(D,t) - \partial_t u(D,t)\right] + \rho(D) G(t) + 
									\phantom{\int_{0}^{D}} \right.$} \\
									&{\LARGE$\left.  \hat{H}(t)\int_{0}^{D} \rho(D-\sigma) \partial_t u(\sigma,t) d\sigma\right\}$, $\quad 
									\rho(s)=k[0 \ \ I]e^{As}[0 \ \ I]^T, \quad A=\begin{pmatrix}
		                                    0 & 0 \\ I & 0 \\
	                                      \end{pmatrix}$}\\
								  &{\LARGE$\textcolor{blue}{\text{Trajectory Generation}:} ~S(t) = a
									\cos(\omega D) \sin(\omega t)$}\\
		&\\
		\hline
  \hline 
		{\LARGE\textbf{Constant Delay}}     &\\       &
									{\LARGE$\textcolor{blue}{\text{PDE}:} ~\partial_{t}\alpha(x,t)=\partial_{x}\alpha(x,t), \quad x\in[0,D]$} \\
		              &{\LARGE$\textcolor{blue}{\text{BC (Dirichlet)}:} ~\textcolor{black}{U(t)=\frac{c}{s+c}\left\{k 
									\left[G(t)+\hat{H}(t)\int_{0}^{D} u(\sigma,t) d\sigma \right]\right\}}$}\\
								  &{\LARGE$\textcolor{blue}{\text{Trajectory Generation}:} ~\textcolor{black}{S(t) = a\sin(\omega (t+D))}$}\\
		&\\
		\hline
  \hline
		{\LARGE\textbf{Variable Delay}}      &\\       &
									{\LARGE$\textcolor{blue}{\text{PDE}:} ~\partial_{t}\alpha(x,t)= \pi(x,t)\partial_{x}\alpha(x,t),
									\quad x\in[0,1], \quad \pi(x,t)=\frac{1+x[\frac{d(\phi^{-1}(t))}{dt}-1]}{\phi^{-1}(t)-t}$} \\
		              &{\LARGE$\textcolor{blue}{\text{BC (Dirichlet):}} ~\textcolor{black}{U(t)=\frac{c}{s+c}\left\{k 
									\left[G(t)+\hat{H}(t)\int_{0}^{1} u(\sigma,t)\left(\phi^{-1}(t)-t\right) d\sigma \right]\right\}}$}\\
								  &{\LARGE$\textcolor{blue}{\text{Trajectory Generation}:} 
									~\textcolor{black}{S(t) = a\sin(\omega t)},\quad \phi(t):= t-D(t)$}\\ 
									&{\LARGE$\textcolor{blue}{\text{Demondulation}:} ~M(t) = \frac{2}{a}\sin(\omega (t-D(t)))\,, ~N(t) = -\frac{8}{a^2}\cos(2\omega (t-D(t)))$}\\
									
		&\\
		\hline
  \hline
		{\LARGE\textbf{Distributed Delay}}     &\\        &
									{\LARGE$\textcolor{blue}{\text{PDE}:}  ~\partial_{t}\alpha(x,t)=\partial_{x}\alpha(x,t), \quad x\in[0,D], 
									\quad y=Q\left( \int_{0}^{D} \Theta(t-\sigma) d\beta(\sigma)  \right)$} \\
		              &{\LARGE$\textcolor{blue}{\text{BC (Dirichlet):}} ~\textcolor{black}{U(t)=\frac{c}{s+c}\left\{k 
									\left[G(t)+\hat{H}(t)\int_{0}^{D} (1-\beta(\sigma))u(D-\sigma,t) d\sigma \right]\right\}}$}\\
								  &{\LARGE$\textcolor{blue}{\text{Trajectory Generation}:} 
									~\textcolor{black}{S(t) = \frac{a}{\gamma(\omega)}\int_{0}^{D}\sin(\omega (t+\xi))d\beta(\xi)}$}\\							
		&\\
		\hline
  \hline
	\end{tabular}
	}
\label{tabelademerda}
\end{table*}


As for the \textcolor{black}{transport and diffusion cases} explored in the NES sections, the term $u(x,t)$ which appears in $U(t)$ of Table~\ref{tabelademerda}, is the state of the infinite-dimensional system corresponding to a copy of the PDE model of  the actuator dynamics. 

In \cite{PDE_cascades_SCL}, differently from what has been done which has dealt with PDEs at the input into an unknown map, and in which we have already advanced from transport PDEs to reaction-advection-diffusion PDEs to wave PDEs according to Table~~\ref{tabelademerda}, the authors consider one last configuration in which the input pathway to the map contains a cascade of PDEs from distinct classes. 
There, we deal with PDEs with input delays such as, for example, the notorious problem of a wave PDE with an input delay where, if the delay is left uncompensated, an arbitrarily short delay destroys the closed-loop stability. Then, we move forward to an even more challenging class of problems for parabolic-hyperbolic cascades of PDEs, focusing on a heat equation at the input of a wave PDE.  
The treatment of such systems with PDE-PDE cascades is again performed by means of boundary control. 
Local exponential stability and convergence to a small neighborhood of the unknown extremum point are  guaranteed by using a backstepping transformation and averaging in infinite dimensions. 

PDE-PDE cascades have a great deal in common with PDE-ODE cascades. For instance, a cascade of a delay into a PDE is just a much generalized version of an integrator with an input delay. 
However, while in a delay-integrator cascade the design can be pursued within the predictor feedback framework, with backstepping just employed for an interpretation and for analysis, in delay-PDE cascades a predictor for a PDE is too complicated of a mathematical object to be of value. Instead, the design is pursued entirely by the backstepping approach. Similarly, while the heat-integrator cascade became familiar, 
and was dealt with through the backstepping design, a more general backstepping design is applied to a heat-wave PDE cascade in a part of reference \cite{PDE_cascades_SCL}.

\begin{figure*}[ht]
\centering
\includegraphics[width=12cm]{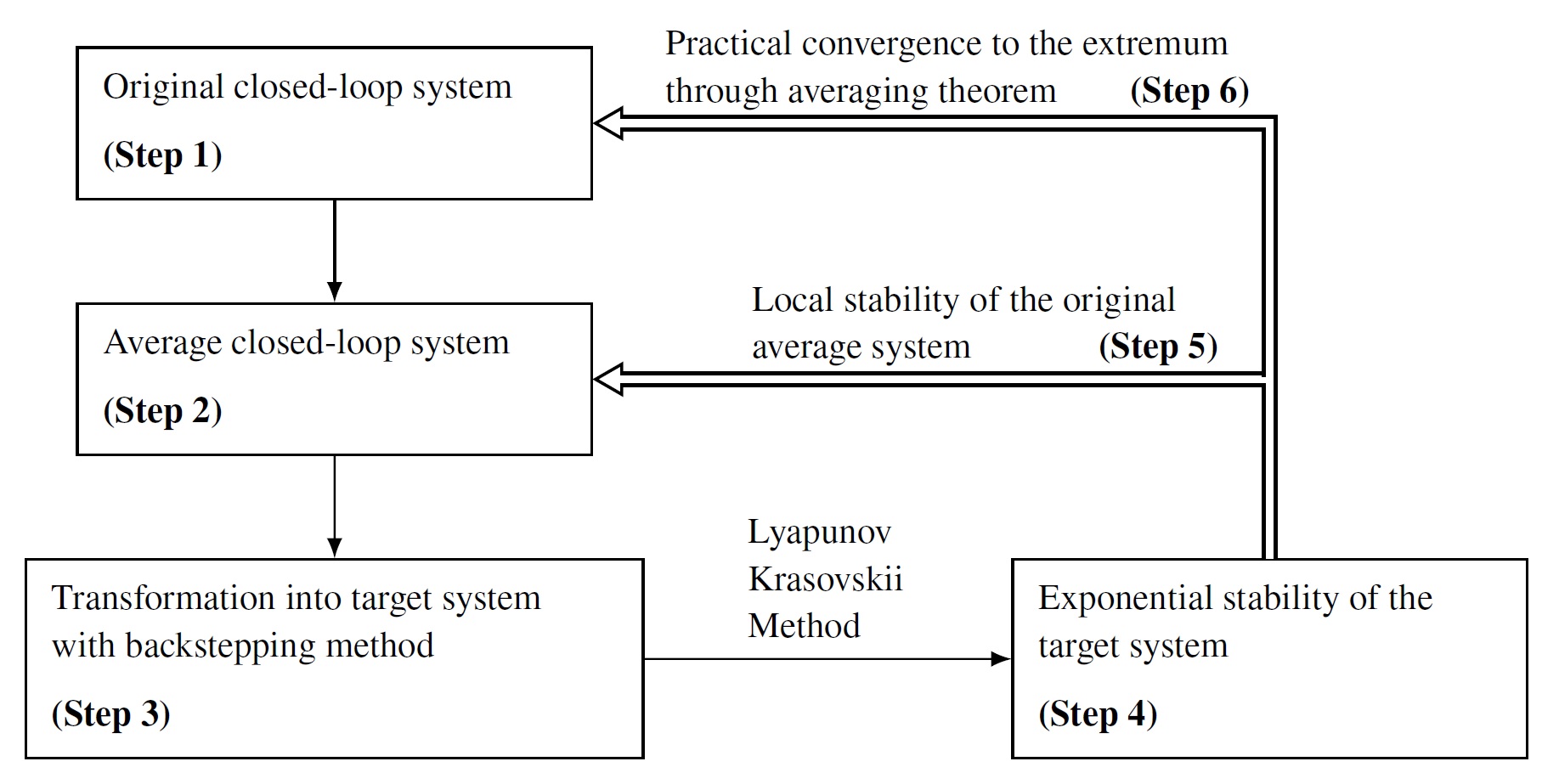}
\caption{\textcolor{black}{Structure of the stability proof for the closed-loop system.}} \label{ch9.fig:structure_proof}
\end{figure*}

Although our goal in this paper is to avoid numbing the reader with lengthy proofs, it is worth providing a general picture on how we carry out the steps to prove the stability results of the ES feedback loop in the presence of PDEs of Fig.~\ref{ch12.fig:cascade}. Fig.~\ref{ch9.fig:structure_proof} shows the structure of the proof, divided into six main steps. We seize this opportunity to highlight that our analysis presents a carefully constructed sequence of analytical steps, a predictor-based infinite-dimensional backstepping transformation, a synthesist of a Lyapunov functional (rather than small-gain analysis), and computation of a Lyapunov estimate, for the overall infinite-dimensional system with nonlinearities, stochastic perturbations, and distributed delays. The analysis process involving so many steps has a large number of possible permutations---all of which \textit{but one} would be wrong. We show how to properly sequence the steps of averaging, backstepping, and Lyapunov functional analysis, to prove stability. This ``analysis pathway'' will serve the needs of future researchers who deal with stochastic extremum seeking under delays. The complete details can be found in the authors' book \cite{OK:2022}. 
%

\section{Applications}

For the distinct cases involving ES plus PDEs in Table~~\ref{tabelademerda}, we dedicated the remainder of the paper to some select applications:

\textbf{1)} \textit{Traffic Control} for linearized Lighthill-Whitham-Richards (LWR) macroscopic PDE models transformed into constant delays \cite{YKOK:2021,Tiago}; 
\textbf{2)} \textit{Optimal Oil Drilling Control} with ES for wave models \cite{AK:2019}; 
\textbf{3)} \textit{Deep-Sea Cable-Actuated Source Seeking} modeled by wave PDEs with Kelvin-Voigt damping \cite{L-CSS:2023}; 
\textbf{4)} \textit{Additive manufacturing} modeled by the \textit{Stefan PDE} \cite{KDK:2019,TDS2022}; 
\textbf{5)} \textit{Bioreactors} considering ES for models described by parabolic PDEs (reaction-diffusion equations), see \cite{DOCHAIN:2006,DOCHAIN:2008};  
\textbf{6)} \textit{Light-Source Seeking} with infinite-dimensional models represented by Euler-Bernoulli beams equations \cite[Chapter~8]{krstic2008boundary}; and  
\textbf{7)} \textit{Neuromuscular Electrical Stimulation} (NMES) problem for ES under time-varying delays \cite{PAZ:2019,ORDK:2018}. 


\section{Traffic Congestion Control with a Downstream Bottleneck}

This section develops boundary control for freeway traffic with a downstream bottleneck \cite{YKOK:2021}. Traffic on a freeway segment with capacity drop at outlet of the segment is a common phenomenon that leads to traffic bottleneck problem. The capacity drop can be caused by lane-drop, hills, tunnel, bridge or curvature on the road. If incoming traffic flow remains unchanged, traffic congestion forms upstream of the bottleneck because the upstream traffic demand exceeds its capacity. Therefore, it is important to regulate the incoming traffic flow of the segment to avoid overloading the bottleneck area. Traffic densities on the freeway segment are described with the Lighthill-Whitham-Richards (LWR) macroscopic Partial Differential Equation (PDE) model. To mitigate the traffic congestion upstream of the bottleneck, incoming flow at the inlet of the freeway segment is controlled so that the optimal density that maximizes the outgoing flow is reached. The density and traffic flow relation at the bottleneck area, described with the fundamental diagram, is considered to be unknown. We tackle this problem using Extremum Seeking (ES) Control with delay compensation for LWR PDE \cite{YKOK:2021}. ES control, a non-model based approach for real-time optimization, is adopted to find the optimal density for the unknown fundamental diagram. A predictor feedback control design is proposed to compensate the delay effect of traffic dynamics in the freeway segment. 

\subsection{Problem Statement}

\begin{figure}[ht]
	\centering
	\includegraphics[width=0.4\textwidth]{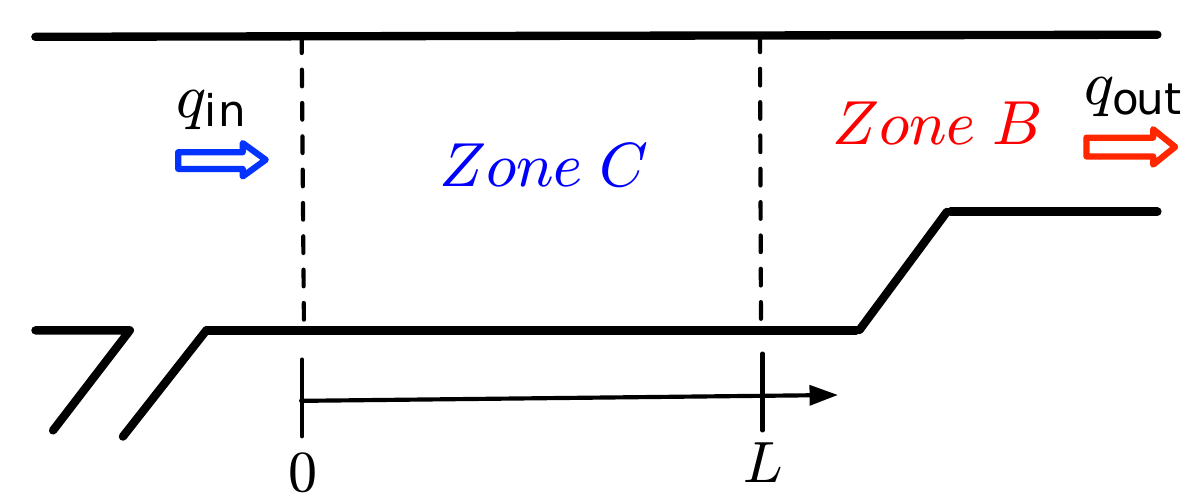}
	\caption{Traffic on a freeway segment with lane-drop.} \label{bottleneck}
\end{figure}

We consider a traffic congestion problem on a freeway-segment with lane drop bottleneck downstream of the segment. The freeway segment upstream of the bottleneck and the lane-drop area are shown in Fig.~\ref{bottleneck} which illustrates the clear ``Zone C'' and the bottleneck ``Zone B'', respectively. The flow is conserved through the clear Zone C to the bottleneck Zone B. The local road capacity is changed due to the lane-drop in Zone B which could be caused by working zone, accidents or lane closure. To prevent the traffic in Zone B overflowing its capacity and then causing congestion in the freeway segment, we aim to find out the optimal density ahead of Zone C that maximizes outgoing flux of Zone B given unknown density-flow relation. Traffic dynamics in Zone C are described with the macroscopic LWR traffic model for the aggregated values of traffic density.

Due to the reduction of lanes in Zone B, the fundamental diagram for the flow and density relation usually changes, which leads to a capacity drop in Zone B. The control objective is to find the optimal input density at inlet of Zone C that drives the measurable output flux of Zone B to its unknown optimal value of an unknown fundamental diagram. 

	\subsection{LWR Traffic Model} 	
	The traffic dynamics in Zone C, upstream of Zone B is described with the first-order, hyperbolic LWR model. 
	 Traffic density $\rho(x,t)$ in Zone C is governed by the following nonlinear hyperbolic PDE, where $x\in[0,L]$, $t\in[0,\infty)$, 
\begin{align}\label{ch*.sys1} 
\partial_t \rho + \partial_x( Q_C(\rho))=&0. 
\end{align}
The fundamental diagram of traffic flow and density function $Q_C(\rho)$ is given by 
$Q_C(\rho) = \rho V(\rho),$
where traffic velocity follows an equilibrium velocity-density relation $V(\rho)$. There are different models that describe the flux and density relationship. A basic and popular choice is Greenshield's model for $V(\rho)$ which is given by 
$
V(\rho)=v_f\left(1-\frac{\rho}{\rho_m}\right), 
$
where $v_f \in \mathbb{R^+}$ is defined as the maximum velocity and $\rho_m \in \mathbb{R^+}$ is the maximum density for Zone C \cite{YKOK:2021}. Then the fundamental diagram of flow and density function $Q_C(\rho)$ is in a quadratic form of density,
\begin{align} 
Q_C(\rho) = -\frac{v_f}{\rho_m}\rho^2 + v_f \rho. \label{ch*.quad}
\end{align}
A critical value of density segregates the traffic into the free flow regime whose density is smaller than the critical value and the congested regime whose density is greater than the critical value. The critical density may be assumed as $\rho_c = \rho_m/2$ for \eqref{ch*.quad} \cite{YKOK:2021}. For the fundamental diagram calibrated with the freeway empirical data, the critical density usually appears at $20\% $ of the maximum value of the density \cite{Dervisoglu},~\cite{Fan}. 

In practice, the quadratic fundamental diagram sometimes does not fit well with traffic density-flow field data. There are several other equilibrium models, \textit{e.g.}, Greenberg model, Underwood model and diffusion model for which the fundamental diagrams are nonlinear functions, see \cite{YKOK:2021} and the references therein. However, by Taylor expansion, a second-order differentiable nonlinear function can be approximated by a quadratic function in the neighborhood of its extremum. The following assumption is made for the nonlinear fundamental diagram. The stability results derived here hold locally for the general form of the fundamental diagram ${\mathcal Q}(\rho)$ that satisfies the following assumption. Here we can adopt other density-flow relations for the fundamental diagram ${\mathcal Q}(\rho)$ but requiring Assumption~\ref{ch*.ass:diagram} below to be satisfied.

\begin{assum} \label{ch*.ass:diagram} 
	The fundamental diagram ${\mathcal Q}(\rho)$ is a smooth function, and it holds that 
$
{\mathcal Q}'(\rho_c) = 0, ~
{\mathcal Q}''(\rho_c) <0. 
$
\end{assum}

Under Assumption \ref{ch*.ass:diagram}, the fundamental diagram can be approximated around the critical density $\rho_{c}$ as follows:
$
		{\mathcal Q}(\rho) = q_c + \frac{{\mathcal Q}^{\prime\prime}(\rho)}{2}(\rho(t)-\rho_c)^2,
$
	where $q_c = {\mathcal Q}(\rho_c)$ is defined as the road capacity or maximum flow, with ${\mathcal Q}^{\prime\prime}(\rho)< 0$.

\subsection{Lane-Drop Bottleneck Control Problem}

\begin{figure}[ht] 
	\centering
	\includegraphics[scale=0.35]{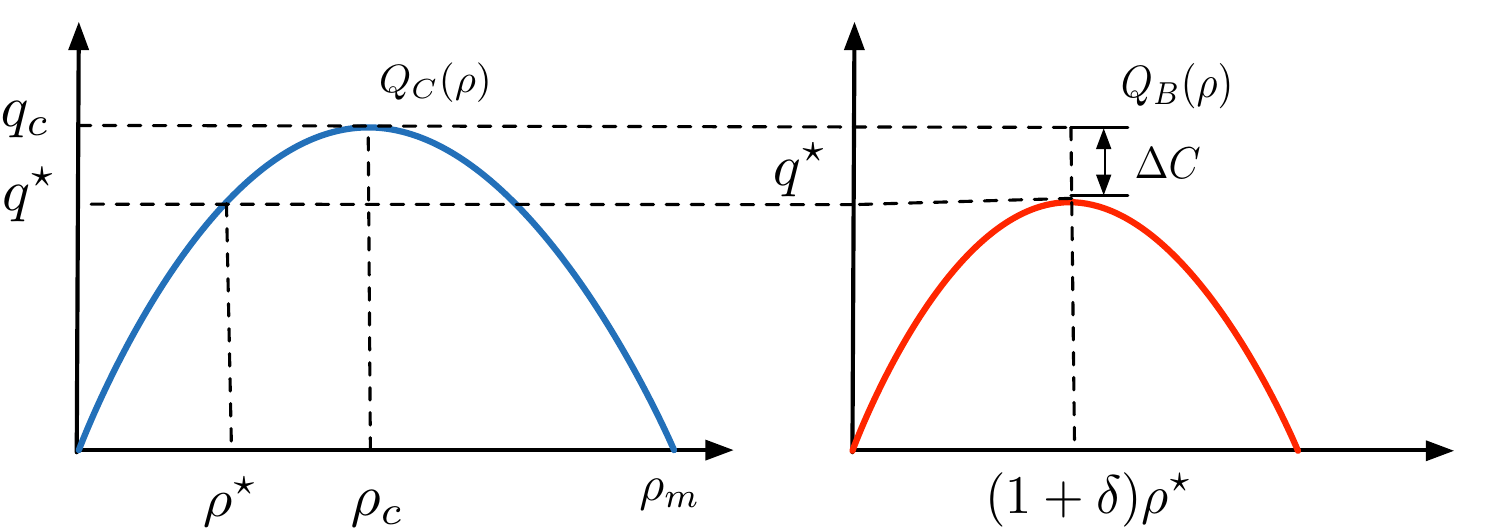}
	\caption{Quadratic fundamental diagram for the clear Zone C and the bottleneck Zone B.}
	\label{ch*.model}
\end{figure}

Due to the reduction of the number of lanes from Zone C to Zone B, we consider the equilibrium density-flow relation of Zone B as shown in Fig.~\ref{ch*.model}, as pointed out in \cite{Treiber}. There is a capacity drop $\Delta C$ of $Q_B$ in Zone B compared to $Q_C$ in Zone C after the congestion has formed upstream of the lane-drop area. The capacity drop caused by a sudden lane-drop is hard to measure in real time and the traffic dynamics of Zone B are affected by the lane-changing and merging activities. Therefore we assume that the fundamental diagram $Q_B(\rho)$ of Zone B is unknown. 

In Fig.~\ref{ch*.model}, the capacity is 
\begin{align}
	\Delta C &= Q_C(\rho_c) - Q_B((1+\delta)\rho^\star),  \\
	q^\star &= Q_C(\rho^\star)= Q_B((1+\delta)\rho^\star),
\end{align}
where $\Delta C$ is unknown. The $\rho^{\star}\in \mathbb{R^+}$ represents the optimal density that keeps Zone C in the free regime while  $(1+\delta)\rho^\star$ reaches the critical density of Zone B so that the discharging flow rate reaches its maximum value $q^\star \in \mathbb{R^+}$. The ratio $\delta$ accounts for the density discontinuity before the outlet in Zone C and after the outlet in Zone B. We assume that $\Delta C$ and $\delta$ are unknown and therefore the optimal density and flow rate $(\rho^\star, q^\star)$ are unknown.

When there is a lane-drop bottleneck presenting downstream, the density at the outlet of Zone C is $\rho(L,t)$, governed by the PDE
in \eqref{ch*.sys1} for $x\in[0,L]$, $t\in[0,\infty)$. The inlet boundary flow is 
$
q_{\textrm{in}}(t) = Q_C(\rho(0,t)).
$
The output measurement of traffic flow in Zone B, $q_{\textrm{out}}(t)$ is given by $Q(\rho)$ with outlet density $\rho(L,t)$,
$
q_{\textrm{out}}(t) = Q(\rho(L,t) ).
$
where the function  $Q(\rho)$ of outlet boundary $x=L$ connecting Zone C and Zone B is defined as follows
\begin{small}
\begin{align}
Q(\rho(L,t))=\begin{cases}
Q_C(\rho(L,t)), \quad\quad\quad\quad\quad\quad\quad\quad ~~\rho(L,t) < \rho^\star,\\
Q_C(\rho^\star) = q^\star = Q_B((1+\delta)\rho^\star), \!\!\!\quad \rho(L,t) = \rho^\star,\\
Q_B((1+\delta)\rho(L,t)),\quad\quad\quad\quad\quad\quad \!\!\!\!\!\!~\rho(L,t) > \rho^\star,
\end{cases}
\end{align}
\end{small}
so that the flow is conserved through the boundary, entering from Zone C to Zone B.
Note that when the optimal density $\rho^\star$ is reached, the flow rate at the outlet of Zone C and the input of Zone C reaches the equilibrium and its maximum value $q^\star$.

The control objective is to design the traffic flow input $q_{\textrm{in}}(t)$ so that the outgoing flow in the lane-drop area Zone B $q_{\textrm{out}}(t)$ is maximized. We aim to find out the optimal outlet density $\rho(L,t)=\rho^{\star}$ that maximizes $q_{\textrm{out}}(t)$ of Zone B and then using the PDE that describes the dynamics of traffic in Zone C to obtain the desirable flow input $q_{\textrm{in}}(t)$ from the inlet of Zone C. Here we approximate $q_{\textrm{out}}(t)$ with a function that satisfies Assumption~\ref{ch*.ass:diagram} and $q_{\textrm{out}}(t)$ can be written as
\begin{align}
	q_{\textrm{out}}(t) = q^\star + \frac{H}{2}(\rho(L,t)-\rho^\star)^2, \label{ch*.static}
\end{align}	
where $H<0$ is the unknown Hessian of the  approximated static map $q_{\textrm{out}}(t)$. 


{Note that we use a static fundamental diagram to model the traffic in the bottleneck Zone B. Therefore, the upstream propagating traffic waves from Zone B to Zone C cannot be captured by our model if Zone B is very congested. Since this result is focused on maximizing the discharging flow rate at the bottleneck area, the ES control seeks the optimal traffic density value in its neighborhood. In bottleneck Zone B, the closer the outlet traffic density $\rho(L,t)$ is to the optimal value $\rho^\star$ where $Q'(\rho) = 0$ is satisfied, the smaller is the propagating characteristic speed of the traffic waves $Q'(\rho)$. Therefore, the spill-back traffic from Zone B to Zone C is negligible in our model.}

In order to find the unknown optimal density at the bottleneck area, we design ES control for the unknown static map $Q(\rho) $ with actuation dynamics governed by a nonlinear hyperbolic PDE in \eqref{ch*.sys1}. Below, we linearize the nonlinear PDE and the traffic dynamics can be represented by the delay effect for the control input design.

\subsection{Linearized Reference Error System} 
We linearize the nonlinear LWR model around a constant reference density $\rho_{r} \in \mathbb{R^+}$, which is assumed to be close to the optimal density $\rho^{\star}$.
Note that the reference density $\rho_{r}$ is in the free regime of $Q(\rho)$ of Zone C, and thus is smaller than the critical density $\rho_{c}$, and therefore the following is satisfied
$
\rho_{r} <\rho_c.
$
Define the reference error density as
\begin{align}
\tilde{\rho}(x,t) = \rho(x,t) - \rho_{r},
\end{align}
and reference flux $q_r$ to be 
$
q_{r} = Q(\rho_r) > 0.
$
By the governing equation \eqref{ch*.sys1} together with \eqref{ch*.quad}, the linearized reference error model is derived as 
\begin{align}
\partial_{t} \tilde{\rho}(x,t)   + u \partial_{x} \tilde{\rho}(x,t)  =& 0,  \label{ch*.trho}\\
\tilde{\rho}(0,t) =& \rho(0,t)- \rho_{r}, \label{ch*.in}
\end{align} 
where the constant transport speed $u$ is given by
$ u = Q^\prime ( \rho)|_{\rho =\rho_r } 
=  V(\rho_r) + \rho_r V^\prime ( \rho)|_{\rho =\rho_r }.
$
The equilibrium velocity-density relation $V(\rho)$ is a strictly decreasing function. 
The reference density $\rho_{r}$ is in the left-half plane of the fundamental diagram $Q_c(\rho)$, which yields the following inequality for the propagation speed $u>0.
$
We define the input density as
$
\varrho(t)  = {\rho}(0,t),
$
and the linearized input at inlet to be  
\begin{align}
\tilde \varrho(t)  =& \varrho(t) - \rho_{r}.
\end{align}
The linearized error dynamics in \eqref{ch*.trho}, \eqref{ch*.in} is a transport PDE with an explicit solution for $t > \frac{x}{u}$ and thus is represented with input density
$
\tilde{\rho}(x,t) = \tilde\varrho\left( t - \frac{x}{u} \right).
$
The density variation at the outlet is
\begin{align}
\tilde{\rho}(L,t) = &\tilde \varrho\left( t - D \right).\label{ch*.orho1}
\end{align}
where the time delay is $D= \frac{L}{u}$.
Therefore, the density at the outlet is given by a delayed input density variation and the reference: 
\begin{align} 
{\rho}(L,t) = & \rho_r + \tilde {\rho}(L,t).\label{ch*.orho2}
\end{align} 
Finally, substituting \eqref{ch*.orho1}, \eqref{ch*.orho2} into the static map \eqref{ch*.static}, we arrive at the following: 
\begin{align}
q_{\textrm{out}}(t) =&  q^\star + \frac{H}{2}( \tilde\varrho\left( t - D \right) + \rho_{r} -\rho^{\star} )^2 \notag\\
= &  q^\star + \frac{H}{2} \left( \varrho\left( t - D \right)  - \rho^{\star} \right)^2.  \label{ch*.out}
\end{align}	
The control objective is to regulate the input $q_{\textrm{in}}(t)$ so that $\varrho\left( t - D \right)$ reaches to an unknown optimal $\rho^\star$ and the maximum of the uncertain quadratic flux-density map $q_{\textrm{out}}(t)$ can be achieved. We can apply the method of ES for static map with delays, originally developed in~\cite{Tiago}. The ES control is designed for finding the extremum of the unknown map.

In practice, control of density at the inlet can be realized with a coordinated operation of a ramp metering and a variable speed limit (VSL) at the inlet, which is widely used in freeway traffic management~\cite{Carlson,Hegyi,Lu,YHZ,Karafyllis,Huan,Yu20}. The controlled density at the inlet is implemented by 
$
\varrho(t) = \frac{q_{\textrm{in}}(t)}{v_c}.
$
where $v_c$ is the speed limit implemented by VSL and $q_{\textrm{in}}(t)$ is actuated by an on-ramp metering upstream of the inlet. Note that the linearized model is valid at the optimal density $\rho^{\star}$ since the reference density is assumed to be chosen near the optimal value.

\subsection{Online Optimization by Extremum Seeking Control}

\begin{figure*}[ht]
	\centering
	\includegraphics[width=0.9\textwidth]{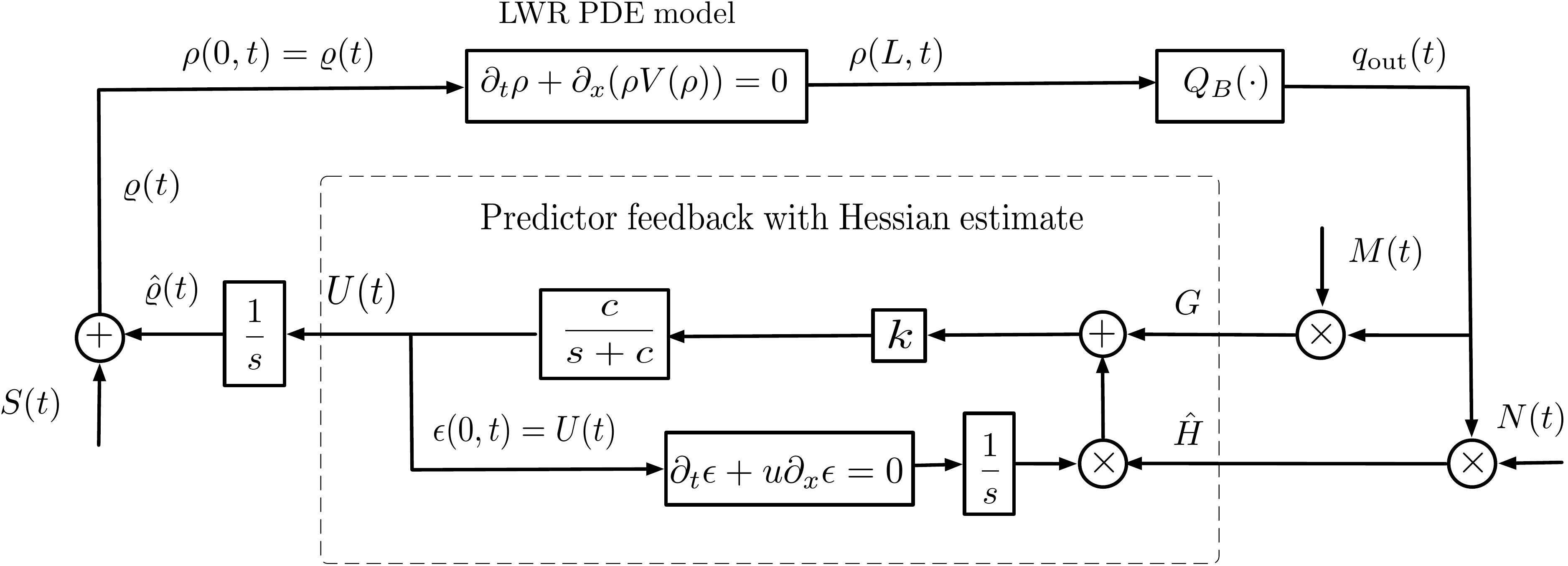}
	\caption{Block diagram for implementation of ES control design for the nonlinear LWR PDE model.}
	\label{ch*.fig:ES}
\end{figure*}

Here, we present the design of ES control with delay by following the procedure in \cite{Tiago}. The block diagram of the delay-compensated ES algorithm applied to LWR PDE model is depicted in Fig.~\ref{ch*.fig:ES}. 

Let $\hat{\varrho}(t)$ be the estimate  of $\rho^\star$, and $e(t)$ be the estimation error defined as
\begin{align}
e(t) = \hat \varrho(t) - \rho^\star, \label{ch*.e}
\end{align}
where $ \hat \varrho(t)$ is an integrator of the predictor-based feedback signal $U(t)$ as 
$
	\dot{\hat{\varrho} }(t) = U(t).
$
From Fig.~\ref{ch*.fig:ES}, the error dynamics can be written as
\begin{align}
\dot e(t-D) = U(t-D), \label{ch*.doterror}
\end{align}
given the delayed estimation error dynamics modeled by $ \epsilon(x,t) = U(t-\frac{x}{u})$.


We introduce the additive perturbation $S(t)$  
\begin{align} \label{ch*.s}
	S(t) = a \sin \left( \omega (t+D) \right)\,,
\end{align}
and the multiplicative demodulation signals $(M(t), N(t))$ given by
\begin{align} 
M(t) = \frac{2}{a} \sin\left( \omega t \right) ,  
 \quad N(t) = - \frac{8}{a^2} \cos\left( 2 \omega t \right),\label{ch*.Nt}
\end{align} 
where $a$ and $\omega$ are respectively the amplitude and frequency of a slow periodic perturbation signal $a \sin (\omega t)$ introduced later. 
Using the demodulation signals, we calculate estimates of the gradient and Hessian of the cost function, denoted as $(G(t), \hat{H}(t))$,
\begin{align} 
G(t) = M(t) q_{\textrm{out}}(t),  \quad \hat{H}(t) = N(t) q_{\textrm{out}}(t),\label{ch*.Hhat} 
\end{align} 
where $\hat{H}(t)$ is the estimate of the unknown Hessian $H$. The averaging of $G(t)$ and $\hat{H}(t)$ yields that 
\begin{align}
G_{\rm av}(t) =H e_{\rm av}(t-D), \quad
\hat H_{\rm av} = (Nq_{\textrm{out}})_{\rm av} = H. \label{ch*.HavH}
\end{align}
Taking the average of \eqref{ch*.doterror}, we have
$\dot e_{\rm av}(t-D) = U_{\rm av}(t-D),$
where $ U_{\rm av}(t)$ is the averaged value for $U(t)$ to be designed later. Substituting the above equation into \eqref{ch*.HavH} yields 
\begin{align}
\dot G_{\rm av}(t) = HU_{\rm av}(t-D). \label{ch*.tD}
\end{align}
The motivation for predictor feedback design is to compensate for the delay by feeding back future states in the equivalent averaged system $G_{av}(t+D)$. 
Given an arbitrary control gain $k>0$, we aim to design
\begin{align}
U_{\rm av}(t) = k G_{\rm av}(t + D), \quad \forall t\geq 0, \label{ch*.tD1}
\end{align}
which requires knowledge of future states.
Therefore we have the following by plugging \eqref{ch*.tD1} into \eqref{ch*.doterror},
\begin{align}
\dot e_{\rm av}(t) =U_{\rm av}(t) = kH e_{\rm av}(t), \quad \forall t \geq D.
\end{align} 
Recalling that $k>0, H<0$, the equilibrium of the average system $e_{\rm av}(t) = 0$ is exponentially stable. Applying the variation of constants formula, $G_{\rm{av}}(t\!+\!D)\!=\!G_{\rm{av}}(t)\!+\!\hat{H}_{\rm{av}}(t) \int_{t-D}^{t}U_{\rm{av}}(\tau) d\tau$ and, from \eqref{ch*.tD1}, one has:
\begin{align}
U_{\rm av}(t) =  k \left( G_{\rm av}(t) + \hat {H}_{\rm av}(t) \int_{t-D}^{t} U_{\rm av} (\tau) d \tau \right),  \label{ch*.Uavt}
\end{align}
which represents the future state $G_{\rm av}(t + D)$ in \eqref{ch*.tD} in terms of the average control signal $U_{\rm av}(\tau)$ for $\tau \in [t-D,t]$. The control input is infinite-dimensional due to its use of history over the past $D$ time units.

For the stability analysis in which the averaging theorem for infinite-dimensional systems is used, we employ a low-pass filter for the above basic predictor feedback controller and then derive an infinite-dimensional and averaging-based predictor feedback given by
\begin{align} 
U(t) = \mathcal{T} \left\{ k \left( G(t) + \hat{H}(t) \int_{t-D}^{t} U (\tau) d \tau \right) \right\} ,  \label{ch*.Ut}
\end{align} 
where $k>0$ is an arbitrary control gain, the Hessian estimate $\hat H(t)$ is updated according to \eqref{ch*.Hhat}, satisfying the average property in \eqref{ch*.HavH}.  $\mathcal{T} \{\cdot\}$ is the low-pass filter operator defined by 
\begin{align} 
\mathcal{T} \left\{ \varphi (t) \right\} = \mathcal{L}^{-1} \left\{ \frac{c}{s+c} \right\} * \varphi(t), 
\end{align} 
where $c \in \mathbb{R^+}$ is the corner frequency, $\mathcal{L}^{-1} $ is the inverse Laplace transformation, and $*$ steads for convolution in time.

Hence, according to \cite[Theorem~1]{YKOK:2021}, we can conclude that 
\begin{align} 
	\lim\limits_{t \to +\infty} \sup|\varrho(t) - \rho^\star| =& \mathcal{O}(a+1/\omega),\\
	\lim\limits_{t \to +\infty} \sup|q_{\rm out}(t) - q^\star| =& \mathcal{O}(a^2+1/\omega^2). \label{ch*.limsup2}
\end{align}

\begin{figure}[ht]
	\centering
	\includegraphics[width=8.1cm]{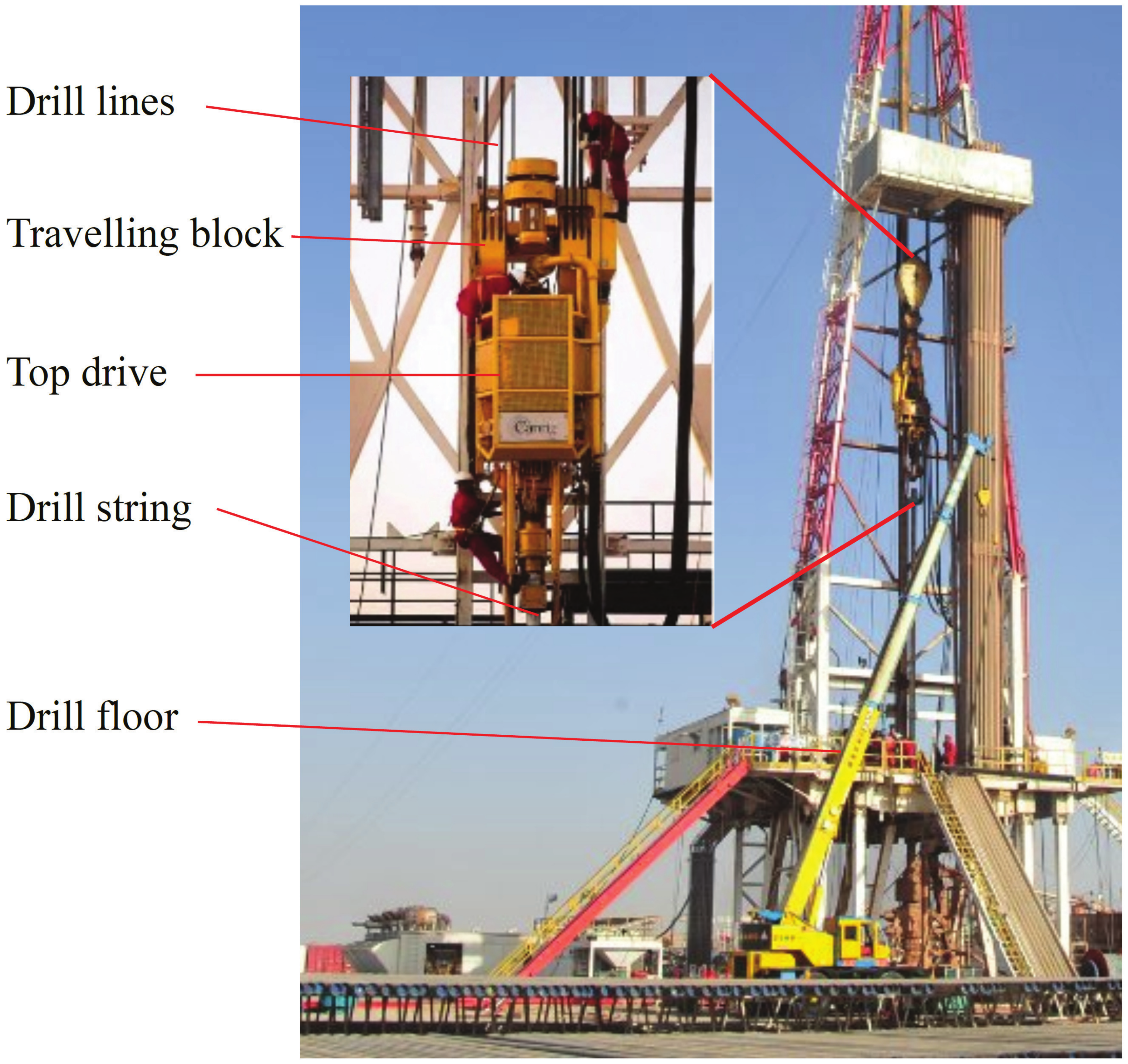}
	\caption{Picture showing the topside of a drilling rig.}
	\label{ch10.fig:petroleo}
\end{figure}
\begin{figure*}[ht]
	\centering
  \includegraphics[width=13cm]{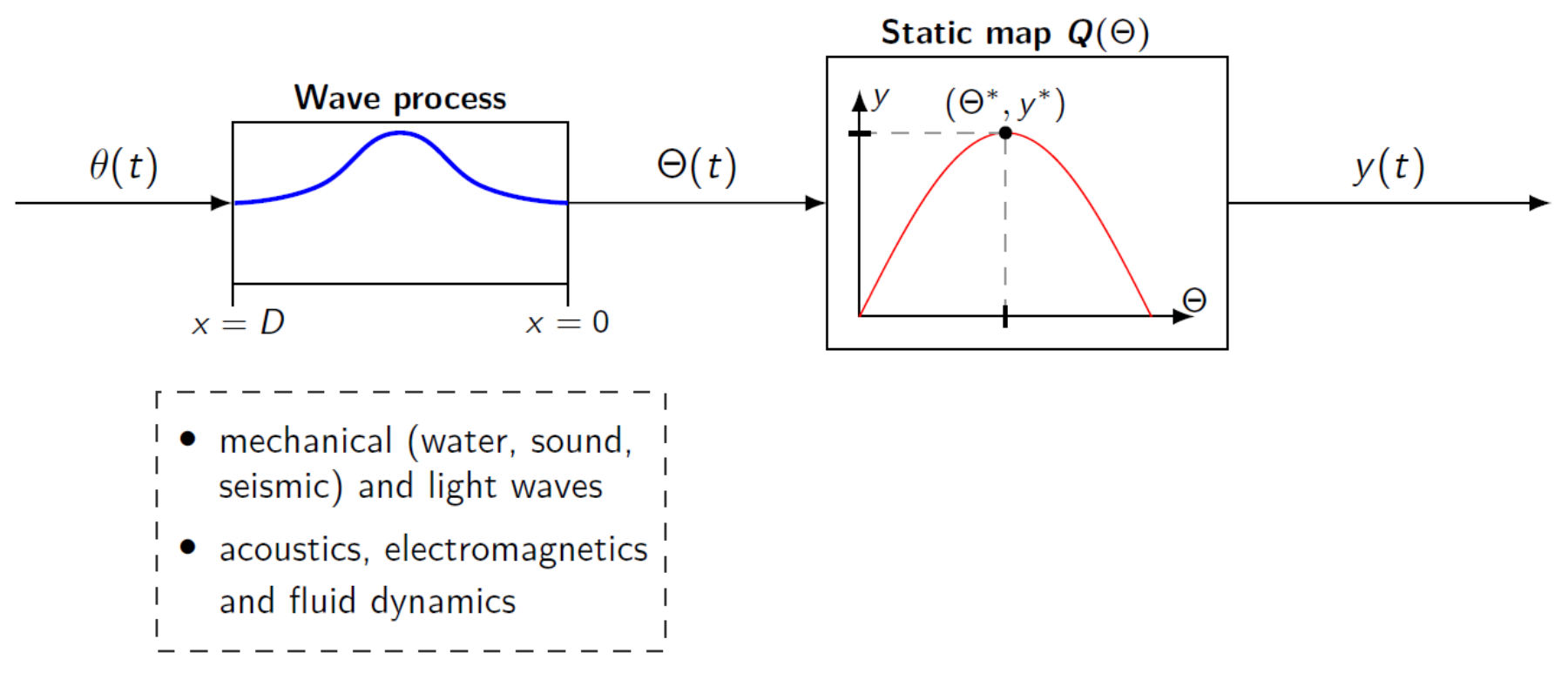}	
	\caption{Cascade of wave PDE with a static map $y(\Theta(t))=Q(\Theta(t))$. The extremum $y(t)=y^*$ is achieved for $\Theta(t)=\Theta^*$. Wave PDEs are used generally to model different sort of processes such as mechanical, acoustics, electromagnetic and fluid dynamics.}
	\label{ch10.fig:cascade}
\end{figure*}

\section{Optimal Oil Drilling Control}


A common type of instability in oil drilling is the friction-induced stick-slip oscillation (see \cite{BK:2014} and references therein), which results in torsional vibrations of the drill-string and can severely damage the drilling facilities
(see Fig.~2 from \cite{AK:2019}). 

The picture in Fig.~\ref{ch10.fig:petroleo} shows a modern land-based drilling rig. The tower operates like the derrick of a crane: the traveling block is connected by several steel drill lines with one attached to the deadline anchor and the other being spooled on a drum controlled by AC induction motors. Another electric motor, called Top Drive, is connected to the travelling block. The Top Drive is used to rotate the drill-string, a set of hundreds of drill pipes (about $30$ feet long each) that conducts the Bore Hole Assembly (BHA). The BHA contains several sensors (pressure, temperature, vibration among others) and the drill bit itself. There are several different types of drill bit design and materials, adequate for drilling different geological formations.

In analogy, the rig operates similarly to a drill press, but with drill bit which is several inches wide ($4$" to $36$" is a common range) and up to several miles long (an onshore well can be as shallow as $200$ yards or as deep as $2$ miles). By rotating this drill-string and using its weight to generate an axial force, the BHA mills the rocks, drilling the well. Because of the small diameter when compared to its length, the drill-string is subject to axial and torsional effects, much like a flexible rod. Because of this elasticity, the force and velocity propagation can be modeled by \textit{wave equations}.

In this particular model, the actuation is the velocity of the travelling block, \textit{i.e.}, the axial velocity of the drill-string on the surface. Although not considered here, the rotational velocity also influences the rate of penetration (ROP) in a real scenario. The model output is the Weight On Hook (WOH) which somewhat models the Weight On Bit (WOB). The WOB estimates the contact between the drill bit and the rock formation and it is the downhole boundary condition to be controlled. In \cite{AK:2019}, the authors have discussed the feasibility of controlling  the hook load to optimize ROP while drilling. 

The key point that enables such an approach is the concept of bit foundering \cite{AK:2019}, \textit{i.e.}, the fact that ROP tapers off (and sometimes starts decreasing) with increasing weight on bit past the foundering point. This makes the static mapping between ROP and weight on bit upwards convex in an interval around the foundering point. This transfers to an upwards convex static
mapping between the equilibrium hook load set point and
feed rate. Consequently, these signals can be used as the plant input and output for the design of a drilling control system. Hence, this physical application motivates our ES scheme for static maps with actuation dynamics described by wave PDEs, as depicted in Fig.~\ref{ch10.fig:cascade}.


\section{Deep-Sea Cable-Actuated Source Seeking}

\begin{figure}[ht]
\centering
\includegraphics[width=6.0cm]{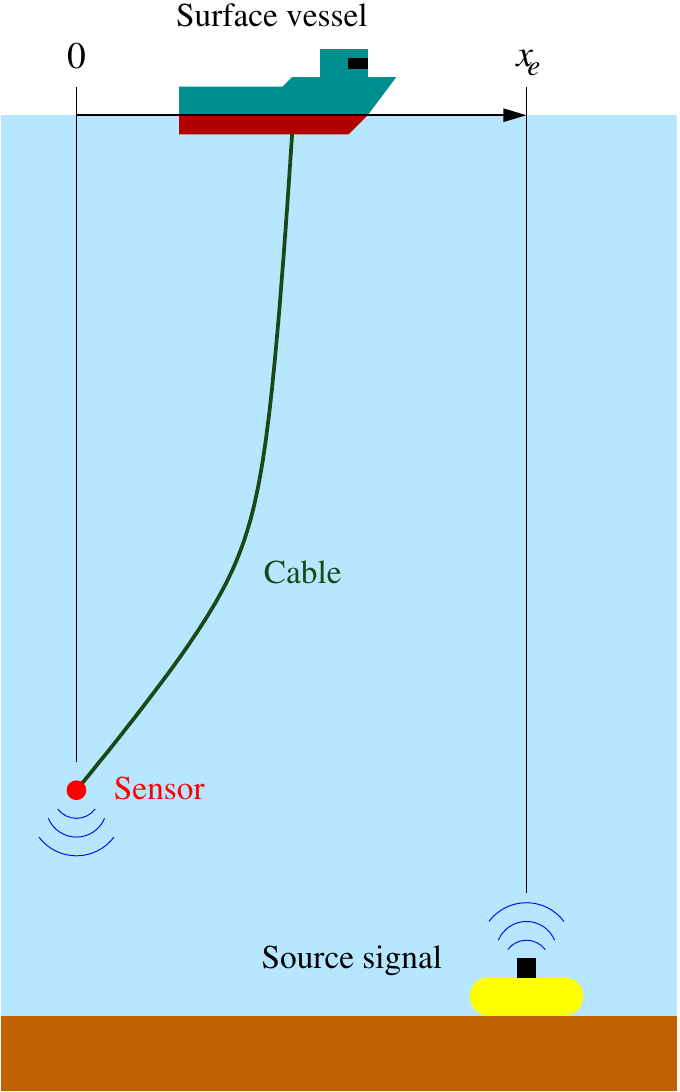}
\vspace{0.5cm}
\caption{\textcolor{black}{Motivating example - underwater search: $x_e$ represents the relative linear position of the source signal with respect to the sensor. The control task aims to drive the sensor to the source signal, meaning that  $x_e(t) \to 0$ (or to a small neighborhood of zero) as $t\to+\infty$.}} 
\label{fig:undersea}
\end{figure}

\textcolor{black}{The application is illustrated in Fig.~\ref{fig:undersea}} 
and involves a deep-sea cable-actuated source seeking. In this scenario, a sensor is suspended on a cable and moved through it from the sea surface using a surface vessel. The sensor operates without position awareness, primarily due to the challenging undersea environment. The task at hand is to locate the source signal as closely as possible. \textcolor{black}{No external fluid flow (e.g., water current) is considered, and the dynamics of the boat is ignored for simplicity \cite{F:2011}.} 


The algorithm proposed here 
is designed to be applicable to such a source-seeking scenario.
\textcolor{black}{The objective of source seeking is to find the source of a signal of an unknown concentration field, which can be chemical, acoustic, electromagnetic, etc. The sensor captures this field, and its strength decays with distance, reaching a maximum at position $x_e$ (relative to the coordinate system of the surface vessel). } Beyond finding the source signal, it also ensures the stabilization of the cable's motion. As the signal source becomes deeper or the cable length increases, the task becomes less demanding for the surface vessel, thanks to the high natural frequency of the longer cable, which reduces the need for rapid vessel motion. Nevertheless, achieving stability in the PDE-compensating ES algorithm may require a more extensive memory.


The cable of this application is represented by a string described by the following PDE model over an interval \textcolor{black}{$x \in [0,D]$:}
\begin{align}
&    \varepsilon \textcolor{black}{\alpha_{tt}} = (1+d\partial_{t})\textcolor{black}{\alpha_{xx}},\label{wavekv_actuator_sea}\\
&    \textcolor{black}{\alpha_{x}(0,t)} = 0\label{boundary_actuator_sea},\\
&    \textcolor{black}{\alpha(0,t)} = \text{measured} \label{boundary_sensor_sea},\\
&    \textcolor{black}{\alpha(D,t)} = \text{controlled}.\label{actuator_sea}
\end{align}
\noindent
\textcolor{black}{Equations (\ref{wavekv_actuator_sea})--(\ref{actuator_sea}) represent the dynamics of a string controlled at the end $x = D$, pinned to the surface vessel, and with a free end at $x = 0$, where the sensor is located.} 
\textcolor{black}{The term $\alpha(x,t)$ in (\ref{wavekv_actuator_sea}) represents the state variable of the PDE dynamics governing the motion of the cable. Equations (\ref{boundary_actuator_sea})--(\ref{actuator_sea}) serve as boundary conditions.} 
%
The constants $\varepsilon$, $d$ and $D$ are positive. \textcolor{black}{The constant $D$ 
physically corresponds to length of the cable.}  The value $1/\varepsilon$ represents the ``stiffness'' of the string, which can be expressed as 
%
$E/\rho$, where $E$ denotes Young's modulus and $\rho$ the density of the material. The term $d\partial_{t}$ models the ``Kelvin-Voigt'' damping, representing the internal material damping, not the damping that arises due to the viscous interaction of the string with the surrounding medium. We assume that this model takes into account \textcolor{black}{a small amount of damping ($d$)}, which is a realistic consideration in any material.  
%
%
We do not rely on the Kelvin-Voigt term as a source of energy dissipation; instead, we use it as a means of enhancing the controllability of the model (\ref{wavekv_actuator_sea})--(\ref{actuator_sea}).

\subsection{Scalar Maps with Actuation PDE Dynamics}
Now, we consider actuation dynamics described by a wave equation containing Kelvin-Voigt damping with $\varepsilon = 1$, $\theta(t) \in \mathbb{R}$ and \textcolor{black}{the sensor} 
$\Theta(t) \in \mathbb{R}$ given by 

\begin{align}
&    \Theta(t) = \alpha(0,t),\label{Theta_p_actuator}\\
&    \partial_{tt}\alpha(x,t) = \partial_{xx}\alpha(x,t) +     d\partial_{xxt}\alpha(x,t),\label{wavekv_actuator}\\
&    \partial_{x}\alpha(0,t) = 0,\label{boundary_actuator}\\
&    \alpha(D,t) = \theta(t),\label{theta_actuator}
\end{align}
\noindent
where $\alpha:[0,D]\times\mathbb{R}_{+}\rightarrow\mathbb{R}$, and $D$ is the known domain length, as mentioned before. \textcolor{black}{The output signal measured with the sensor is represented by the unknown static map} 
\begin{equation}
y(t) = Q(\Theta(t)),
\label{eq:initial_output_static_map}
\end{equation}
\textcolor{black}{with input $\Theta(t)$ in (\ref{Theta_p_actuator}).} 

\begin{figure}[ht]
\begin{center}
\includegraphics[scale=0.44]{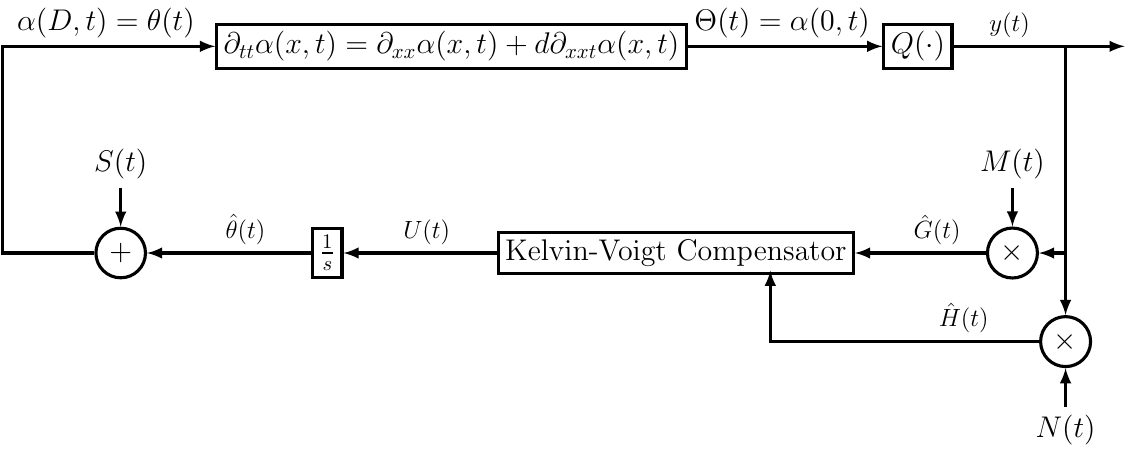}    
\vspace{-1cm}
\caption{Gradient extremum seeking control loop.} 
\label{fig:esc}
\end{center}
\end{figure}

\textcolor{black}{We assume that the unknown nonlinear map is locally quadratic}, such as in (\ref{extrachapter.eq:static_map_CSm}), resulting in (\ref{extrachapter.eq:final_output_static_map_CSm}). 
%
%
%
Adopting the proposed scheme in \cite{Tiago} and combining (\ref{Theta_p_actuator})--(\ref{theta_actuator}) with the ES approach, the closed-loop ES with actuation dynamics governed by the Kelvin-Voigt PDE is illustrated in Fig. \ref{fig:esc}.

\subsection{Trajectory Generation for the Probing Signal}
The perturbation $S(t)$ is adapted from the basic ES scheme in order to accommodate 
actuation dynamics. The trajectory generation problem, as described in \cite[Chapter12]{krstic2008boundary}, is outlined as follows: 
%
%
\begin{align}
&    S(t) := \beta(D,t), \label{perturbation_beta}\\
&    \partial_{tt}\beta(x,t) = \partial_{xx}\beta(x,t) +     d\partial_{xxt}\beta(x,t),\label{wavekv_beta}\\
&    \partial_{x}\beta(0,t) = 0,\label{boundary_beta}\\
&    \beta(0,t) = a\sin{(\omega t)},\label{initial_cond_beta}
\end{align}
\noindent
where $\beta:[0,D]\times\mathbb{R}_{+}\rightarrow\mathbb{R}$. The explicit solution of (\ref{perturbation_beta}) is derived for the reference trajectory
\begin{equation}
\beta(D,t):=\beta^{r}(D,t) = S(t), \; \beta(0,t):=\beta^{r}(0,t) = a\sin{(\omega t)}.
\label{eq:reference_trajectory}
\end{equation}
\noindent
This solution is found by postulating the reference solution  $\beta^{r}(x,t)$ as a power series of the spatial variable with time dependent coefficients: 
$\beta^{r}(x,t) = \sum_{i=0}^{\infty} a_{i}(t)\frac{x^{i}}{i!}$, as in \cite{Laroche1998MotionPF}. The string reference solution is given by \cite{SKSB:2009}
\begin{equation}
\beta^{r}(x,t) = -a\dfrac{j}{2}\left[\cosh{(j\sigma x})e^{j\omega t}-\cosh{(j\bar{\sigma} x})e^{-j\omega t}\right]
\label{eq:reference_sol}
\end{equation}
\noindent
with $\sigma = \dfrac{\omega}{\sqrt{1+j\omega d}}$ and $\bar{\sigma}$ being its complex conjugate. Equation (\ref{eq:reference_sol}) can be written as the purely real function
\begin{equation}
\begin{aligned}
\beta^{r}(x,t) ={} & \dfrac{a}{2}\Big[e^{\hat{\beta}(\omega)x}\sin{(\omega t + \overline{\beta}(\omega)x)}+\\
               & e^{-\hat{\beta}(\omega)x}\sin{(\omega t- \overline{\beta}(\omega)x)}\Big],
\label{eq:reference_solution2}
\end{aligned}
\end{equation}
\noindent
where the real function $\overline{\beta}(\omega)$ and $\hat{\beta}(\omega)$ are defined as 
\begin{equation}
\overline{\beta}(\omega) = \omega\sqrt{\dfrac{\sqrt{1+\omega^{2}d^{2}}+1}{2(1+\omega^{2}d^{2})}},
\label{eq:beta}
\end{equation}
\begin{equation}
\hat{\beta}(\omega) = \omega\sqrt{\dfrac{\sqrt{1+\omega^{2}d^{2}}-1}{2(1+\omega^{2}d^{2})}}.
\label{eq:beta_hat}
\end{equation}
On the other hand, the demodulation signals $M(t)$ and $N(t)$, used for estimating 
the gradient and Hessian, respectively, of the static map by multiplying them by the output $y(t)$, are defined in \cite{GKN:2012} as
\begin{equation}
\hat{H}(t) = N(t)y(t)\;\; \text{with} \;\; N(t) = -\dfrac{8}{a^{2}}\cos{(2\omega t)}.
\label{eq:hessian}
\end{equation}
\begin{equation}
G(t) = M(t)y(t)\;\; \text{with} \;\; M(t) = \frac{2}{a}\sin{(\omega t)}.
\label{eq:gradient_USOpen}
\end{equation}

\subsection{Estimation Errors and Error Dynamics}

Since our objective is to find $\Theta^{*}$, which corresponds to the optimal unknown actuator $\theta(t)$, we introduce the following estimates and the estimation errors

\begin{equation}
\hat{\theta}(t) = \theta(t) - S(t),\;\;\; \hat{\Theta}(t) = \Theta(t) - a\sin{(\omega t)},
\label{eq:estimated}
\end{equation}

\begin{equation}
\tilde{\theta}(t) := \hat{\theta}(t) - \Theta^{*},\;\;\; \vartheta(t) := \hat{\Theta}(t) - \Theta^{*}.
\label{eq:estimated_errors}
\end{equation}

Let $\bar{\alpha}:[0,D]\times\mathbb{R}_{+}\rightarrow\mathbb{R}$ be defined as $\bar{\alpha}(x,t) := \alpha(x,t) - \beta(x,t) - \Theta^{*}$. Manipulating (\ref{Theta_p_actuator})--(\ref{theta_actuator}) and (\ref{perturbation_beta})--(\ref{initial_cond_beta}) with the help of (\ref{eq:estimated}) and (\ref{eq:estimated_errors}), we get:

\begin{align}
&    \vartheta(t) = \bar{\alpha}(0,t), \label{dynamics_baractuator}\\
&    \partial_{tt}\bar{\alpha}(x,t) = \partial_{xx}\bar{\alpha}(x,t) +     d\partial_{xxt}\bar{\alpha}(x,t),\label{wavekv_baractuator}\\
&    \partial_{x}\bar{\alpha}(0,t) = 0,\label{boundary_baractuator}\\
&    \bar{\alpha}(D,t) = \tilde{\theta}(t).\label{theta_baractuator}
\end{align}

\textcolor{black}{The error-dynamics} is obtained by taking the time derivative of (\ref{dynamics_baractuator})--(\ref{theta_baractuator}) and using  $\dot{\tilde{\theta}}\!=\!U(t)$ and \textcolor{black}{$u(x,t)\!=\!\bar{\alpha}_t(x,t)$:} 
\begin{align}
&    \dot{\vartheta}(t) = u(0,t), \label{dynamics_error}\\
&    \partial_{tt}u(x,t) = \partial_{xx}u(x,t) +     d\partial_{xxt}u(x,t),\label{wavekv_error}\\
&    \partial_{x}u(0,t) = 0,\label{boundary_actuator_error}\\
&    u(D,t) = U(t).\label{eq:30}
\end{align}
%

\subsection{Boundary Extremum Seeking Control Law} \label{debinhadoidinha}


\begin{figure}[ht]
\begin{center}
\includegraphics[width=8.3cm]{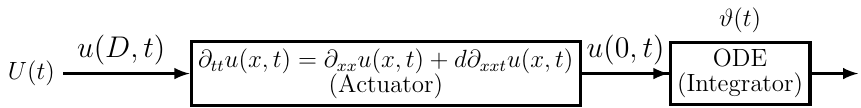}    
\caption{The cascade of the PDE dynamics and the ODE integrator.} 
\label{fig:ode_plant}
\end{center}
\end{figure}

We consider the PDE-ODE cascade shown in Fig. \ref{fig:ode_plant}, and use the backstepping transformation
\begin{equation}
w(x,t) = u_{av}(x,t) - \overline{K}\int_{0}^{x} k(x,\sigma)u_{av}(\sigma,t) \,d\sigma - \overline{K}\vartheta_{av}(t)
\label{eq:backstepping}
\end{equation}
\noindent
to transform the original (\ref{dynamics_error})--(\ref{eq:30}) into the target system
\begin{align}
&    \dot{\vartheta}_{av}(t)=\overline{K}\vartheta_{av}(t)+w(0,t),\quad \overline{K}<0, \label{eq:33} \\
&    w_{tt}=(1+d\partial_t)(w_{xx}-c w),\quad c>0, \label{eq:34}\\
&    w_x(0,t)=0, \label{eq:35}\\
&     w(D,t)=0. \label{eq:36}
\end{align}

The gain kernel PDE $k(x,\sigma)$ comes from the solution of (see \cite[Section~4.2]{krstic2008boundary}) 
\begin{align} \label{socagadadeultimahora}
k_{xx}=k_{\sigma \sigma}+ck\,, \quad k_{\sigma}(x,0)=0\,, \quad k(x,x) = \frac{c}{2}x. 
\end{align}
\textcolor{black}{The solution to the PDE in (\ref{socagadadeultimahora}) is obtained through a summation of successive approximation series \cite[Section~4.4]{krstic2008boundary}:} 
\begin{equation}
    k(x,\sigma)=-cx \frac{I_1\Big(\sqrt{c(x^2-\sigma^2)}\Big)}{\sqrt{c(x^2-\sigma^2)}}
\end{equation}
%
%
%
%
%
\noindent
and, from (\ref{eq:backstepping}) and (\ref{eq:36}), the average control law is given by:
\begin{small}
\begin{equation}
u_{av}(D,t)=\overline{K}\vartheta_{av}(t)-\overline{K}\int_0^DcD\frac{I_1\Big(\sqrt{c(D^2-\sigma^2)}\Big)}{\sqrt{c(D^2-\sigma^2)}}u_{av}(\sigma,t) d\sigma,
\label{eq:control_backstepping_1}
\end{equation}
\end{small}
$\!\!$where $I_1$ is \textcolor{black}{the modified Bessel function \cite[Appendix~A.2]{krstic2008boundary}.}  
%
Thus, introducing a result of \cite{GKN:2012}, the averaged version of the gradient and Hessian estimate are calculated as 
\begin{equation}
G_{av}(t) = H\vartheta_{av}(t), \;\;\; \hat{H}_{av}(t) = H.
\label{eq:grad_hessian_estimate}
\end{equation}

From (\ref{eq:30}) and (\ref{eq:control_backstepping_1}), choosing $\overline{K} \!=\! KH$ with $K\!>\!0$ and plugging the average gradient and Hessian estimates (\ref{eq:grad_hessian_estimate}), we obtain
\begin{small}
\begin{equation}
U_{av}(t) = KG_{av}(t)- KH\int_0^DcD\frac{I_1\Big(\sqrt{c(D^2-\sigma^2)}\Big)}{\sqrt{c(D^2-\sigma^2)}}u_{av}(\sigma,t) d\sigma.
\label{eq:control_backstepping_2}
\end{equation}
\end{small}

We introduce a low-pass filter to obtain the non-average controller 
\begin{equation}
\begin{split}
U(t) &= \dfrac{\overline{c}}{s + \overline{c}} \Bigg\{ K \Bigg[ G(t) - \Hat{H}(t) \\
& \quad \times \int_0^D \! cD \frac{I_1\Big(\sqrt{c(D^2 - \sigma^2)}\Big)}{\sqrt{c(D^2 - \sigma^2)}} u(\sigma,t) \, d\sigma \Bigg] \Bigg\}.
\end{split}
\label{eq:control_backstepping_3}
\end{equation}

\noindent
with $\overline{c}\rightarrow +\infty$ sufficiently large. 

\section{Materials Phase Change PDE ES-Control: From Fixed Domain to Moving Boundary} 

We next present the design and analysis of the ES for static maps with input governed by a PDE of the diffusion type defined on a time varying spatial domain described by an ODE. We compensate for the average-based actuation dynamics by a controller via backstepping transformation for the moving boundary, which is utilized to transform the original coupled PDE-ODE into a target system whose exponential stability of the average equilibrium of the average system is proved. 

\subsection{One-phase Stefan Problem}


\begin{figure}[ht]
\begin{center}
\includegraphics[width=8cm]{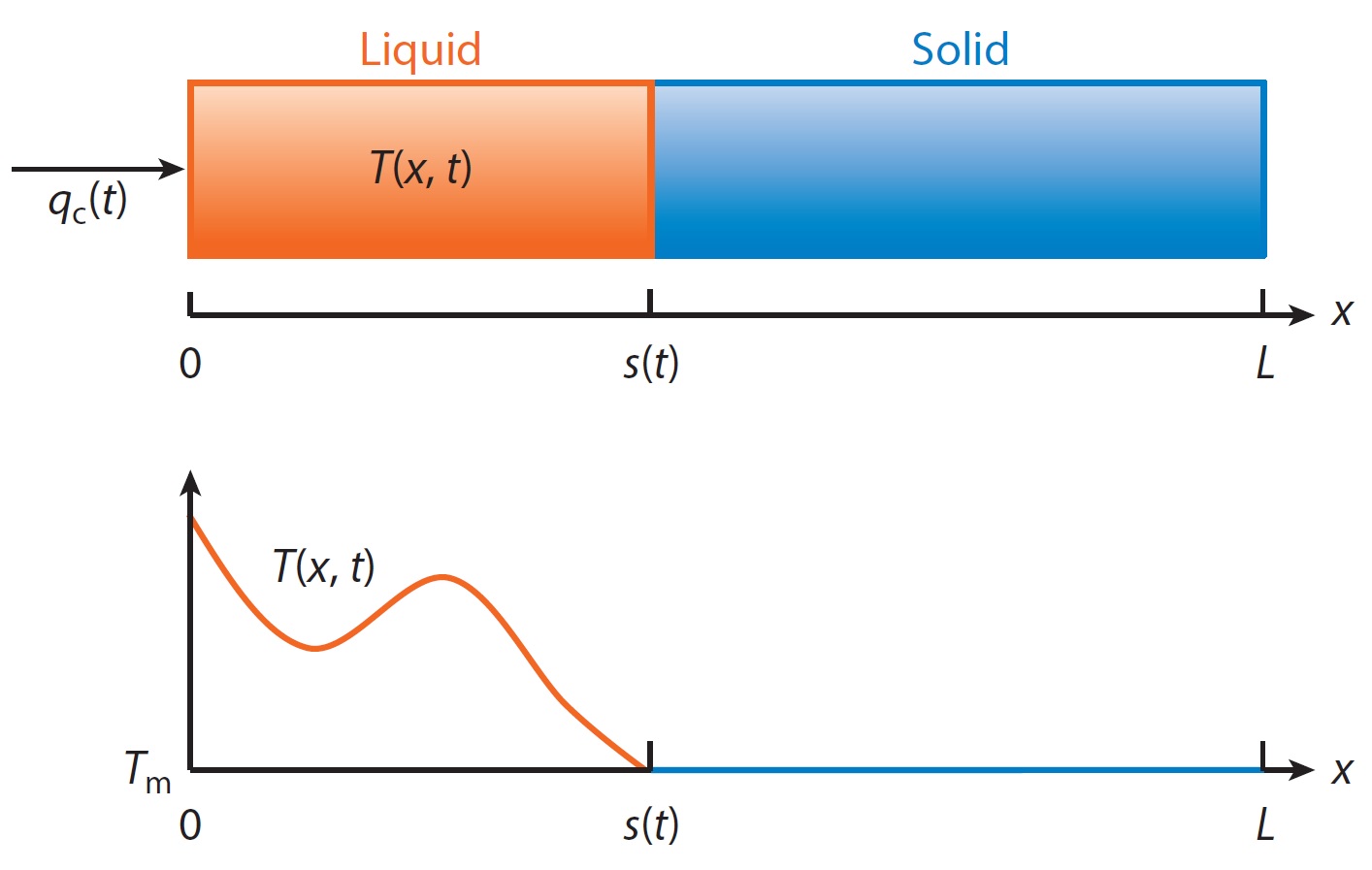}    
\caption{Schematic of one-phase Stefan problem \cite{Shumon_green_book:2020,AR_KK:2022}.
The temperature profile in the solid phase is assumed to be a uniform melting temperature.}
\label{extrachapter.fig:stefan}
\end{center}
\end{figure}

\begin{figure*}[ht]
\begin{center}
\includegraphics[width=14.0cm]{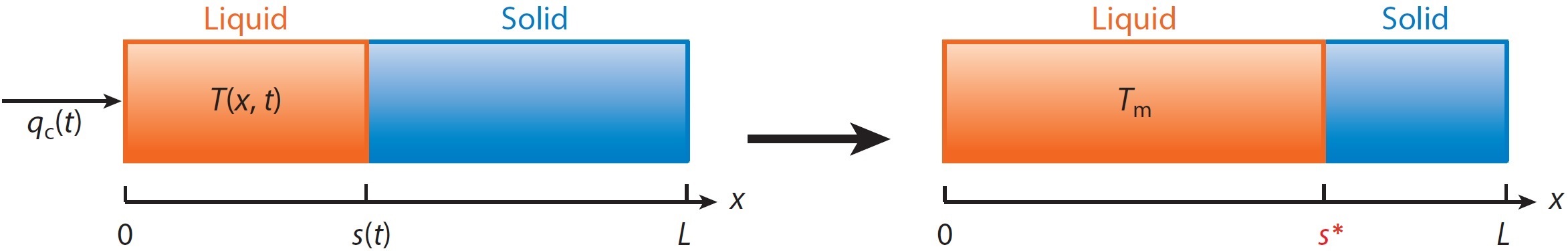}    
\caption{Control objective of the Stefan problem. We aim to design a heat flux input $q_c(t)$ such that the interface position $s(t)$ is driven to the setpoint position $s^*$.}
\label{extrachapter.fig:stefan_objectives}
\end{center}
\end{figure*}

\begin{figure}[ht]
    \begin{center}
    \includegraphics[width=8cm]{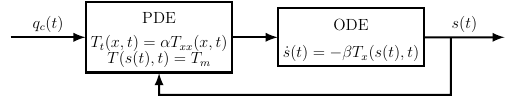} 
    \caption{The cascade of the PDE dynamics and the ODE system.}
    \label{extrachapter.fig:ode_plant}
    \end{center}
\end{figure}

The physical model which describes the 1-D Stefan problem in a pure one-component material of length $L$ is described in Fig.~\ref{extrachapter.fig:stefan}. The domain $[0,L]$ is divided into two sub-domains $[0,s(t)]$ and $[s(t),L]$ which represents the liquid phase and the solid phase, respectively. The system is controlled by the heat flux $q_{c}(t)$ at $x=0$, because we are dealing with a Neumann boundary actuation as shown below:
\begin{small}
\begin{align}
    T_{t}(x,t) &= \alpha T_{xx}(x,t),\quad x\in (0,s(t)),\quad \alpha = \dfrac{k}{\rho C_{p}}\label{extrachapter.T1}\\
    -kT_{x}(0,t) &= q_{c}(t)\label{extrachapter.T2}\\
    T(s(t),t) &= T_{m}\label{extrachapter.T3}\\
    \dot{s}(t) &= -\beta T_{x}(s(t),t),\quad \beta = \dfrac{k}{\rho \Delta H^{*}}, \label{extrachapter.T4}
\end{align}
\end{small}
$\!\!$where $T(x,t)$, $T_{m}$, $q_{c}(t)$, $k$, $\rho$, $C_{p}$ and $\Delta H^{*}$ are the distributed temperature of the liquid phase, melting temperature, manipulated heat flux, liquid heat conductivity, liquid density, liquid heat capacity and latent heat of fusion, respectively. Equations (\ref{extrachapter.T2}) and (\ref{extrachapter.T3}) are the boundary conditions of the system and (\ref{extrachapter.T4}) is the Stefan condition, which describes the dynamics of the moving boundary. Fig.~\ref{extrachapter.fig:ode_plant} shows the block diagram of the PDE-ODE cascade represented by equations (\ref{extrachapter.T1})-(\ref{extrachapter.T4}).


\subsection{Actuation dynamics and output signal} \label{extrachapter.ESC}

For the sake of simplicity, we consider actuation dynamics which are described by a heat equation with $\alpha,\beta,k = 1$, $\theta(t) \in \mathbb{R}$ and the propagated actuator $\Theta(t) \in \mathbb{R}$ given by 
\begin{align}
    \dot{\Theta}(t) = \dot{s}(t) &= -\alpha_{x}(s(t),t)),\quad x \in (0,s(t))\label{extrachapter.Theta_p_actuator}\\
    \partial_{t}\alpha(x,t) &= \partial_{xx}\alpha(x,t)\label{extrachapter.stefan_actuator}\\
    \alpha(s(t),t) &= 0\label{extrachapter.boundary_actuator}\\
    -\partial_{x}\alpha(0,t) &= \theta(t)\label{extrachapter.theta_actuator},
\end{align}
where $\alpha:[0,s(t)]\times\mathbb{R}_{+}\rightarrow\mathbb{R}$ is $\alpha(x,t) = T(x,t) - T_{m}$ and $s(t) = \Theta(t)$ is the unknown interface represented as the moving boundary. The output is measured by the unknown static map with input (\ref{extrachapter.Theta_p_actuator}), according to (\ref{extrachapter.eq:initial_output_static_map_CSm}).  
The ES goal is to optimize an unknown static map Q($\cdot$) using a real-time optimization control with optimal unknown output $y^{*}$ and optimizer $\Theta^{*}$ as well as measurable output $y$ and input $\theta$. Consequently, the control objectives of the Stefan problem are achieved, \textit{i.e.}, $\lim\limits_{t\rightarrow \infty}s(t) = s^*$ and $\lim\limits_{t\rightarrow \infty}T(x,t) = T_m\,, \forall x \in [0, s^*]$, as illustrated in Fig.~\ref{extrachapter.fig:stefan_objectives}. 

\begin{figure}[ht]
\begin{center}
    \includegraphics[width=8cm]{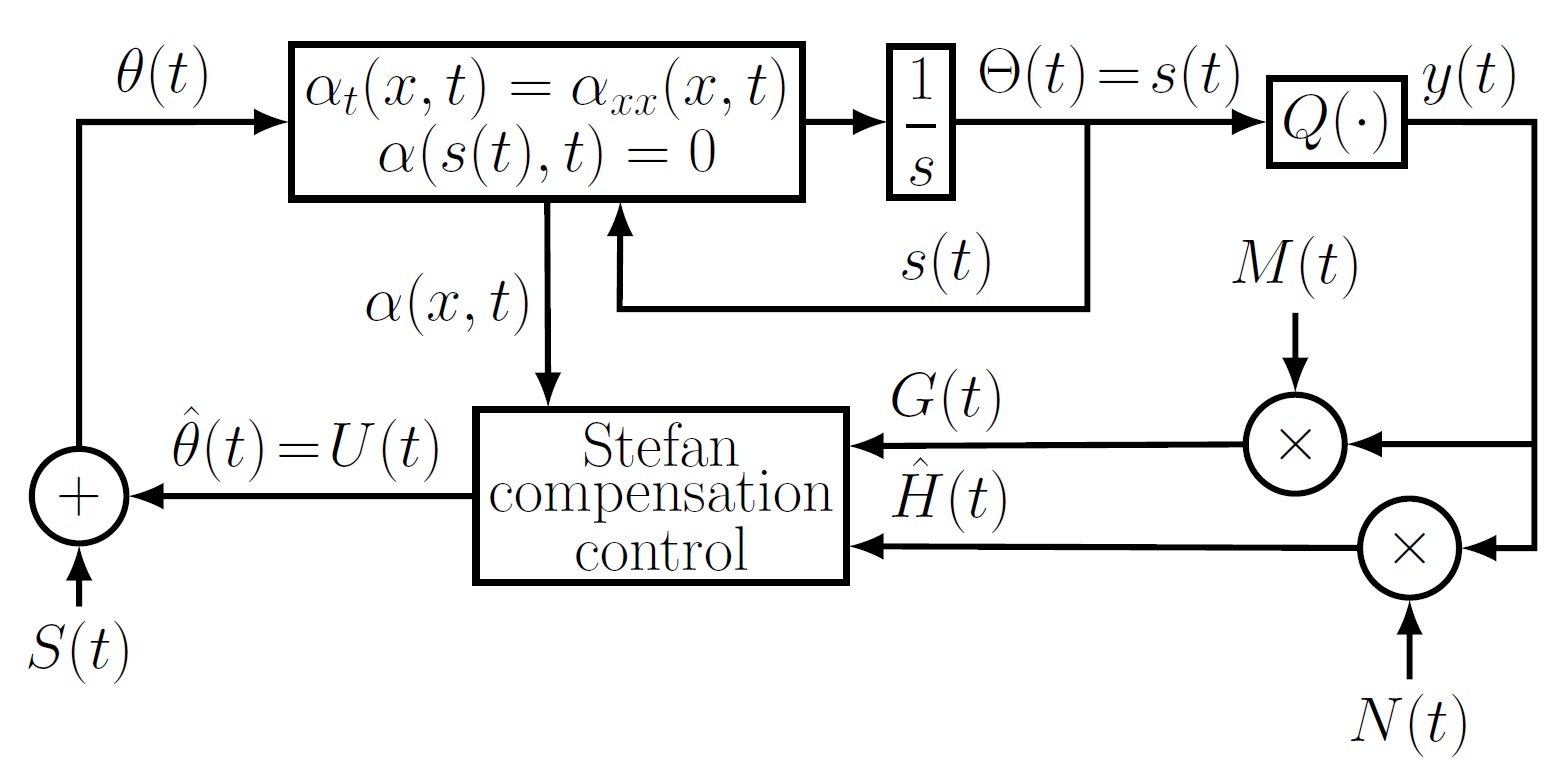} 
	\caption{Extremum seeking control loop applied to the one-phase Stefan problem.}
	\label{extrachapter.fig:esc}
\end{center}
\end{figure}

On the other hand, the unknown nonlinear map is locally quadratic, given by (\ref{extrachapter.eq:static_map_CSm}), such that the output of the static map is simply 
(\ref{extrachapter.eq:final_output_static_map_CSm}). 
Adapting the proposed scheme in \cite{Tiago} and combining (\ref{extrachapter.Theta_p_actuator})-(\ref{extrachapter.theta_actuator}) with the ES approach, the closed-loop ES with actuation dynamics is shown in Fig.~\ref{extrachapter.fig:esc}.


The demodulation signal $N(t)$ which is used to estimate the Hessian of the static map by multiplying it with the output $y(t)$ of the static map is defined in \cite{GKN:2012} as
\begin{equation}
\hat{H}(t) = N(t)y(t)\;\; \text{with} \;\; N(t) = -\dfrac{8}{a^{2}}\cos{(2\omega t)}
\label{extrachapter.eq:hessian}
\end{equation}
whereas the signal $M(t)$ is used to estimate the gradient of the static map as follows:
\begin{equation}
    G(t) = M(t)y(t)\;\; \text{with} \;\; M(t) = \frac{2}{a}\sin{(\omega t)}.
\label{extrachapter.eq:ditherSM}
\end{equation}

\subsection{Additive Probing Signal}

The perturbation $S(t)$ is adapted from the basic ES to the case of PDE actuation dynamics. The trajectory generation problem as in \cite{krstic2008boundary} is described as follows:
\begin{align}
    S(t) &:= -\partial_{x}\beta(0,t),\quad x \in (0,s(t)) \label{extrachapter.perturbation_beta}\\
    \partial_{t}\beta(x,t) &=          \partial_{xx}\beta(x,t)\label{extrachapter.stefan_beta}\\
    \beta(s(t),t) &= 0\label{extrachapter.boundary_beta}\\
    \beta_{x}(s(t),t) &= -a\omega\cos{(\omega t)}\label{extrachapter.initial_cond_beta},
\end{align}
where $\beta:[0,s(t)]\times\mathbb{R}_{+}\rightarrow\mathbb{R}$. The explicit solution of (\ref{extrachapter.perturbation_beta}) is found respectively for the reference trajectory and the reference solution postulated by a power series \cite{dunbar2003motion}:
\begin{align}
    &s(t) = a\sin{(\omega t)} \label{extrachapter.pertubation_s(t)}\\
    &\beta(x,t) = \sum_{i=0}^{\infty} \label{extrachapter.power_series} \dfrac{a_{i}(t)}{i!}[x-s(t)]^{i}.
\end{align}
We can calculate the first coefficients of the power series replacing the boundary conditions (\ref{extrachapter.boundary_beta}) and (\ref{extrachapter.initial_cond_beta}) at (\ref{extrachapter.power_series}), such that
\begin{equation}
    a_{0}(t) = 0,\;\;\;\;\; a_{1} = -\dot{s}(t).
\end{equation}
The general expression $a_{i}(t) = \dot{a}_{i-2}(t) - a_{i-1}(t)\dot{s}(t)$ is obtained by substituting (\ref{extrachapter.power_series}) in (\ref{extrachapter.stefan_beta}). We provide here the analytic expression of the first four coefficients of the series (\ref{extrachapter.power_series}) so that one can see how the successive derivatives of $s(t)$ appear: 
\begin{align}
    &a_{2}(t) = \dot{s}(t)^{2}\\
    &a_{3}(t) = \ddot{s}(t) - \dot{s}(t)^{3}\\
    &a_{4}(t) = \ddot{s}(t)^{2} + \ddot{s}(t)\dot{s}(t) + \dot{s}(t)^{4}.
\end{align}
The trajectory generation solution which provide all terms of the power series (\ref{extrachapter.power_series}) is given by \cite{hill1967parabolic}
\begin{equation}
    \beta(x,t) =  \sum_{i=0}^{\infty} \dfrac{1}{(2i)!}\dfrac{\partial^{i}}{\partial t^{i}}[x-s(t)]^{2i}. \label{extrachapter.beta_solution}
\end{equation}
Although (\ref{extrachapter.beta_solution}) is not an explicit expression, choosing suitable values for $a$ and $\omega$ in (\ref{extrachapter.pertubation_s(t)}), the series converges with few iterations of the infinite sum, getting the desirable sinusoidal signal $s(t)$ in the output of the integrator.

According to (\ref{extrachapter.perturbation_beta}), we take the spatial derivative of (\ref{extrachapter.beta_solution}) and substitute $x=0$, thus arriving at the final expression of
\begin{equation}
    S(t) = -\sum_{i=0}^{\infty} \dfrac{1}{(2i-1)!}\dfrac{\partial^{i}}{\partial t^{i}}[-a\sin{(\omega t)}]^{2i-1}.
\end{equation}

\subsection{Estimation Errors and PDE-Error Dynamics}

Since our objective is to find $\Theta^{*}$, which corresponds to the optimal unknown actuator $\theta(t)$, we introduce the following estimation errors: 
\begin{equation}
\hat{\theta}(t) = \theta(t) - S(t),\;\;\; \hat{\Theta}(t) = \Theta(t) - a\sin{(\omega t)},
\label{extrachapter.eq:estimated}
\end{equation}
\begin{equation}
\tilde{\theta}(t) := \hat{\theta}(t) - \Theta^{*},\;\;\; \vartheta(t) := \hat{\Theta}(t) - \Theta^{*},
\label{extrachapter.eq:estimated_errors}
\end{equation}
recalling that $\Theta(t) := s(t)$. Combining $\hat{\Theta}(t)$ in (\ref{extrachapter.eq:estimated}) and (\ref{extrachapter.eq:estimated_errors}), we obtain the relation between the propagated estimation error $\vartheta(t)$, the propagated input $\Theta(t)$ and the optimizer of the static map $\Theta^{*}$: 

\begin{equation}
    \Theta(t) - \Theta^{*} = \vartheta(t) + a\sin{(\omega t)}. \label{extrachapter.vasin}
\end{equation}

Let us define

\begin{equation}
    u(x,t) = \alpha(x,t) - \beta(x,t) \label{extrachapter.u},
\end{equation}
\begin{equation}
    \hat{\theta}(t) = U(t). \label{extrachapter.utheta}
\end{equation}

By (\ref{extrachapter.Theta_p_actuator})-(\ref{extrachapter.theta_actuator}) and (\ref{extrachapter.perturbation_beta})-(\ref{extrachapter.initial_cond_beta}), with the help of (\ref{extrachapter.eq:estimated}) and (\ref{extrachapter.eq:estimated_errors}), we have our original system:
\begin{align}
    \dot{\vartheta}(t) &= -u_{x}(s(t),t),\quad x \in (0,s(t)) \label{extrachapter.dynamics_error}\\
    u_{t}(x,t) &= u_{xx}(x,t) \label{extrachapter.wavekv_error}\\
    u(s(t),t) &= 0\label{extrachapter.boundary_actuator_control}\\
    -u_{x}(0,t) &= U(t).\label{extrachapter.controller_error}
\end{align}

\begin{figure*}[ht]
\begin{center}
\includegraphics[width=12.0cm]{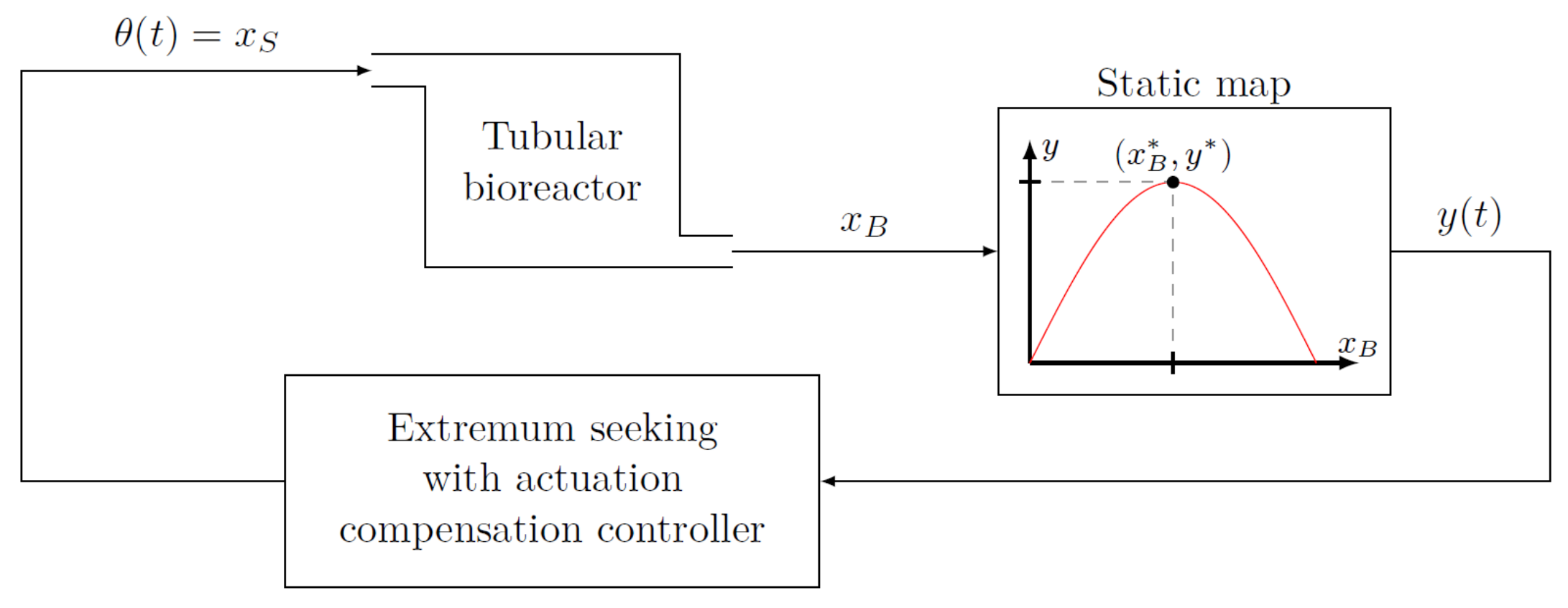}
\caption{Extremum seeking scheme for tubular bioreactors.}\label{fig:bioreactor}
\end{center}
\end{figure*}
\begin{figure*}[ht]
\begin{center}
\includegraphics[width=12.0cm]{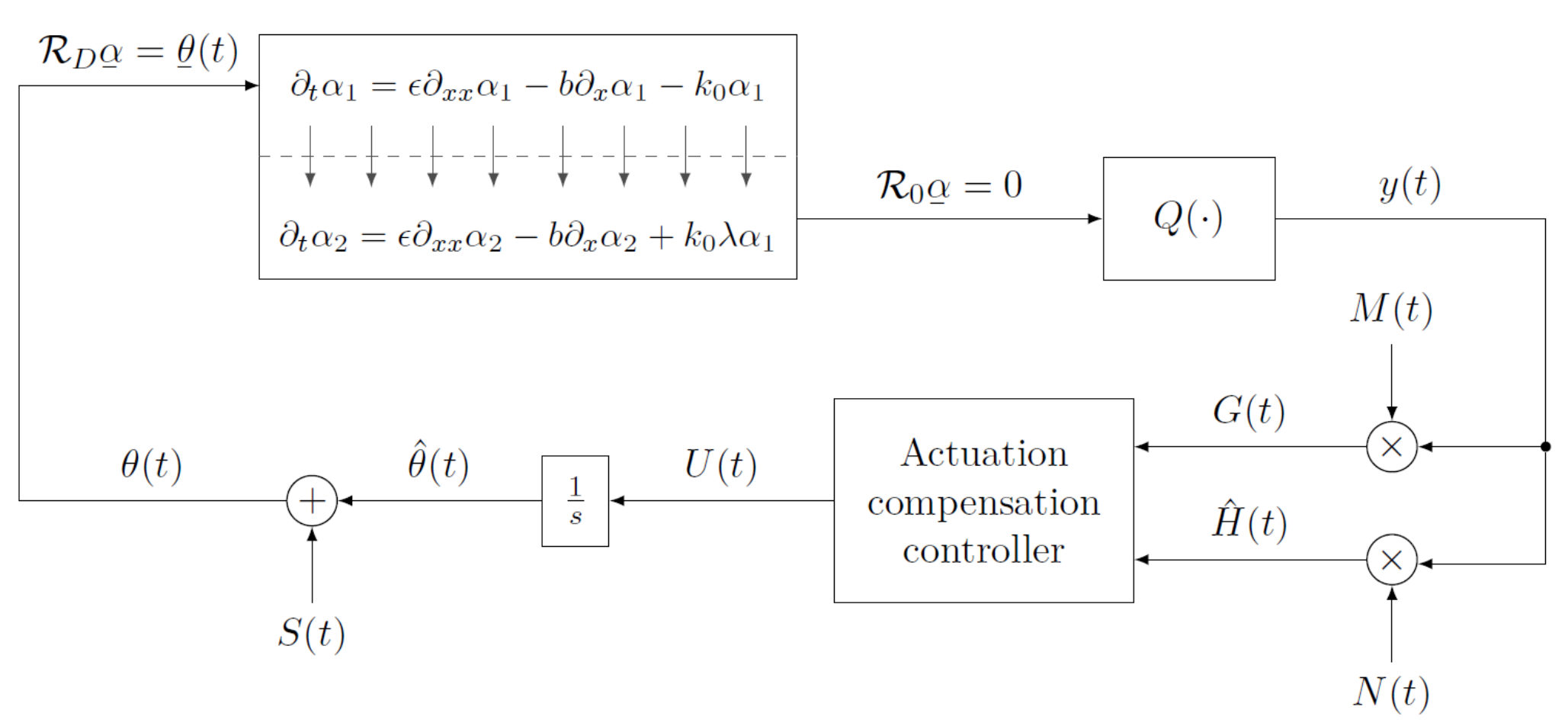}
\caption{Basic gradient ES scheme for tubular bioreactors modeled by RAD-PDEs.}\label{fig:bioreactor_ES}
\end{center}
\end{figure*}

\subsection{Stefan Compensation Control Law} \label{extrachapter.controller}


We consider the PDE-ODE cascade (\ref{extrachapter.dynamics_error})-(\ref{extrachapter.controller_error}) and use the backstepping transformation 
\begin{equation}
\begin{split}
w(x,t) &= u(x,t) - \bar{K}\int_{x}^{s(t)} (x-\sigma)u(\sigma,t)\,dy+\\
& -\bar{K}(x-s(t))\vartheta(t)
\label{extrachapter.eq:control_backstepping}
\end{split}
\end{equation}
with $\bar{K}>0$ as an arbitrary controller gain. Equation (\ref{extrachapter.eq:control_backstepping}) transforms (\ref{extrachapter.dynamics_error})-(\ref{extrachapter.controller_error}) into the target system:
\begin{align}
    \dot{\vartheta}(t) &= -\bar{K}\vartheta(t) - w_{x}(s(t),t),\,\,\,\,x \in (0,s(t)) \label{extrachapter.dynamics_target1}\\
    w_{t}(x,t) &= w_{xx}(x,t) + \bar{K}\dot{s}(t)\vartheta(t) \label{extrachapter.stefan_target}\\
    w_{x}(0,t) &= 0\label{extrachapter.boundary_target1}\\
    w(s(t),t) &= 0.\label{extrachapter.controller_target1}
\end{align}
The compensation controller can be obtained by taking the derivative of (\ref{extrachapter.eq:control_backstepping}) with respect to $t$ and $x$, respectively, along the solution (\ref{extrachapter.dynamics_error})-(\ref{extrachapter.controller_error}), and substituting $x=0$:
\begin{equation}
    U(t) = -\bar{K}\left( \vartheta(t) + \int_{0}^{s(t)} u(x,t)\,dx\right) \label{extrachapter.control_law}.
\end{equation}


Since we have no measurement on $\vartheta(t)$, (\ref{extrachapter.control_law}) is not applicable directly. Thus, introducing a result of \cite{GKN:2012}, the average version of the gradient and Hessian estimates are calculated by

\begin{equation}
G_{\rm av}(t) = H\vartheta_{\rm av}(t), \;\;\; \hat{H}_{\rm av}(t) = H.
\label{extrachapter.eq:grad_hessian_estimate}
\end{equation}

Averaging (\ref{extrachapter.control_law}), choosing $\bar{K} = KH$ with $K<0$, and plugging in the average gradient and Hessian estimates (\ref{extrachapter.eq:grad_hessian_estimate}), we obtain
\begin{equation}
    U_{\rm av}(t) = -KG_{\rm av}(t) - KH\int_{0}^{s_{\rm av}(t)} u_{\rm av}(x,t)\,dx. \label{extrachapter.average_control_law}
\end{equation}

We introduce a low-pass filter to the controller with the purpose of applying the average theorem for infinite-dimensional systems \cite{HL:1990}, such that
\begin{equation}
U(t) = \dfrac{c}{s+c}\Bigg\{K\Bigg[G(t)+
\hat{H}(t) \int_{0}^{s(t)} u(x,t)\,dx \Bigg]\Bigg\},
\label{extrachapter.eq:filter_control_law}
\end{equation}
for $c>0$ sufficiently large.

\section{Bioreactors}


Next, we present an application idea for this novel control concept: a tubular bioreactor. We consider the optimization problem in Fig.~\ref{fig:bioreactor} of a tubular bioreactor, where the goal is to operate the bioreactor at the unknown optimal product rate, \textit{e.g.},  growth of biomass by determining the optimal input, hence the substrate concentration, \textit{e.g.}, glucose concentration.
The product concentration $x_B$ is generally not measurable, unlike the product rate.

Since the static map, which determines the product rate depending on the product concentration $x_B$, is not known or only approximately known, we can apply our introduced control concept. A simple model of a tubular bioreactor is presented by Winkin et al. \cite{Winkin2000349}, where the chemical reaction of a reactant $R$ and a product $P$, given by $	R \rightarrow \lambda P$, 
where the stoichiometric coefficient $\lambda$ of the reaction is considered. In its simplest form, the model is linear, since the nonlinear reaction term is simplified to a linear kinetic model only depending on the reactant concentration. The control-loop structure to reach and operate the tubular bioreactor at the optimal (highest product rate) is shown in Fig.~\ref{fig:bioreactor_ES}. 

The boundary operators of the tubular reactor in Fig.~\ref{fig:bioreactor_ES} are defined as 
\begin{align}
\begin{split}
	\mathcal{R}_D\underline{\alpha} &= \begin{bmatrix}
	\epsilon \partial_x\alpha_1(D,t)-b\alpha_1(D,t) \\
	\epsilon \partial_x\alpha_2(D,t)-b\alpha_2(D,t)
	\end{bmatrix} = \begin{bmatrix}
	-b\theta(t)  \\
	0
	\end{bmatrix}\,, \\
	\mathcal{R}_0\underline{\alpha} &= \begin{bmatrix}
	\epsilon \partial_x\alpha_1(0,t) \\
	\epsilon \partial_x\alpha_2(0,t)
	\end{bmatrix} = \begin{bmatrix}
	0  \\
	0
	\end{bmatrix}\,.
\end{split}
\label{eq:bioreactor_boundary_op}
\end{align}
In a first step, the average-based infinite-dimensional controller to compensate the tubular reactor process, which are the actuator dynamics in this case, has to be derived. Note that the dynamics of the reactor has to be known to apply this control concept. Baccoli et al. \cite{baccoli2015boundary} considers a similar dynamics, \textit{i.e.}, coupled linear parabolic PDEs, especially reaction-diffusion equations, with boundary input in one variable. But we additionally have an ODE cascade, arising from the integrator in the control loop, which has to be stabilized. 
There seems to be no work that considers a coupled parabolic PDE-ODE cascade with two parabolic PDEs coupled in domain plus boundary input in one variable and boundary measurement in the other variable. The controller derivation and further derivations like the perturbation signal $S(t)$ for this application would go beyond the scope of this paper. 
We emphasize by invoking the averaging theorem for infinite-dimensional systems in ``Averaging Theorem for General Infinite-Dimensional Systems \cite{HL:1990}''  
would work, since the publication \cite{Winkin2000349} showed the analytic semigroup property of the operator, which describes the coupled parabolic PDE system in Fig.~\ref{fig:bioreactor_ES}. Furthermore, the exponential stability proof of the average system will follow the same steps, but with an extended Lyapunov-Krasovskii functional and more calculation steps. Possibly there will be some restrictions on the system parameters $\epsilon,\ b,\ k_0$ and $\lambda$.

\section{Infinite-Dimensional Light-Source Seeking}

While the wave equation is the most appropriate ``point of entry'' into the realm of hyperbolic
PDEs, beam equations are considered a physically relevant benchmark for control of
hyperbolic PDEs and structural systems in general.

The simplest beam model is the Euler–Bernoulli model
\begin{align}
&    \textcolor{black}{\alpha_{tt}} +  {\alpha_{xxxx}} = 0,\label{wavekv_actuator_light}\\
&    \textcolor{black}{\alpha_{x}(0,t)} = \alpha_{xx}(0,t) = \alpha_{xxx}(0,t) = 0  \quad \text{(free end condition)} \label{boundary_actuator_light},\\
&    \textcolor{black}{\alpha(0,t)} =  \text{clamped end condition} \label{boundary_sensor_light}.
\end{align}
The obvious difference between the PDEs (\ref{wavekv_actuator_sea})--(\ref{actuator_sea}) and (\ref{wavekv_actuator_light})--(\ref{boundary_sensor_light}) 
is in the number of spatial derivatives—the wave equation is second order in $x$, whereas the Euler-Bernoulli beam
model is fourth order in $x$. One consequence of this difference is that a wave equation requires
one boundary condition per end point (see (\ref{boundary_sensor_sea}) or (\ref{actuator_sea})), whereas the Euler-Bernoulli
beam model requires two boundary conditions per end point; see (\ref{boundary_actuator_light}) or (\ref{boundary_sensor_light}). A more
important difference is in the eigenvalues. Both the beam and the string models have all of
their eigenvalues on the imaginary axis. However, while the string eigenvalues are equidistant, the beam eigenvalues get further and further apart
as they go up the imaginary axis. This difference in the eigenvalue pattern is a consequence of the difference in the number of derivatives in $x$.

The reader might wonder how these differences translate into control. Is it obvious that a beam is more difficult to control than a string?
The answer is not clear and is not necessarily ``yes.'' While the presence of higher derivatives clearly generates some additional issues to
deal with in the control design, the wave equation has its own peculiarities that one should not underestimate. For example,
controllability results for beams are valid on arbitrary short time intervals, whereas for strings such results hold only over time intervals that are lower
bounded in proportion to the ``wave propagation spee'' of the string (which physically corresponds to ``string tension''). Also, it is almost intuitively
evident that keeping a string from vibrating may not be easier than keeping a beam from vibrating. 

\begin{figure}[ht]
\centering
\includegraphics[width=8.0cm]{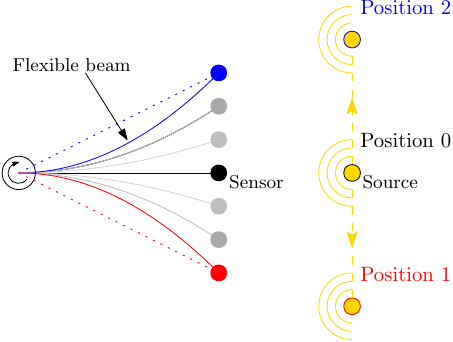}
\caption{\textcolor{black}{Light-source seeking with a coupled flexible beam.}} 
\label{fig:underlight}
\end{figure}

Fortunately both the \textit{backstepping design} and the \textit{trajectory generation problem} for the Euler–Bernoulli model were already studied in
\cite[Section~8.2 and Example~12.8]{krstic2008boundary}, respectively. As before, in the source seeking context, we can assume $\alpha(0,t)=\Theta$ and $\alpha(D,t)=\theta$, where $D$ would be the length of the flexible beam. The illustrated system in Fig.~\ref{fig:underlight} is composed of a flexible and inextensible Euler–Bernoulli beam clamped to a rotating actuator hub with a light sensor mass at its free-end. The control task aims to control the rotating actuator in order to drive the sensor to the moving light-source signal, while compensating the vibration generated during the beam displacement. Differently from the standard output regulation perspective, the paradigm here is the maximization of the signal perceived by the light sensor, while position regulation of the end-effector is also indirectly achieved. 

\section{Neuromuscular Electrical Stimulation}

\begin{figure}[ht]
\begin{center}
\includegraphics[width=8.2cm]{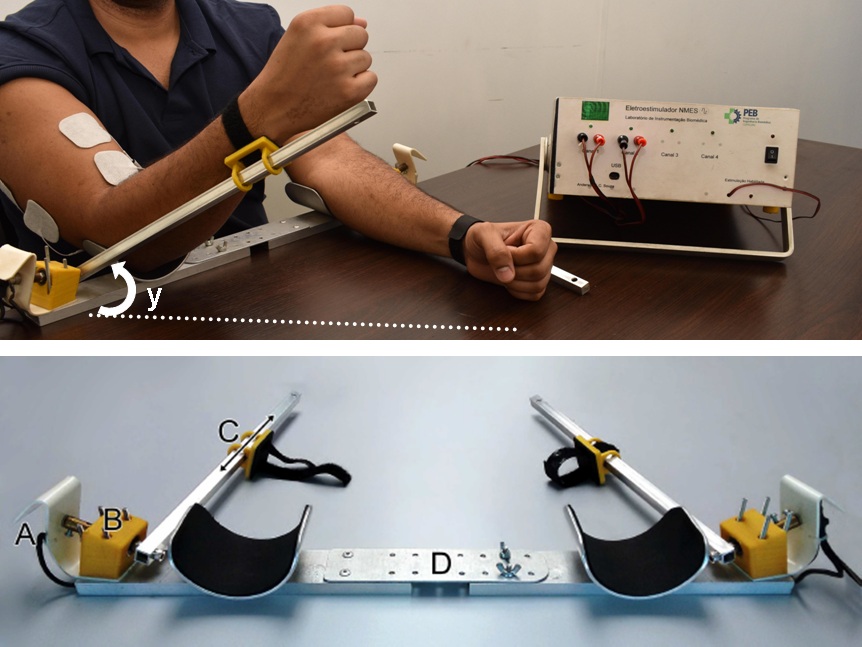}
\caption{Mechanical apparatus for NMES experimental tests. The point $A$ in the image indicates a goniometer
(simple potentiometer) linked to a steel axis $B$ allowing angular displacement readings. Letter $C$ shows that the wrist has an attachment with linear freedom of movement along the aluminum square rod, while $D$ points out that there is an adjustment for the lateral distance of the elbows. In the picture on the top, the controlled joint angle, denoted by $y$, and the NMES equipment are presented.}
\label{ch6.fig4_new}
\end{center}
\end{figure}

The recent manuscript \cite{PAZ:2019} on neuromuscular electrical stimulation (NMES) brings a promising application of ES under nonconstant delays. The authors have proposed a stochastic proportional-derivative-integral (PID) automatic tuner via ES and applied it to precise tracking of a flexion-extension reference for NMES and motor relearning for rehabilitation. Pictures of the mechanical apparatus as well as
the electrical stimulation device for NMES experimental tests are shown in Fig.~\ref{ch6.fig4_new}. The experimental results are innovative since, unlike the referred literature, stroke patients were recruited for the successful tests rather than only healthy subjects. Remarkably, the patients are using such a device at the Public Hospital Universit\'{a}rio Clementino Fraga in Rio de Janeiro, Brazil.

In general, NMES devices are applied to clinical work in an open-loop fashion and their parameters must be set at the beginning of the therapy, not facilitating the clinical practice. The levels of electrical stimulation follow pre-calibrated profiles, requiring the presence of a practitioner to modify the stimulation parameters. This requires protocols aiming to enhance muscle contraction together with the execution of the intended contractions. The downside of this procedure is that the device always returns the same portion of electrical assistance to the patient, if there is no therapist intervention. In addition, the open-loop devices are not prepared to promote a proper association by means of some feedback error between the subject's intended movement and the artificial activation provided by the NMES system.

\begin{figure*}[ht]
	\begin{center}
		\centering{	\includegraphics[width=13.5cm]{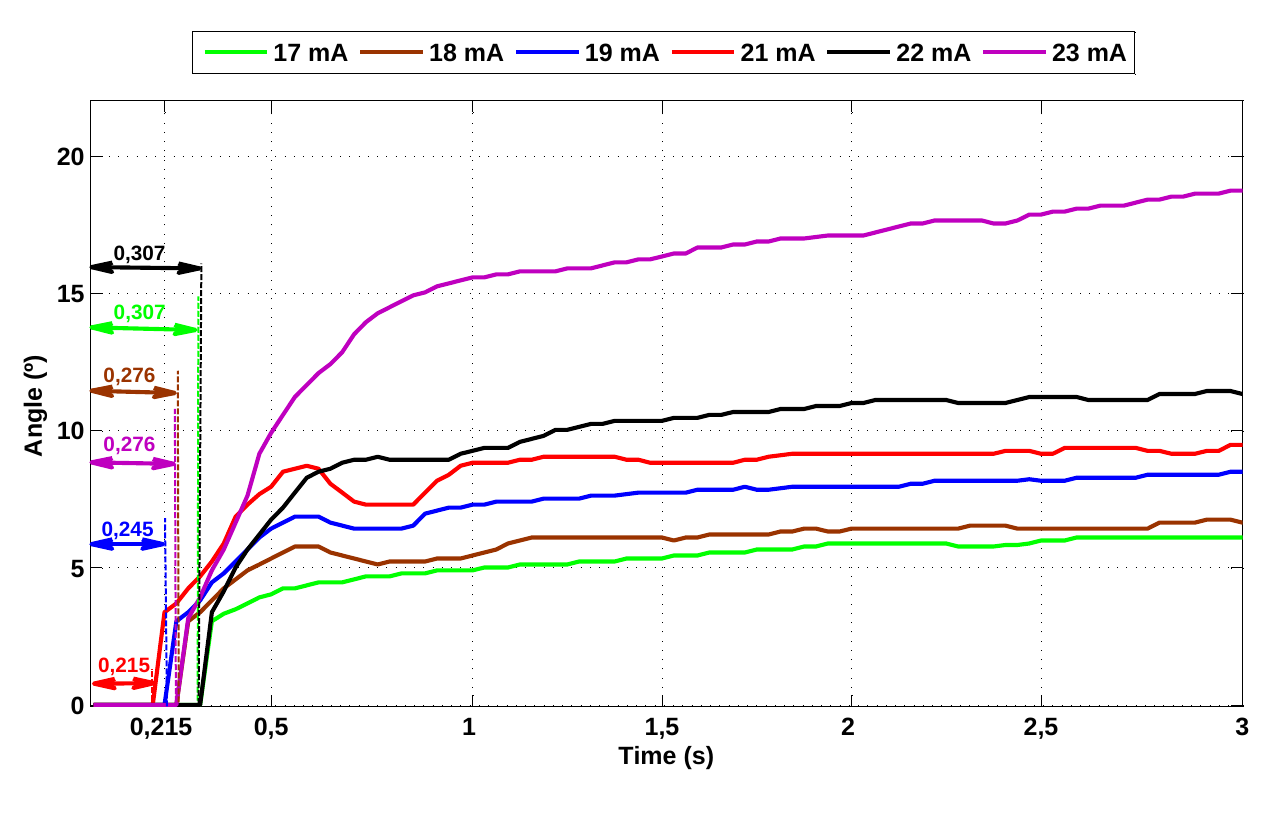}}
				\vspace{-0.5cm}
		\caption{\textcolor{black}{Open-loop responses for different step-current-inputs of a stroke patient with delays of order $300$ms. Extracted from \cite{PAZ:2019}.}}
		\label{delay_jorge}
	\end{center}
\end{figure*}

In this sense, closed-loop strategies are adequate to generate the NMES electrical current amplitudes based on the angular displacement (or the measurement of some other variable) related to the upper limbs. Although PID controllers have been explored in different engineering applications, the true limitation is that a PID controller is designed for linear systems, but the neuromuscular plant, which is being controlled, is nonlinear, time-varying and subject delays. For instance, it is possible to note a time delay in the subject responses shown in Fig.~\ref{delay_jorge}. The electromechanical-neuromuscular delay is in general of a time-varying nature and distinct for each subject. Moreover, it is worth mention that clinicians' knowledge of control systems is limited. Therefore, their expertise in tuning controllers is limited as well. Furthermore, in NMES applications, each patient is unique and requires a particular set of PID parameters. Since it may be difficult to find proper parameters for each patient, better procedures or a more intelligent-adaptive controller are indeed well motivated.

On the other hand, manual tuning is a time-consuming task and analytical methods are based on an exaggerated knowledge of the plant, requiring particular experimental validations to the identification of an acceptable plant model. However, a precise plant model in NMES is not known, and very long identification procedures are not desirable with the patients. This adverse environment of modeling inspires the application of adaptive-robust control methodologies and automatic tuning techniques.

In this context, the model-free PID tuner via multi-parameter stochastic ES in \cite{PAZ:2019} inspired by its earlier deterministic version \cite{KK:2006} is shown to be importantly fruitful. The proposed PID tuner eliminates the initial off-line tests with patients since the control gains are automatically computed in order to minimize a cost function according to the tracking error $e(t):=y(t)-r(t)$ between the elbow's angle of the patient's arm $y(t)$ and the reference trajectory $r(t)$. This research is highly successful since the stochastic algorithm provides faster (better transient) responses, which is perfect for a self-tuning, fatigue resistant control method for neuromuscular-based therapies. Moreover, the parameters of the stochastic ES are simpler to tune since the orthogonality assumption on the dither vector signals of multi-parameter deterministic ES impose additional obstacles in adjusting the frequencies of the sinusoidal perturbations. Finally, deterministic ES may restrict the region of convergence of the algorithm, and the adaptation using a periodic-deterministic perturbation for learning may be rather poor and unusual in some model-free optimization frameworks. Stochastic perturbations overcome those obstacles as well.

\begin{figure}[ht]
	\hspace{-0.5cm}
	\includegraphics[scale=0.31]{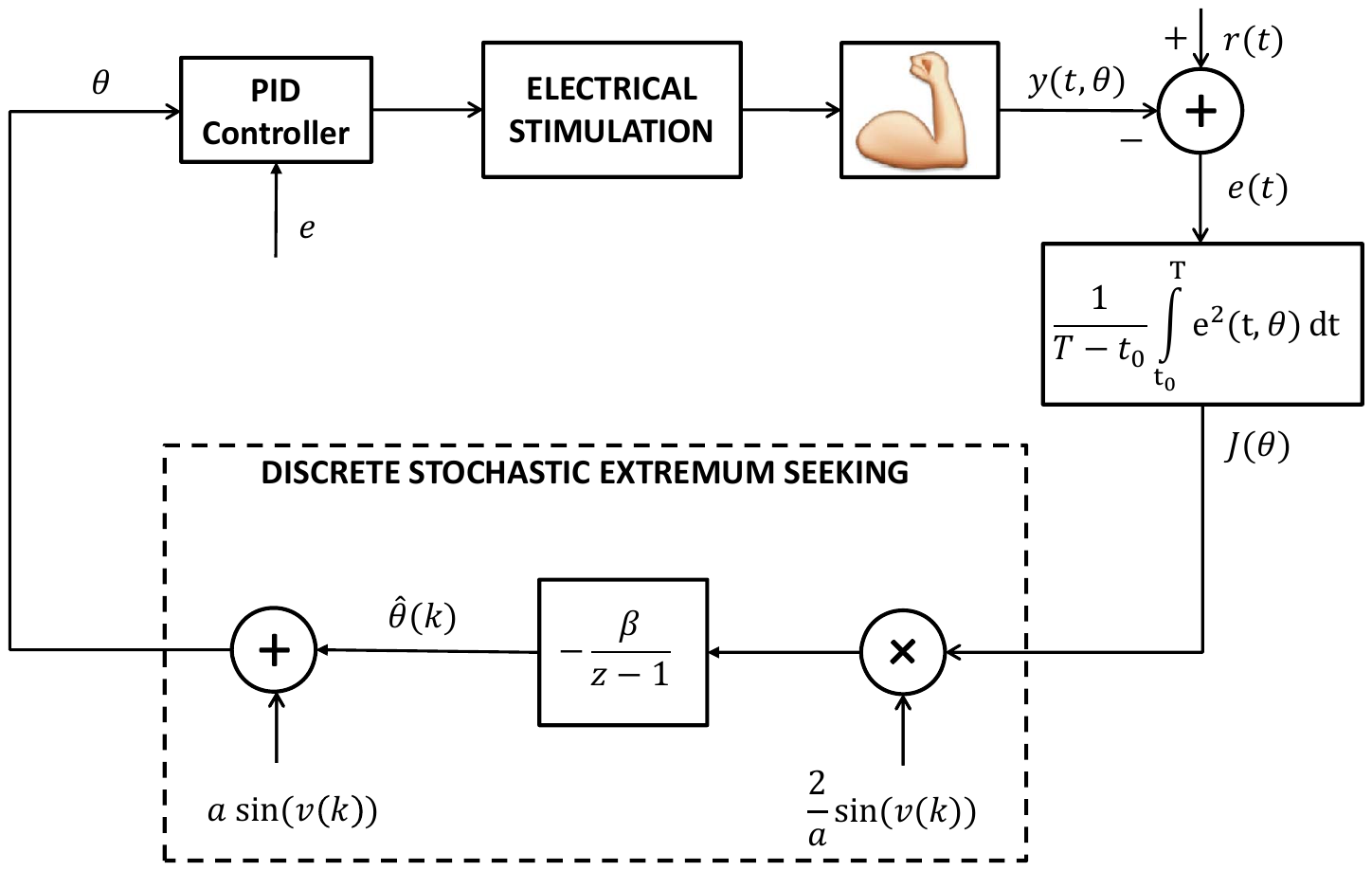}
	\vspace{-0.9cm}
	\caption{Block diagram of the closed-loop system for NMES using discrete-time stochastic ES, where $\nu(k)$ is the stochastic perturbation vector modeled by Gaussian white noise signals, and the cost function $J(\theta):=\frac{1}{T-t_0} \int_{t_0}^{T} e^2(t,\theta)dt$ is defined over a period between two time instants $t_0$ and $T$. The PID control law is $u(t) = K_p e(t) + K_i \int_{0}^{t} e(\tau)~ d\tau + K_d \frac{de(t)}{dt}$, with $\theta:=[K_p\,, K_i\,, K_d]^T$ being the PID parameters (gains) to be adapted by means of the ES algorithm.}
	\label{Block Diagram scale sem centering}
\end{figure}

\newpage
Specifically, ES minimizes a cost function which quantifies the performance of the PID controller and iteratively modifies the arguments 
of the cost function (the PID parameters) \textcolor{black}{so that its output reaches} a local minimum. According to the block diagram in Fig.~\ref{Block Diagram scale sem centering}, the time-domain implementation of the discrete-time stochastic ES
algorithm is given by: 
\begin{equation} \label{adaptation}
\hat{\theta}_\textrm{i}(k+1)=\hat{\theta}_\textrm{i}(k)-\beta \frac{2}{a} \sin(\nu_\textrm{i}(k)) J(\hat{\theta}(k)+a \sin(\nu(k)))\,,
\end{equation}
where $k$ is the discrete iteration number, the step size $\beta>0$ is sufficiently small,
the subscript $\textrm{i}\!=\!1, 2, 3$ indicates the $i$-th entry
of a vector, $\nu(k)=[\nu_1(k) \ \nu_2(k) \ \nu_3(k)]^T$ and $\sin(\nu(k))=[\sin(\nu_1(k)) \ \sin(\nu_2(k)) \ \sin(\nu_3(k))]^T$. The elements of the stochastic Gaussian perturbation vector $\nu(k)$ are sequentially and mutually independent such that $E\{\nu(k)\}=0$, $E\{\nu_\textrm{i}^2(k)\}=\sigma^2_\textrm{i}$ and $E\{\nu_\textrm{i}(k) \ \nu_\textrm{j}(k)\}=0$, $\forall \textrm{i} \neq \textrm{j}$, with $E\{\cdot\}$ denoting the 
expectation of a signal. In addition, it is assumed that the probability density function of the perturbation vector is
symmetric about its mean.

\begin{figure*}[ht]
	\begin{center}
		\includegraphics[width=11.5cm]{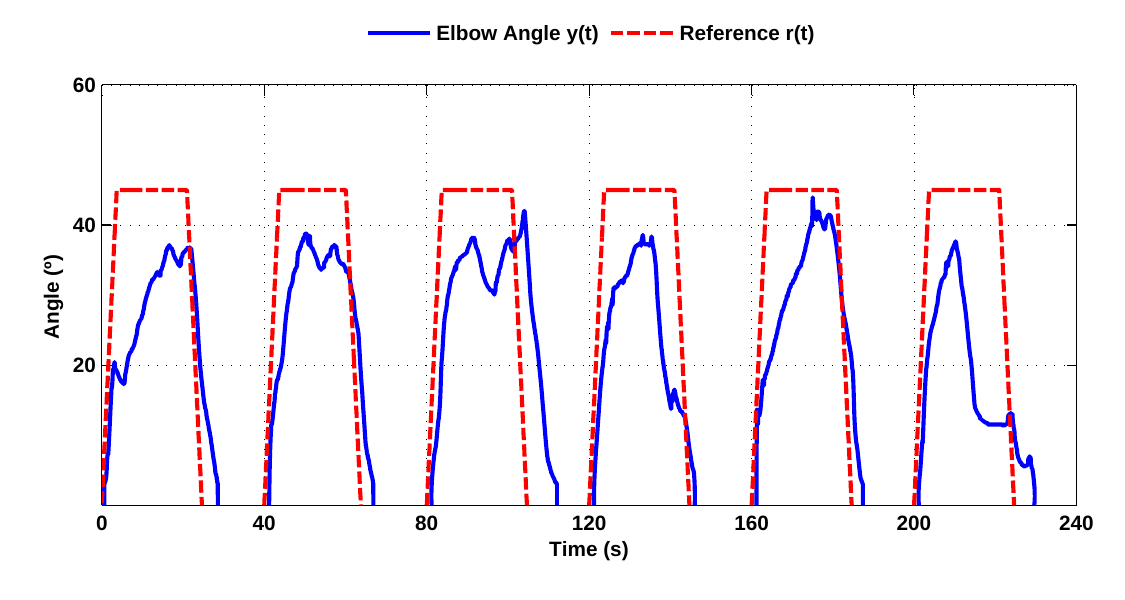}
	  \vspace{-0.5cm}
		\caption{\textcolor{black}{The graphic portrays the angular elbow joint movement performed by the stroke patient without the help of the proposed NMES controller. It can be seen that the stroke subject is not able to actively contract his arm to the final flexion position.}}\label{no_stim2}
	\end{center}
      \end{figure*}

\begin{figure*}[ht]
	\begin{center}
		\includegraphics[width=12.0cm]{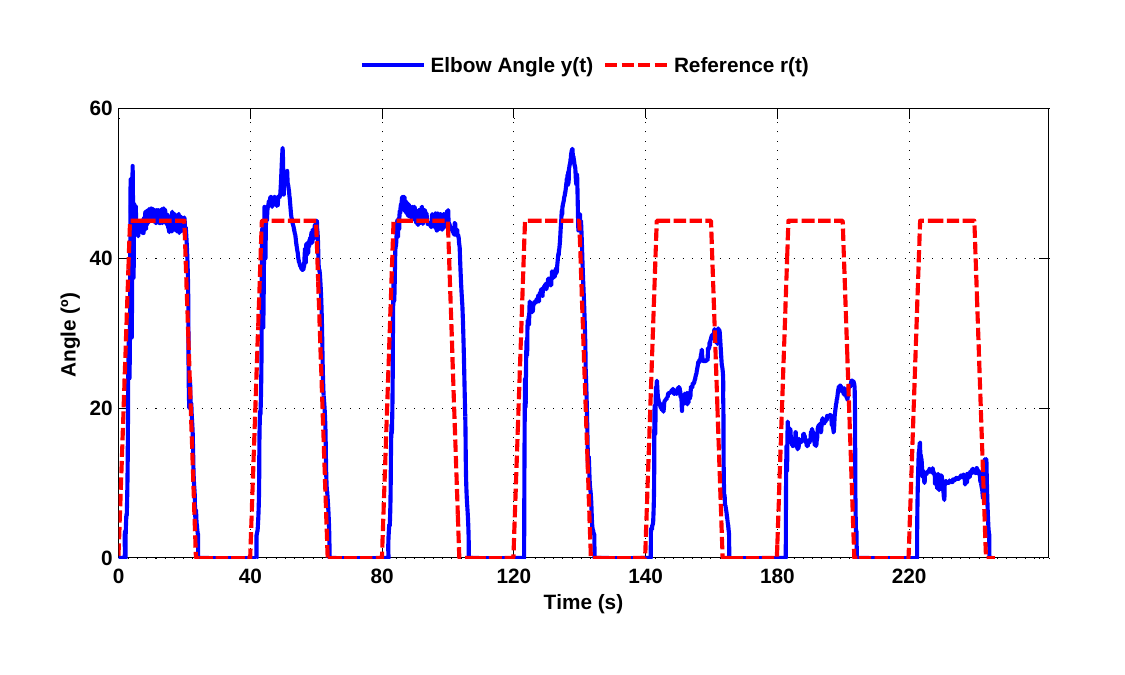}
		\vspace{-0.5cm}
		\caption{\textcolor{black}{Output responses for a stroke patient: PID with fixed gains
		($K_p=1\,, K_i=1\,, K_d=1$) not guaranteeing an acceptable trajectory tracking after $120$ seconds.}}\label{fixo3}
	\end{center}
\end{figure*}

\begin{figure*}[ht]
	\begin{center}
		\includegraphics[width=15.0cm]{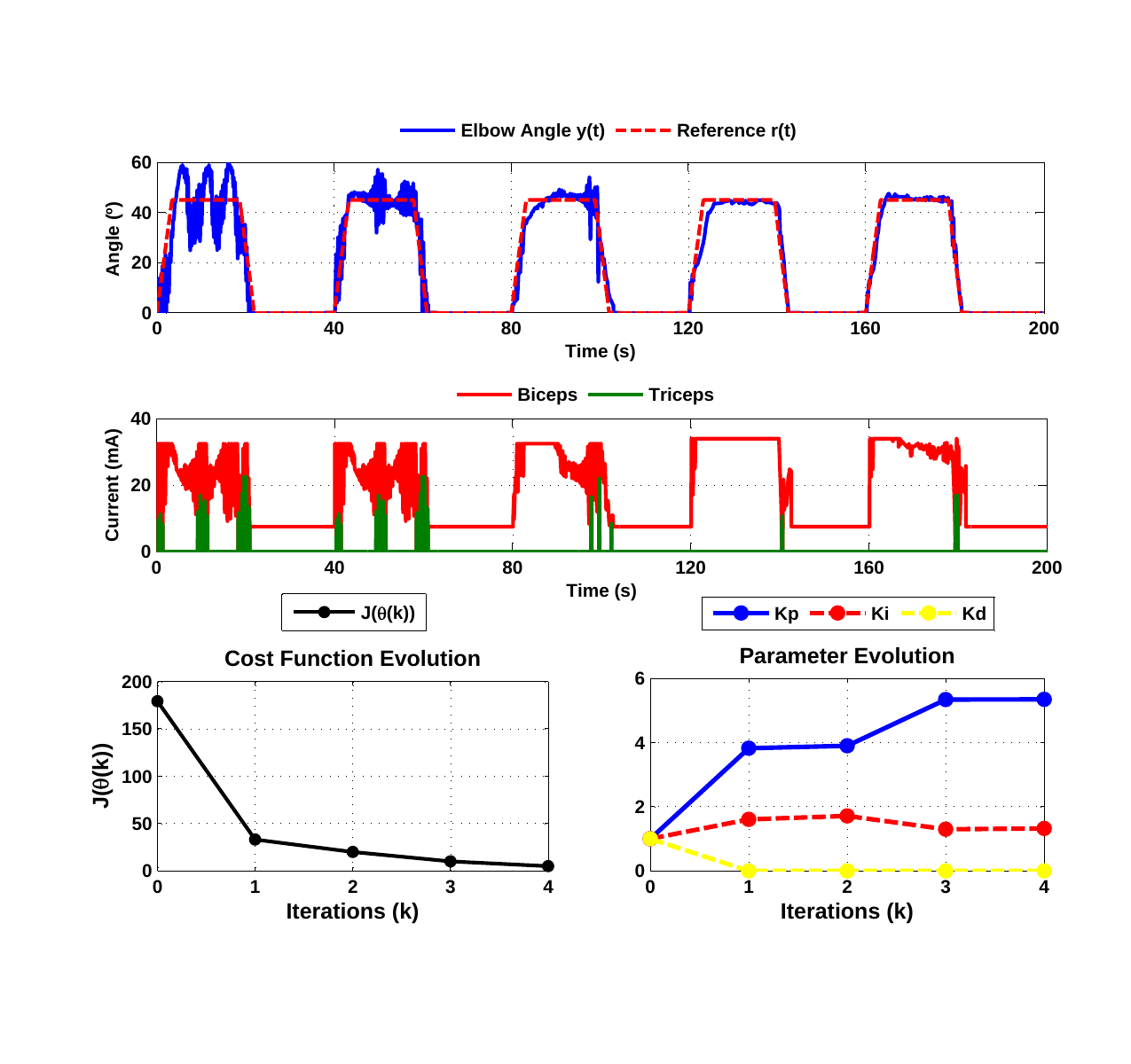}
		\vspace{-1.6cm}
		\caption{\textcolor{black}{Stochastic extremum-seeking based PID control for NMES of a stroke patient.}}\label{stocha1}
	\end{center}
\end{figure*}

Figs.~\ref{no_stim2} to \ref{stocha1} show the advantages for the closed-loop responses for a stroke patient, originally presented in \cite{PAZ:2019}. In the clinical scenario, Fig.~\ref{no_stim2} illustrates that even if the patients know the trajectory of the movement to be performed, they cannot execute it by themselves. Fig.~\ref{fixo3} also highlights that a fixed-gain PID scheme is not able to control adequately the stroke patient as the number of cycles/iterations increase, due to the time-varying nature of the neuromuscular system under delays. Indeed, PID control with fixed gains is not able to bring satisfactory results in long-running tests. Moreover, a unique fixed gain tuning is not applicable for different individuals. On the other hand, the ES adaptive approach for simple adaptation of PID controller parameters is model-free having the interesting ability of controlling on multiple subjects without tediously tuning the designer or practitioner. On the contrary, the response curves in Fig.~\ref{stocha1} ratify the improved behavior of the adaptive PID control scheme over a fixed-gains PID controller even in this adversarial scenario for NMES.
%

Although the delay discussion was not the focus on \cite{PAZ:2019}, it was evidenced there and in previous publications \cite{MDOD:2016,AKFS:2015,SGD:2011} that they may represent a significant challenge in NMES. It motivates the application of predictor feedback developed here for delay compensation in extremum seeking algorithms plus PID control or even other techniques \cite{AKFS:2015,AS:2017}.

\section{Conclusion}  \label{sec:concl}
This paper has generalized the ES results obtained in \cite{FKB:2012} to a wider class of infinite-dimensional systems governed by homogeneous and heterogeneous PDEs, rather than being restricted to noncooperative games governed by ODEs. We have introduced a non-model based approach via extremum seeking and boundary control to find, in a distributed way, the Nash equilibria of noncooperative games with unknown quadratic payoff functions and the players acting through PDE dynamics. A player could stably attain its Nash equilibrium by measuring only the value of its payoff (no other information about the game is needed).

We did not consider  
just one kind of PDE dynamics in the players' decision variables. The challenge of PDE dynamics from mixed PDE classes since games are decentralized and heterogeneous---each player has its own particular dynamics, possibly of a different nature.  

In this scenario, we were able to develop a result for heterogeneous games with distinct transport {and} diffusion (heat) PDEs to be simultaneously compensated---or homogeneous games with PDEs of different transport and diffusion coefficients---in the action paths of each player, where the players estimate only the diagonal entries of the Hessian matrix due to the players' own payoffs. We were able to dominate sufficiently small off-diagonal terms using a small-gain argument for the average system. Convergence to a small neighborhood of the Nash equilibrium is achieved, even in the presence of transport-heat PDEs.  


The introduction to boundary control for noncooperative games was presented in this paper for both hyperbolic-transport and parabolic-heat types of PDEs. There is no strong reason why the  exposition provided could not have been conducted on some of the other classes of PDEs. However, transport and heat PDEs are particularly convenient because they are at the same time sufficiently simple and sufficiently general to serve as a design template using which the reader can pursue extensions of the Nash equilibrium seeking design to other classes of PDEs \cite{krstic2008boundary,TAC:2020} and its integration to other real-world applications described by a game-theoretic framework \cite{ZB:2024,LPZB:2022,ZB:2015,ICBS:2001}.


We have explored a variety of potential engineering applications where these advanced strategies could be effectively implemented. These included areas such as fluid dynamics, structural mechanics as well as advanced energy and transportation systems. The aim was to demonstrate how the theoretical foundations can lead to innovative solutions in real-world engineering challenges, paving the way for future research.

\newpage
\appendices

\section{Proof of Theorem~\ref{ch16.theorem.16.1.vacina}} \label{appendix_A}
The proof of Theorem~\ref{ch16.theorem.16.1.vacina} follows steps similar to those employed to prove the results about extremum seeking under diffusion PDEs in \cite[Theorem~1]{FKKO:2018} or even \cite[Theorem~1]{Tiago} for pure delays (transport PDEs). In this sense, we will simply point out the main differences for the case of games (not classical extremum seeking), instead of giving a full independent proof.

While in \cite[Theorem~1]{FKKO:2018}  and \cite[Theorem~1]{Tiago} it was possible to prove local exponential stability of the average closed-loop system using a Lyapunov functional, 
a different approach is adopted here for the Nash equilibrium seeking in noncooperative games. We will show that it is possible to guarantee the local exponential stability for the average closed-loop system (\ref{saco1M_NC_CP14})--(\ref{saco3M_NC_CP14}) by means of a small-gain analysis.


First, consider the equivalent parabolic PDE-ODE representation  (\ref{saco1M_NC_CP14})--(\ref{saco3M_NC_CP14}) rewritten for each Player $i$, $\forall i\in\{1\,,\ldots\,,N\}$:
\begin{eqnarray}
\dot{\bar{G}}^{\rm{av}}_i(t) &=& H_{ii}^{i}k_i \bar{G}^{\rm{av}}_i(t)+\epsilon H_{ii}^{i}k_i \phi_i^{\rm{av}}(1,t)\,, \label{saco1M_NC_i_CP14}\\
\partial_t u^{\rm{av}}_i(x,t)&=&D_i^{-2}\partial_{xx} u^{\rm{av}}_i(x,t)\,, \quad x\in (0\,,1)\,, \label{saco2M_NC_i_CP14}\\
\partial_x u^{\rm{av}}_i(0,t)&=&0\,, \\
u^{\rm{av}}_i(1\,,t)&=& k_i \bar{G}^{\rm{av}}_i(t)+\epsilon k_i \phi^{\rm{av}}_i(1,t) \label{saco3M_NC_i_CP14}\,.
\end{eqnarray}
%

Hence, the average closed-loop system (\ref{saco1M_NC_i_CP14})--(\ref{saco3M_NC_i_CP14}) satisfies all the assumptions \textbf{(A1)} to \textbf{(A7)} of the Small-Gain Theorem \cite[Theorem~8.2, p.~205]{KK:2018} for the parabolic PDE-ODE loops with $p(z)\!=\!1$, $r(z)\!=\!D_i^{2}$,
$q(z)\!=\!0$, $F(\bar{G}^{\rm{av}}_i\,,u^{\rm{av}}\,, 0)\!=\!H_{ii}^{i}k_i \bar{G}^{\rm{av}}_i\!+\!\epsilon H_{ii}^{i}k_i \phi_i^{\rm{av}}(1)$, $g(x\,,\bar{G}^{\rm{av}}_i\,,u^{\rm{av}})\!=\!0$, $f(x,t)\!=\!0$, $\varphi_0(0\,,u_i^{\rm{av}}\,,\bar{G}^{\rm{av}}_i)\!=\!b_1 u_i^{\rm{av}}(0,t)$, $b_1 \!<\! 0$, $b_2\!=\!1$, 
$\varphi_1(0\,,u^{\rm{av}}\,,\bar{G}^{\rm{av}}_i)\!=\!k_i \bar{G}^{\rm{av}}_i\!+\!\epsilon k_i \phi^{\rm{av}}_i(1)$, $a_1=1$, $a_2=0$,  $L=\max(|H_{ii}^{i}|k_i\,, \frac{1}{\sqrt{3}}\epsilon |H_{ii}^{i}|k_i k_H D_j^2)$, $K_0=1$, $B_0=C_0=0$, $\gamma_0$ is of order
$\mathcal{O}(\epsilon)$, $K_1=\frac{1}{\sqrt{3}} \epsilon k_i k_H D_j^2$, $B_1=k_i$,
$C_1=0$, $\gamma_1$ is of order
$\mathcal{O}(1)$, $K_2=B_2=0$ and $i\neq j$. The constant $k_H>0$ is defined in the next just after inequality (\ref{Elizabeth_spam_folder}). Assumption \textbf{(A6)} of \cite[Theorem~8.2, p.~205]{KK:2018} holds with $M=1$, 
$\gamma_3=\frac{1}{\sqrt{3}}\epsilon |H_{ii}^{i}| k_i k_H D_j^2$, $\sigma=|H_{ii}^{i}|k_i$ as it can be readily verified by means of the variations of constants formula 
\begin{align}
\bar{G}^{\rm{av}}_i(t)&=\exp(-|H_{ii}^{i}|k_i t)\bar{G}^{\rm{av}}_i(0) \nonumber \\
&\quad+\int_{0}^{1} \exp(-|H_{ii}^{i}|k_i(t+s))\epsilon H_{ii}^{i}k_i \phi_i^{\rm{av}}(1,s)ds\,, \nonumber
\end{align}
and from the application of the Cauchy-Schwarz inequality to the term $\phi^{\rm{av}}(1,t)$ in equation (\ref{eq:phii_CP14}):
\begin{small}
\begin{eqnarray} \label{Elizabeth_spam_folder}
\phi_i^{\rm{av}}(1,t) \!\!\!\!&\leq&\!\!\!\! \sum_{j\neq i}^{N} |H_{ij}^{i}| D_j^2
\left(\!\int_{0}^{1} \!\!\!(1-\tau)^2 d\tau  \right)^{\frac{1}{2}}\!\!\!\!\times\!
\left(\!\int_{0}^{1} \!\![u^{\rm{av}}_j(\xi,t)]^{2} d\xi \right)^{\frac{1}{2}}  \nonumber \\
\!\!\!\!&\leq&\!\!\!\! \frac{1}{\sqrt{3}}  k_H \sum_{j\neq i}^{N}  D_j^2
\left(\!\int_{0}^{1} [u^{\rm{av}}_j(\xi,t)]^2 d\xi \right)^{\frac{1}{2}}\,, 
\end{eqnarray}
\end{small}
$\!\!$since $|H_{ij}^{i}| \!\!<\!\! k_H \!\!<\!\! \frac{1}{\epsilon}|H_{ii}^{i}|$ according to \textcolor{black}{Assumptions~\ref{ch.13.Assumption 1.} and \ref{ch13.Assumption 2.}},
where $k_H$ is a positive constant of order $\mathcal{O}(1)$.  
It follows that the small-gain condition \cite[Inequality (8.3.24)]{KK:2018}: 
\begin{eqnarray} \label{cuzao}
&\max(\gamma_0K_0\,,\gamma_1K_1)+\sigma^{-1}K_2<1\,, \nonumber \\
&\gamma_3 \max(\gamma_0B_0\,,\gamma_1B_1)+\gamma_3\sigma^{-1}B_2<1
\end{eqnarray}
holds provided $0<\epsilon<1$ is sufficiently small. 
Therefore, if 
such a small-gain condition holds, then \cite[Theorem~8.2, p.~205]{KK:2018} allows us to conclude that there exist constants $\delta\,,\Delta>0$ such that for every $u_0^{\rm{av}} \in C^{0}([0\,,1])$, $\bar{G}^{\rm{av}}_{0} \in \mathbb{R}^{n}$, the unique generalized solution of this initial-boundary value problem, with $u^{\rm{av}}(x,0)=u^{\rm{av}}_0$ and $\bar{G}^{\rm{av}}(0)=\bar{G}^{\rm{av}}_{0}$, satisfies the following estimate:  
%
\begin{eqnarray} \label{bocao_CP14}
|\bar{G}^{\rm{av}}(t)|\!+\!\|u^{\rm{av}}(t)\|_{\infty}\!\leq\! 
\Delta(|\bar{G}^{\rm{av}}_{0}|\!+\!\|u^{\rm{av}}_0\|_{\infty})\exp(-\delta t)\,.
\end{eqnarray}
%
Therefore, we conclude that the origin of the average closed-loop system (\ref{saco1M_NC_CP14})--(\ref{saco3M_NC_CP14}) is exponentially stable under the assumption 
of $0\!<\!\epsilon\!<\!1$ being sufficiently small. Then, from (\ref{eq:hatGi_CP14}) and (\ref{transformacao_cu_NC_CP14}), we conclude the same results in the norm
\begin{eqnarray}\label{periodic_solution_NCG_infinite_dimension_CP14}
\left(\sum_{i=1}^{N}\left[\vartheta_{i}^{\rm{av}}(t)\right]^2  + \int_{0}^{D_i} \left[u_{i}^{\rm{av}}(x,t)\right]^2dx \right)^{1/2} \end{eqnarray}
since $H$ is non-singular, \textit{i.e.}, $|\vartheta_{i}^{\rm{av}}(t)|\leq |H^{-1}||\hat{G}^{\rm{av}}(t)|$. 

\textcolor{black}{As developed in \cite[Theorem~1]{FKKO:2018}, the next steps of the proof would be the application of the local averaging theory for infinite dimensional systems in \textcolor{black}{\cite[Sec.~2]{HL:1990}} (see also ``Averaging Theorem for General Infinite-Dimensional Systems''), showing that the periodic solutions satisfy (\ref{periodic_solution_NCG_CP14}) for $\omega$ sufficiently large, and then the conclusion of the attractiveness of the Nash Equilibrium $\Theta^*$ according to (\ref{limsup1_NCG_CP14}).} The final residual set for the error $\theta(t)\!-\!\theta^*$ in
(\ref{limsup2_CP14}) depends on $|a|e^{\max(D_i)\sqrt{\frac{\omega}{2}}}$ due to the amplitude of  $S_i(t)$ in (\ref{eq:Si_CP14}). 

\section{Proof of Theorem~\ref{scheme2}}
\label{appendix_B}
%

The proof of Theorem~\ref{scheme2} follows the same structure of that of Theorem~\ref{ch16.theorem.16.1.vacina}, but now considering the equivalent hyperbolic PDE-ODE representation  
for each Player $i$, $\forall i\in\{1\,,\ldots\,,N\}$:
\begin{eqnarray}
\dot{\bar{G}}^{\rm{av}}_i(t) &=& H_{ii}^{i}k_i \bar{G}^{\rm{av}}_i(t)+\epsilon H_{ii}^{i}k_i \phi_i^{\rm{av}}(1,t)\,, \label{saco1M_NC_i}\\
\partial_t u^{\rm{av}}_i(x,t)&=&D_i^{-1}\partial_x u^{\rm{av}}_i(x,t)\,, \quad x\in ~ \textcolor{black}{(0\,,1)}\,, \label{saco2M_NC_i}\\
u^{\rm{av}}_i(1\,,t)&=& k_i \bar{G}^{\rm{av}}_i(t)+\epsilon k_i \phi^{\rm{av}}_i(1,t) \label{saco3M_NC_i}\,,
\end{eqnarray}
where $H_{ii}^{i}<0$, $k_i>0$, $0<\epsilon<1$,  $D_i^{-1}>0$, and
\begin{align}
\phi_{i}(1,t)&:=-\sum_{j\neq i}H_{ij}^{i}\int_{0}^{1} D_ju_j(\xi,t) d\xi\,.  \label{eq:phii} 
\end{align}
In this case, the average closed-loop system (\ref{saco1M_NC_i})--(\ref{saco3M_NC_i}) satisfies both assumptions \textbf{(H1)} and \textbf{(H2)} of the Small-Gain Theorem \cite[Theorem~8.1, p.~198]{KK:2018} for the hyperbolic PDE-ODE loops rather than \cite[Theorem~8.2, p.~205]{KK:2018} employed in the proof of Theorem~\ref{ch16.theorem.16.1.vacina}, this latter oriented for parabolic PDE-ODE interconnections---see ``Small-Gain Theorem for ODE and Hyperbolic PDE Loops''. 

\textcolor{black}{It follows that the small-gain condition in \cite[Theorem~8.1, p.~198]{KK:2018} holds provided $0\!<\!\epsilon\!<\!1$ is sufficiently small.}   
Therefore, if 
such a small-gain condition holds, then \cite[Theorem~8.1, p.~198]{KK:2018} allows us to conclude that 
%
the origin of the average closed-loop system 
(\ref{saco1M_NC_i})--(\ref{saco3M_NC_i}) is exponentially stable. 
Then, from (\ref{eq:hatGAv}) and the following transformation \cite{a14} 
%
\begin{align} 
\bar{G}_{i}^{\rm{av}}(t)&=\hat{G}_{i}^{\rm{av}}(t)+ \sum_{j=1}^{N}\epsilon_{ij}^{i}H_{ij}^{i}\int_{t-D_j}^{t} U_j^{\rm{av}}(\tau) d\tau \nonumber\\
&= \hat{G}_{i}^{\rm{av}}(t) + \sum_{j=1}^{N}\epsilon_{ij}^{i}H_{ij}^{i}\int_{0}^{D_j} u_j^{\rm{av}}(\xi,t) d\xi\,, \label{transformacao_cu_NC}
\end{align}
we can conclude the same results in the norm
\begin{eqnarray}\label{periodic_solution_NCG_infinite_dimension}
\left(\sum_{i=1}^{N}\left[\tilde{\theta}_{i}^{\rm{av}}(t-D_i)\right]^2  + \int_{0}^{D_i} \left[u_{i}^{\rm{av}}(\tau)\right]^2d\tau \right)^{1/2}\,.
\end{eqnarray}
%

The application of the local averaging theory for functional differential equations in \cite{HL:1990} (see also ``Averaging Theorem for Functional Differential Equations'') shows that the periodic solutions satisfy inequality 
(\ref{periodic_solution_NCG}) and leads to the conclusion of the attractiveness of the Nash Equilibrium $\theta^*$ according to (\ref{limsup1_NCG}). 

\section{Proof of Theorem~\ref{duoply_heterogeneous_queirozBday}}
\label{appendix_C}
%
%
%
First, consider the equivalent hyperbolic/parabolic PDE-ODE representation  (\ref{ch15.saco1M_NC})--(\ref{ch15.saco3M_NC}) rewritten for each Player
$P_i$, $i\in\{1\,,2\}$:
\begin{eqnarray}
\dot{\bar{G}}^{\rm{av}}_i(t) &=& H_{ii}^{i}k_i \bar{G}^{\rm{av}}_i(t)+\epsilon H_{ii}^{i}k_i \phi_i^{\rm{av}}(1,t)\,, \label{ch15.saco1M_NC_i}\\
\partial_t u^{\rm{av}}_1(x,t)&=&D_1^{-1}\partial_{x} u^{\rm{av}}_1(x,t)\,, \quad ~x\in (0\,,1)\,, \label{ch15.saco2M_NC_i_delay}\\
\partial_t u^{\rm{av}}_2(x,t)&=&D_2^{-2}\partial_{xx} u^{\rm{av}}_2(x,t)\,, \!\quad x\in (0\,,1)\,, \label{ch15.saco2M_NC_i}\\
\partial_x u^{\rm{av}}_i(0,t)&=&0\,, \label{ch15.saco4M_NC_i} \\
u^{\rm{av}}_i(1\,,t)&=& k_i \bar{G}^{\rm{av}}_i(t)+\epsilon k_i \phi^{\rm{av}}_i(1,t) \label{ch15.saco3M_NC_i}\,.
\end{eqnarray}
%
\sloppy
For Player $P_1$, the average closed-loop system, given by (\ref{ch15.saco1M_NC_i})--(\ref{ch15.saco2M_NC_i_delay}) and (\ref{ch15.saco4M_NC_i})--(\ref{ch15.saco3M_NC_i}),
satisfies both assumptions \textbf{(H1)} and \textbf{(H2)} of the Small-Gain Theorem \cite[Theorem~8.1, p.~198]{KK:2018} for the hyperbolic PDE-ODE loop with $n\!=\!1$, $x(t)\!=\!\bar{G}_1^{\rm{av}}(t)$, $F(\bar{G}_1^{\rm{av}}(t)\,,u_1^{\rm{av}}(x,t)\,, v(t))\!=\!H_{11}^{1}k_1 \bar{G}_1^{\rm{av}}(t)\!+\!\epsilon v(t)$ with $v(t)\!=\! H_{11}^{1}k_1 \phi_1^{\rm{av}}(1,t)$, $c\!=\!1/D_1$, $a(x)\!=\!0$, $g(x\,,\bar{G}^{\rm{av}}_1(t)\,,u_1^{\rm{av}}(x,t))\!=\!0$, $f(x,t)\!=\!0$, $\varphi(d(t)\,,u_1^{\rm{av}}(x,t)\,,\bar{G}^{\rm{av}}_1(t))\!=\!k_1 \bar{G}^{\rm{av}}_1(t)\!+\!\epsilon d(t)$ where $d(t)\!=\! k_1 \phi^{\rm{av}}_1(1,t)$, $N\!=\!k_1$, $L\!=\!|H_{11}^{1}|k_1$,  $B=0$, $\gamma_2=k_1$, $b_2=\epsilon$ and $A=\gamma_1=0$. 

Notice that Assumption \textbf{(H1)} holds with $M=1$, 
$\gamma_3\!=\!0$, $\sigma=|H_{11}^{1}|k_1$ and $b_3>0$ being an appropriate constant of order $\mathcal{O}(\epsilon)$ as it can be readily verified by means of the variation-of-constants formula (for $i=1$)
\begin{align} \label{ch15.cdc2021denovo}
\bar{G}^{\rm{av}}_i(t)&=\exp(-|H_{ii}^{i}|k_i t)\bar{G}^{\rm{av}}_i(0) \nonumber \\
&\quad+\int_{0}^{1} \exp(-|H_{ii}^{i}|k_i(t+s))\epsilon H_{ii}^{i}k_i \phi_i^{\rm{av}}(1,s)ds\,,
\end{align}
and from the application of the Cauchy-Schwarz inequality to the term $\phi_1^{\rm{av}}(1,t)$ in equation (\ref{ch15.eq:phii}): 
\begin{small}
\begin{eqnarray}
\phi_1^{\rm{av}}(1,t) \!\!\!\!&\leq&\!\!\!\! |H_{12}^{1}| D_2^2
\left(\!\int_{0}^{1} \!\!\!(1-\tau)^2 d\tau  \right)^{\frac{1}{2}}\!\!\!\!\times\!
\left(\!\int_{0}^{1} \!\![u^{\rm{av}}_2(\xi,t)]^{2} d\xi \right)^{\frac{1}{2}}  \nonumber \\
&\leq& \frac{1}{\sqrt{3}}  k_H D_2^2
\left(\!\int_{0}^{1} [u^{\rm{av}}_2(\xi,t)]^2 d\xi \right)^{\frac{1}{2}}\,, \label{ch15.cdc2021quesacodecu}
\end{eqnarray}
\end{small}
$\!\!$since $|H_{ij}^{i}| \!<\! k_H \!<\! \frac{1}{\epsilon}|H_{ii}^{i}|$ according to Assumption~\ref{ch17.Assumption 1.nene},
where $k_H$ is a positive constant of order $\mathcal{O}(1)$.  
%
Therefore, if 
$0<\epsilon<1$ is sufficiently small, 
then \cite[Theorem~8.1, p.~198]{KK:2018} allows us to conclude that there exist constants $\delta_1\,,\Delta_1>0$ such that for every $u_{1,0}^{\rm{av}} \in C^{0}([0\,,1])$, $\bar{G}^{\rm{av}}_{1,0} \in \mathbb{R}$, the unique generalized solution of this initial-boundary value problem, with $u_1^{\rm{av}}(x,0)=u^{\rm{av}}_{1,0}$ and $\bar{G}^{\rm{av}}_{1}(0)=\bar{G}^{\rm{av}}_{1,0}$, satisfies the following estimate, $\forall t \geq 0$:  
%
\begin{eqnarray} \label{ch15.bocao_delay}
|\bar{G}^{\rm{av}}_{1}(t)|\!+\!\|u_1^{\rm{av}}(t)\|_{\infty} \!\!\!\!&\leq&\!\!\!\!
\Delta_1(|\bar{G}^{\rm{av}}_{1,0}|+\|u^{\rm{av}}_{1,0}\|_{\infty})\exp(-\delta_1 t) + \nonumber \\
&+& \bar{\gamma}_1 \epsilon \max_{0\leq s \leq t} \left(\|u_2^{\rm{av}}(s)\|_{\infty}\right)\,,
\end{eqnarray}
%
for some adequate constant $\bar{\gamma}_1>\max(|H_{11}^1|k_1\,,k_1)$.

For Player $P_2$, the average closed-loop system (\ref{ch15.saco1M_NC_i}) and (\ref{ch15.saco2M_NC_i})--(\ref{ch15.saco3M_NC_i}) satisfies all the assumptions \textbf{(A1)} to \textbf{(A7)} of the Small-Gain Theorem \cite[Theorem~8.2, p.~205]{KK:2018} for the parabolic PDE-ODE loop with $n\!=\!1$, $p(x)\!=\!1$, $r(x)\!=\!D_2^{2}$,
$q(x)\!=\!0$, $F(\bar{G}^{\rm{av}}_2(t)\,,u_2^{\rm{av}}(x,t)\,,v(t))=H_{22}^{2}k_2 \bar{G}^{\rm{av}}_2(t)+\epsilon v(t)$,
$v(t)=H_{22}^{2}k_2 \phi_2^{\rm{av}}(1,t)$, $g(x\,,\bar{G}^{\rm{av}}_2(t)\,,u_2^{\rm{av}}(x,t))=0$, $f(x,t)=0$, 
$\varphi_0(d(t)\,,u_2^{\rm{av}}(x,t)\,,\bar{G}^{\rm{av}}_2(t))\!=\!
b_1 u_2^{\rm{av}}(0,t)$, $b_1 \!<\! 0$, $b_2\!=\!1$, 
$\varphi_1(d(t)\,,u_2^{\rm{av}}(x,t)\,,\bar{G}^{\rm{av}}_2(t))\!=\!
k_2 \bar{G}^{\rm{av}}_2(t)\!+\!\epsilon d(t)$, $d(t)=k_2 \phi^{\rm{av}}_2(1,t)$,
$a_1=1$, $a_2=0$, 
%
%
$L=|H_{22}^{2}|k_2$, 
%
%
$K_0=|b_1|$, $B_0=C_0=0$, $\gamma_0$ is of order
$\mathcal{O}(1)$, 
%
$K_1=0$,  
%
$B_1=k_2$, $C_1=\epsilon$, $\gamma_1$ is of order
$\mathcal{O}(1)$ and $K_2=B_2=0$.  
Assumption \textbf{(A6)} of \cite[Theorem~8.2, p.~205]{KK:2018} holds with $M=1$, 
%
%
$\gamma_3\!=\!0$, $\sigma\!=\!|H_{22}^{2}|k_2$ and $b_3\!>\!0$ being of order $\mathcal{O}(\epsilon)$, as it can be readily verified by means of the variation-of-constants formula (\ref{ch15.cdc2021denovo}), with $i\!=\!2$,     
%
and from the application of the Cauchy-Schwarz inequality  
to the term $\phi_2^{\rm{av}}(1,t)$ in equation (\ref{ch15.eq:phii}):
\begin{small}
\begin{eqnarray}
\phi_2^{\rm{av}}(1,t) &\leq& |H_{21}^{2}|
\left(\int_{0}^{1} D_1^2 d\tau \right)^{\frac{1}{2}}\times
\left(\int_{0}^{1} [u^{\rm{av}}_1(\xi,t)]^{2} d\xi \right)^{\frac{1}{2}}  \nonumber \\
&\leq& k_H D_1
\left(\int_{0}^{1} [u^{\rm{av}}_1(\xi,t)]^2 d\xi \right)^{\frac{1}{2}}\,, 
\end{eqnarray}
\end{small}
$\!\!\!$with the same $k_H\!>\!0$ defined just after (\ref{ch15.cdc2021quesacodecu}). 
%
%
Hence, it follows that the small-gain condition in \cite[Inequality (8.3.24)]{KK:2018}
holds provided $0\!<\!\epsilon\!<\!1$ is sufficiently small. 
Thus, 
\cite[Theorem~8.2, p.~205]{KK:2018} allows us to conclude that there exist constants $\delta_2\,,\Delta_2\!>\!0$ such that 
%
%
\begin{eqnarray} \label{ch15.bocao}
|\bar{G}^{\rm{av}}_{2}(t)|\!+\!\|u_2^{\rm{av}}(t)\|_{\infty} \!\!\!\!&\leq&\!\!\!\!
\Delta_2(|\bar{G}^{\rm{av}}_{2,0}|\!+\!\|u^{\rm{av}}_{2,0}\|_{\infty})\exp(-\delta_2 t) + \nonumber \\
\!\!\!\!&+&\!\!\!\! \bar{\gamma}_2 \epsilon \max_{0\leq s \leq t} \left(\|u_1^{\rm{av}}(s)\|_{\infty}\right)\,,
\end{eqnarray}
%
$\forall t \geq 0$, for some adequate constant $\bar{\gamma}_2>\max(|H_{22}^2|k_2\,,k_2)$, $u_2^{\rm{av}}(x,0)=u^{\rm{av}}_{2,0}$ and
$\bar{G}^{\rm{av}}_{2}(0)=\bar{G}^{\rm{av}}_{2,0}$.  
Since inequalities (\ref{ch15.bocao_delay}) and (\ref{ch15.bocao}) are similar to those found in \cite[Theorem~11.2, p.~269]{KK:2018}---see inequalities (11.2.23) and (11.2.24)---we can finally invoke \cite[Theorem~11.5, p.~277]{KK:2018}, under the condition of $0<\epsilon<1$ sufficiently small, to conclude 
%
\begin{eqnarray} \label{ch15.bocao_final}
\!\!\!\!\!\!\!\!\!\!\!\!\!\!\!\!\!\!&&|\bar{G}^{\rm{av}}_{1}(t)|\!+\!|\bar{G}^{\rm{av}}_{2}(t)|\!+\!\|u_1^{\rm{av}}(t)\|_{\infty}\!+\!\|u_2^{\rm{av}}(t)\|_{\infty} \!\leq\! \nonumber \\
\!\!\!\!\!\!\!\!\!\!\!\!\!\!\!\!\!\!&&\Delta(|\bar{G}^{\rm{av}}_{1,0}|\!+\!|\bar{G}^{\rm{av}}_{2,0}|\!+\!\|u^{\rm{av}}_{1,0}\|_{\infty}\!+\!\|u^{\rm{av}}_{2,0}\|_{\infty})
\exp(-\delta t)\,, 
\end{eqnarray}
%
$\forall t \!\geq\! 0$, for some $\delta>0$ and $\Delta>0$. 

Therefore, we conclude that the origin of the average closed-loop system (\ref{ch15.saco1M_NC})--(\ref{ch15.saco3M_NC}) is exponentially stable under the assumption 
of $0\!<\!\epsilon\!<\!1$ being sufficiently small. Then, from (\ref{ch15.eq:hatGAv}), (\ref{ch15.transformacao_cu_NC_delay}) and (\ref{ch15.transformacao_cu_NC}), we conclude the same results in the norm
\begin{eqnarray}\label{ch15.periodic_solution_NCG_infinite_dimension}
\left(\sum_{i=1}^{2}\left[\vartheta_{i}^{\rm{av}}(t)\right]^2  + \int_{0}^{D_i} \left[u_{i}^{\rm{av}}(x,t)\right]^2dx \right)^{1/2}\,. \end{eqnarray}
%


\textcolor{black}{After applying the  averaging theory for infinite dimensional systems in \textcolor{black}{\cite[Sec.~2]{HL:1990}} (see ``Averaging Theorem for General Infinite-Dimensional Systems''), we can show that (\ref{ch15.periodic_solution_NCG}) and (\ref{ch15.limsup1_NCG}) are indeed satisfied. The final residual sets for the errors $\theta_i(t)-\theta_i^*$ in (\ref{ch15.limsup2}) and (\ref{ch15.limsup3}) depend on $a_1$ and $a_2e^{D_2\sqrt{\frac{\omega}{2}}}$ due to the amplitude of the additive dithers $S_i(t)$ in
(\ref{ch15.eq:Si}), for $i\in\{1\,,2\}$.} 


\section{Averaging Theorem for Functional Differential Equations \cite{HL:1990}}
Consider the  delay system %
\begin{eqnarray}\label{delay_FDE}
\dot{x}(t)&=&f(t/\epsilon,x_t)\,, \quad \forall t \geq 0\,,
\end{eqnarray}
where  $\epsilon$ is a
real parameter, $x_t(\Theta) = x (t+\Theta)$ for $-r\leq \Theta \leq 0$, and $f : \mathbb{R}_{+} \times \Omega \to \mathbb{R}^n$ is a continuous functional from a neighborhood $\Omega$ of $0$ of the supremum-normed Banach space $X = C([-r, 0]; \mathbb{R}^n)$ of continuous functions from $[-r, 0]$ to $\mathbb{R}^n$. Assume that $f(t,\varphi)$ is periodic in $t$ uniformly with respect to $\varphi$ in compact subsets of $\Omega$ and that $f$ has a continuous Fr\'{e}chet  derivative $\partial f (t,\varphi)/\partial \varphi$ in $\varphi$ on $\mathbb{R}_{+} \times \Omega$. If $y = y_0\in \Omega$ is an exponentially stable equilibrium for the average system 
%
\begin{eqnarray}\label{delay_FDE_average}
\dot{y}(t)&=&f_0(y_t)\,, \quad \forall t\geq 0\,,
\end{eqnarray}
where $f_0(\varphi)=\lim_{T\to \infty}\frac{1}{T} \int_{0}^{T} f(s,\varphi) ds$,
then, for some $\epsilon_0 > 0$ and $0 <\epsilon \leq \epsilon_0$, there is 
a unique  periodic solution $t \mapsto x^*(t,\epsilon)$ of \textcolor{black}{\textnormal{(\ref{delay_FDE})}} with the properties of being continuous in $t$ and $\epsilon$, satisfying  
$|x^*(t, \epsilon) - y_0| \leq \mathcal{O}(\epsilon)$, for $t \in \mathbb{R}_{+}$,
and such that there is $\rho>0$ so that
%
if $x(\cdot;\varphi)$ is a solution of \textcolor{black}{\textnormal{(\ref{delay_FDE})}} with $x(s) = \varphi$ and $|\varphi - y_0| < \rho$,
then $|x(t)-x^*(t,\epsilon)| \leq C e^{-\gamma(t-s)}$,
for $C>0$ and $\gamma>0$.

\section{Averaging Theorem for General Infinite-Dimensional Systems \cite{HL:1990}}

Consider the infinite-dimensional system, defined in the Banach space $\mathcal{X}$
\begin{align}
\dot{z} = \mathcal{A}z+J(\omega t, z)
\label{apx:inf_sys}
\end{align}
with $z(0)=z_0 \in \mathcal{X}$ and the operator $\mathcal{A}:\mathcal{D}(\mathcal{A}) \rightarrow \mathcal{X}$ generates an analytic semigroup. Moreover the nonlinearity $J:\mathbb{R}_+ \times \mathcal{X} \rightarrow \mathcal{X}$ with $t \mapsto J(\omega t,z)$ is Fr\'{e}chet differentiable in $z$, strongly continuous and periodic in $t$ uniformly with respect to $z$ in a compact subset of $\mathcal{X}$. Along with \eqref{apx:inf_sys}, the average system 
\begin{align}
\dot{z}_{\text{av}} = \mathcal{A}z_\text{av}+J_0(z_{\text{av}})
\label{apx:avg_sys}
\end{align}
with $J_0(z_{\text{av}}) = \lim\limits_{T\rightarrow \infty} \frac{1}{T} \int_{0}^{T}J(\tau,z_\text{av})d\tau$ is considered. Suppose that $z_{av} = 0 \in D \subset \mathcal{X}$ is an exponentially stable equilibrium point of the average system \eqref{apx:avg_sys}. Then for some $\bar{\omega} > 0$ and $\omega > \bar{\omega}$, we have the following: 
\begin{itemize}
	\item[a)] there exists a unique exponentially stable periodic solution $t \mapsto \bar{z}(t,1/\omega)$, continuous in $t$ and $1/\omega$, with $\|\bar{z}(t,1/\omega)\| \le \mathcal{O}(1/\omega)$ for $t>0$;
	\item[b)] with $\|z_0 - z_{\text{av}}(0)\| \le \mathcal{O}(1/\omega)$, the solution estimate of \eqref{apx:inf_sys} is given by 
	\begin{align}
	\|z(t)-z_{\text{av}}\|\le \mathcal{O}(1/\omega),\quad t>0;
	\end{align}
	\item[c)] for $\|z_0\| \le \mathcal{O}(1/\omega)$, and the stable manifold theorem, it holds 
	\begin{align}
	\|z(t)-\bar{z}(t,1/\omega)\| \le Ce^{-\gamma t},\quad t>0,
	\end{align}
	for some $C,\ \gamma>0$.
\end{itemize}

\section{Small-Gain Theorem for ODE and Hyperbolic PDE Loops \cite[Theorem~8.1, p.~198]{KK:2018}}
Consider generalized solutions of the following initial-boundary value problem
\begin{align}
& \dot{x}(t)=F(x(t)\,,u(z,t)\,,v(t))\,, \quad \forall t \geq 0\,, \label{caralho_1} \\
& u_t(z,t)+cu_z(z,t)=a(z)u(z,t)+g(z\,,x(t)\,,u(z,t))+ \nonumber \\ 
& f(z,t)\,, \quad \forall (z,t) \in  [0\,,1] \times \mathbb{R}_{+}\,, 
\label{caralho_2} \\
& u(0,t)=\varphi(d(t)\,, u(z,t)\,,x(t))\,, \quad \forall t \geq 0\,, \quad u(z,0)=u_0\,, \nonumber\\ 
& \quad x(0)=x_0\,. \label{caralho_3}
\end{align}
The state of the system \textcolor{black}{\textnormal{(\ref{caralho_1})--(\ref{caralho_3})}} is $(u(z,t),x(t))\in C^{0}([0,1]\times \mathbb{R}_{+}) \times \mathbb{R}^n$, while \textcolor{black}{the other variables} $d\in C^{0}(\mathbb{R}_+;\mathbb{R}^q)$,
$f\in C^{0}([0\,,1] \times \mathbb{R}_+)$ and
$v\in C^{0}(\mathbb{R}_+\,;\mathbb{R}^m)$ are external inputs. We assume that $(0\,,0) \in C^{0}([0\,,1])\times \mathbb{R}^n$ is an equilibrium point for the input-free system, \textit{i.e.}, $F(0\,,0\,,0)=0$,
$g(z\,,0\,,0)=0$, and $\varphi(0\,,0\,,0)=0$. Now, we assume that the ODE subsystem satisfies the ISS property:\\
\textbf{(H1)} There exist constants $M\,, \sigma>0$, $b_3\,, \gamma_3\geq0$, such that for every $x_0\in\mathbb{R}^n$, $u\in C^{0}([0\,,1] \times \mathbb{R}_{+})$ and $v\in C^0(\mathbb{R}_{+}\,;\mathbb{R}^{m})$ the unique solution $x \in C^{1}(\mathbb{R}_{+}\,;\mathbb{R}^{n})$ of (\ref{caralho_1}) with $x(0)=x_0$ satisfies the following estimate, $ \forall t \geq 0$:
\begin{align}
|x(t)| \leq M |x_0| \exp(-\sigma t) + 
\max_{0\leq s \leq t}(\gamma_3 \|u(s)\|_{\infty}+b_3|v(s)|)\,.
\end{align}
We next need to estimate the static gain of the interconnections. To this purpose, we employ the following further assumption.\\
\textbf{(H2)} There exist constants
$b_2\,,\gamma_1\,,\gamma_2\,,A\,,B \geq 0$ such that the following growth conditions hold for every $x\in C^{1}(\mathbb{R}_{+};\mathbb{R}^{n})$,
$u \in C^{0}([0\,,1]\times\mathbb{R}_{+})$ and $d \in C^{0}(\mathbb{R}_{+};\mathbb{R}^{q})$:
\begin{eqnarray}
|g(z\,,x\,,u)| \!\!&\leq&\!\! A \|u\|_{\infty} + \gamma_1 |x|\,, \quad
\forall z \in [0\,,1]\,,\\
|\varphi(d\,,u\,,x)| \!\!&\leq&\!\! B \|u\|_{\infty} + \gamma_2 |x|
+ b_2 |d|\,.
\end{eqnarray}
Let $c>0$\footnote{If the scalar $c<0$ is considered, the direction of convection must be reversed such that the boundary $u(0,t)$ is replaced by $u(1,t)$ and vice versa.} be a given constant and $a \in C^{0}([0\,,1])$ be a given function. \textcolor{black}{Consider the mappings as} $F:\mathbb{R}^n \times C^{0}([0\,,1]) \times \mathbb{R}^m \to \mathbb{R}^n$, $g:[0\,,1] \times \mathbb{R}^n \times C^{0}([0\,,1]) \to \mathbb{R}$, $\varphi: \mathbb{R}^q \times C^{0}([0\,,1]) \times \mathbb{R}^n \to \mathbb{R}$ being continuous mappings with $F(0\,,0\,,0)=0$ for which there exist constants $L>0$,
$\bar{N}\in \textcolor{black}{\lbrack 0\,,1 \lbrack}$ such that the inequalities $\max_{0\leq z \leq 1}(|g(z\,,x\,,u)-g(z\,,y\,,w)|)+|F(x\,,u\,,v)-F(y\,,w\,,v)|\leq 
L|x-y|+L\|u-w\|_{\infty}$, $|\varphi(d\,,u\,,x)-\varphi(d\,,w\,,y)|\leq
\bar{N}|x-y|+N\|u-w\|_{\infty}$, hold for all $u\,, w \in C^{0}([0\,,1])$, $x\,,y \in \mathbb{R}^n$, $v\in\mathbb{R}^{m}$, $d\in\mathbb{R}^{q}$. Suppose that assumptions \textbf{(H1)} and \textbf{(H2)} hold and that the following small-gain condition is satisfied:
\begin{small}
\begin{align}
&(\gamma_1 \gamma_3\!+\!A)c^{-1} \max_{0\leq z \leq 1}\left(p(z) \int_{0}^{z}\!\frac{1}{p(l)}dl \right)\!+\!(\gamma_2 \gamma_3\!+\!B)\max_{0\leq z \leq 1}(p(z)) \nonumber \\
&\!+\!2 \sqrt{(\gamma_1 \gamma_3\!+\!A)c^{-1}(\gamma_2 \gamma_3\!+\!B)\!\!\max_{0\leq z \leq 1}\!\!(p(z))\!\! \max_{0\leq z \leq 1}\!\left(p(z)\int_{0}^{z}\!\frac{1}{p(l)}dl \right)} \nonumber \\
&\!<\!1\,,
\end{align}
\end{small}
\!with $p(z)\!:=\!\exp\left(c^{-1}\int_{0}^{z}a(w)dw\right)$ for $z\!\in\![0\,,1]$ \textcolor{black}{[recall (8.2.11) and (8.2.14)] in \cite[Section~8.2]{KK:2018}.} Then, there exist constants $\delta\,, \Theta\,, \gamma\!>\!0$ such that for every $u_0\!\in\! C^{0}([0\,,1])$, $x_0 \!\in\! \mathbb{R}^{n}$, $d \!\in\! C^{0}(\mathbb{R}_{+}\,;\mathbb{R}^{q})$ with $u_0(0)\!=\!\varphi(d(0)\,,u_0\,,x_0)$, $f \!\in\! C^{0}([0\,,1] \times \mathbb{R}_{+})$, and $v \!\in\! C^{0}(\mathbb{R}_{+}\,;\mathbb{R}^{m})$ the unique generalized solution of the initial-boundary value problem
(\ref{caralho_1}), (\ref{caralho_2}), (\ref{caralho_3}) satisfies the following estimate:
\begin{align}
& |x(t)|+\|u(t)\|_{\infty} \leq \Theta(|x_0|+\|u_0\|_{\infty})\exp(-\delta t) \nonumber\\
& +\gamma\left[\!\max_{0\leq s\leq t}(|v(s)|) \!+\! \max_{0\leq s\leq t}(\|f(s)\|_{\infty}) \!+\! \max_{0\leq s\leq t}(|d(s)|) \right], \!\!\!\!\quad\! \nonumber \\ 
& \forall t \geq 0. 
\end{align}
Analogous small-gain results for parabolic PDE-ODE loops and parabolic-hyperbolic PDE loops can be found in \cite[Theorem~8.2, p.~205]{KK:2018} and \cite[Theorem~11.2(11.5), p.~269(277)]{KK:2018}, respectively. \textcolor{black}{Such theorems are used in the Appendix for the proofs of Theorems~\ref{ch16.theorem.16.1.vacina}, \ref{scheme2} and \ref{duoply_heterogeneous_queirozBday}.}

\section{Notation, Norms and Terminology}
The $2$-norm of the state vector $X(t)$ for a finite-dimensional system described by an ODE is denoted by single bars, $|X(t)|$.
In contrast, norms of functions (of $x$) are denoted by double bars. By default, $\|\cdot\|$ denotes the spatial $L_2[0,D]$-norm, \textit{i.e.},
$\|\cdot\|=\|\cdot\|_{L_2[0,D]}$, where we drop the index $L_2([0,D])$ if not otherwise specified. Since the state variable $u(x,t)$ of the infinite-dimensional system governed by a PDE is a function of two arguments, 
we should emphasize that taking a norm in one of the variables makes the norm a function of the other variable, as adopted
in \cite{K:2009,K:2010}. For example, the $L_2[0,D]$-norm of $u(x,t)$ in $x \in [0,D]$ is $\|u(t)\|=\left(\int_{0}^{D} u^2(x,t) dx \right)^{1/2}$, whereas the $L_\infty[0,D]$-norm is defined by $\|u(t)\|_{L_{\infty}[0,D]}=\|u(t)\|_{\infty}=\sup_{x\in[0,D]}|u(x,t)|$. Moreover, the $\mathcal{H}_1$-norm is given by $\|u(t)\|^2_{\mathcal{H}_1}=\|u(t)\|^2_{L_2}+\|u_x(t)\|^2_{L_2}$.

We denote the partial derivatives of a function $u(x,t)$ as $\partial_x u(x,t) = \partial u(x,t)/\partial x,\ \partial_t u(x,t) = \partial u(x,t)/\partial t$. We conveniently use the compact form $u_x(x,t)$ and $u_t(x,t)$ for the former and the latter, respectively.

Consider a generic nonlinear system described by $\dot{x}=f(t,x,\epsilon)$, where $x\in \mathbb{R}^n$, $f(t,x,\epsilon)$ is periodic in $t$ with
period $T$, \textit{i.e.}, $f(t+T,x,\epsilon)= f(t,x,\epsilon)$. Thus, for $\epsilon>0$ sufficiently small, we can obtain its average model given by
$\dot{x}_{\rm{av}}=f_{\rm{av}}(x)$, with $f_{\rm{av}}(x)=1/T \int_{0}^{T}f(\tau,x,0)d\tau$, where $x_{\rm{av}}(t)$ denotes the average version of the state $x(t)$ \cite{khalil1996nonlinear}.

As defined in \cite{khalil1996nonlinear}, a vector function $f(t,\epsilon)\!\in\!\mathbb{R}^n$ is said to be of order $\mathcal{O}(\epsilon)$ over an interval $[t_1,t_2]$, if $\exists k,\bar{\epsilon}: |f(t,\epsilon)| \!\le\! k\epsilon, \forall \epsilon \!\in\! [0,\bar{\epsilon}]\ \text{and}\ \forall t \in [t_1,t_2]$. In most cases we do not provide explicit values for the constants $k$ and $\bar{\epsilon}$, in which case $\mathcal{O}(\epsilon)$ can be interpreted as an order of magnitude relation for sufficiently small $\epsilon$.

The term ``$s$'' stands either for the Laplace variable or the differential operator ``$d/dt$'', according to the context. For a transfer function $H_0(s)$ with a generic input $u$, pure convolution
$h_0(t)*u(t)$, with $h_0(t)$ being the impulse response of $H_0(s)$, is also denoted by $H_0(s)u$, as done in \cite{ioannou1996robust}. 
The maximum and minimum eigenvalues of a square matrix $A$ are denoted by $\lambda_{\rm{max}}(A)$ and $\lambda_{\rm{min}}(A)$, respectively. 

\textcolor{black}{
The definition of the Input-to-State Stability (ISS) for ODE-based as well as PDE-based systems are assumed to be as provided in
\cite{Sontag_Wang:1995} and \cite{KK:2018}, respectively.}

\textcolor{black}{
Let $A\subseteq\mathbb{R}^{n}$ be an open set. By $C^{0}(A;\Omega)$, we denote the class of continuous functions on $A$, which take values
in $\Omega\subseteq\mathbb{R}^{m}$. By $C^{k}(A;\Omega)$, where $k\geq1$ is an integer, we denote the class of functions on $A\subseteq\mathbb{R}^{n}$ with continuous derivatives of order $k$, which take values in $\Omega\subseteq\mathbb{R}^{m}$. In addition,
$C([a\,,b];\mathbb{R}^{n})$ is the Banach space of continuous functions mapping the interval $[a\,,b]$ into $\mathbb{R}^{n}$, see
\cite[Chapter~2]{HL:1993}. Alternatively, $\mathcal{C}^n(\mathcal{X})$ denotes an $n$-times continuously differentiable function on the domain $\mathcal{X}$. In addition, $\mathbb{R}_+$ stands for the domain of positive real numbers including $0$.}

According to \cite{HL:1993,F:2014}, we assume the usual definitions for any delayed-system $\dot{x}(t)=f(t,x_{t})$, $t\geq t_{0}$ and
$x(t_{0}+\Theta)=\xi(\Theta)\,,  \Theta \in [-D_{\max}\,,0]$, 
where $t_{0}$ is an arbitrary initial time instant $t_{0}\geq0$, $x(t)\!\in\!\mathbb{R}^{N}$ is the state vector, $D_{\rm{max}}\!>\!0$ is the maximum time delay allowed, 
the history function of the delayed state is given by $x_{t}(\Theta)=x(t+\Theta) \in C([-D_{\max}\,,0];\mathbb{R}^{N})$, and 
the functional initial condition $\xi$ is also assumed to be continuous on $[-D_{\max}\,,0]$. Without loss of generality, we take $t_0=0$ throughout the paper.

\section{The Basic Idea of Nash Equlibrium Seeking (NES) in a Two-Player Game}




Let players P1 and P2 represent two firms that produce the same good, have dominant control over a market, and compete for profit by setting their prices $u_1$ and $u_2$, respectively. The profit of each firm is the product of the number of units sold and the profit per unit, which is the difference between the sale price and the marginal or manufacturing cost of the product. In mathematical terms, the profits are modeled by
\begin{equation} \label{Intro_z.eq:1}
    J_i(t)=s_i(t)\big(u_i(t)-m_i \big),
\end{equation}
where $s_i$ is the number of sales, $m_i$ the marginal cost, and $i \in \{1,2\}$ for P1 and P2. Intuitively, the profit of each firm will be low if it either sets the price very low, since the profit per unit sold will be low, or if it sets the price too high, since then consumers will buy the other firm's product. The maximum profit is to be expected to lie somewhere in the middle of the price range, and it crucially depends on the price level set by the other firm.

To model the market behavior, we assume a simple, but quite realistic model, where for whatever reason, the consumer prefers the product of P1, but is willing to buy the product of P2 if its price $u_2$ is sufficiently lower than the price $u_1$. Hence, we model the sales for each firm as
\begin{equation} \label{Intro_z.eq:2}
    s_1(t)=S_d-s_2(t), \quad s_2(t)=\frac{1}{p}\big ( u_1(t)-u_2(t) \big ), 
\end{equation}
where the total consumer demand $S_d$ is held fixed for simplicity, the preference of the consumer for P1 is quantified by $p>0$, and the inequalities $u_1>u_2$ and $(u_1-u_2)/p < S_d$ are assumed to hold.

Substituting (\ref{Intro_z.eq:2}) into (\ref{Intro_z.eq:1}) yields expressions for the profits $J_1(u_1,u_2)$ and $J_2(u_1,u_2)$ that are both quadratic functions of the prices $u_1$ and $u_2$, namely,
\begin{equation} \label{Intro_z.eq:4}
    J_1=\frac{-u_1^2+u_1u_2+(m_1+S_dp)u_1-m_1u_2-S_dpm_1}{p},
\end{equation}
\begin{equation} \label{Intro_z.eq:5}
    J_2=\frac{-u_2^2+u_1u_2+m_2u_1+m_2u_2}{p},
\end{equation}
and thus, the Nash equilibrium is easily determined to be
\begin{equation} \label{Intro_z.eq:6}
    u_1^\ast=\frac{1}{3}(2m_1+m_2+2S_dp)\,,
\end{equation}
\begin{equation} \label{Intro_z.eq:7}
    u_2^\ast=\frac{1}{3}(m_1+2m_2+S_dp)\,.
\end{equation}
To make sure the constraints $u_1>u_2$, $(u_1-u_2)/p < S_d$ are satisfied by the Nash equilibrium, we assume that $m_1-m_2$ lies
in the interval $(-S_dp,2S_dp)$. If $m_1=m_2$, this condition is
automatically satisfied.

For completeness, we provide here the definition of a Nash
equilibrium $u^\ast=[u_1^\ast\,, u_2^\ast]^T$ in an 2-player game:
\begin{equation} \label{Intro_z.eq:8}
    J_i(u_i^\ast,u_{-i}^\ast) \ge J_i(u_i,u_{-i}^\ast)\,, \qquad \forall u_i \in U_i,\; i \in \{1\,, 2\}\,,
\end{equation}
where $J_i$ is the payoff function of player $i$, $u_i$ its action, $U_i$ its action set, and $u_{-i}$ denotes the action of the other player. Hence, no player has an incentive to unilaterally deviate its action from $u^\ast$. 
In the duopoly example, $U_1=U_2=\mathbb{R}_+$, where $\mathbb{R}_+$ denotes the set of positive real numbers.

    \begin{figure} 
        \centering
        \includegraphics[width=0.25\textwidth]{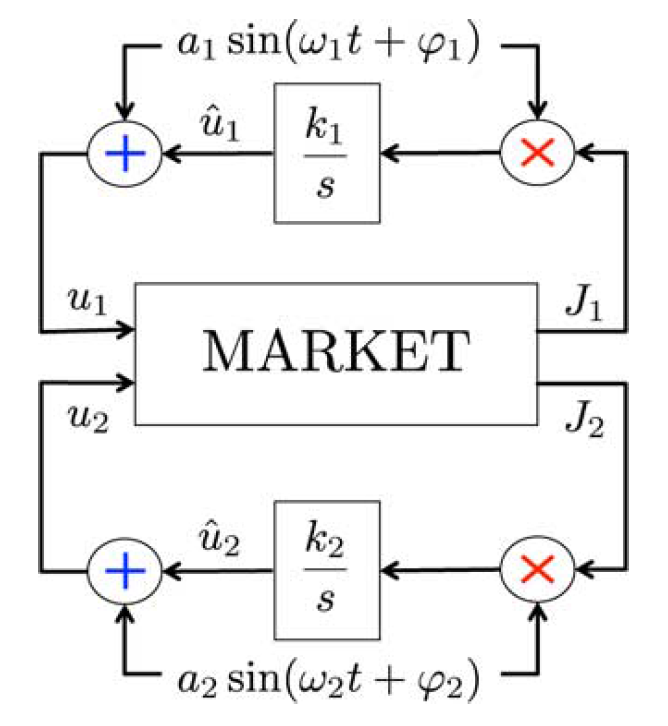}
        \caption{\hspace{-0 cm} Deterministic Nash seeking schemes applied by players in a duopoly market structure.}
        \label{fig_intro_z_Nash_Market}
    \end{figure}

To attain the Nash strategies (\ref{Intro_z.eq:6})--(\ref{Intro_z.eq:7}) without any knowledge of modeling information, such as the consumer's preference $p$, the total demand $S_d$, or the other firm's marginal cost or price, the firms implement a non-model based real-time optimization strategy, \textit{e.g.}, deterministic extremum seeking with sinusoidal perturbations, to set their price levels. Specifically, P1 and P2 set their prices, $u_1$ and $u_2$ respectively, according to the time-varying strategy \textcolor{black}{(Fig.~\ref{fig_intro_z_Nash_Market}):}
\begin{equation} \label{Intro_c.eq:00}
    \dot{\hat{u}}_i(t)=k_i\mu_i(t)J_i(t),
\end{equation}
\begin{equation} \label{Intro_c.eq:01}
    u_i(t)=\hat{u}_i(t)+\mu_i(t),
\end{equation}
where $\mu_i(t)=a_i\sin(\omega _i t + \varphi_i)$, $k_i$, $a_i$, $\omega_i > 0$ and $i \in \{1,2 \}$. Further, the frequencies are of the form
\begin{equation} \label{Intro_c.eq:02}
    \omega_i = \omega \overline{\omega}_i,
\end{equation}
where $\omega$ is a positive real number and $\overline{\omega}_i$ is a positive rational number. 

In contrast, the firms are also guaranteed to converge to the Nash equilibrium when employing the standard parallel action update scheme \cite[Proposition $4.1$]{Basar:1999}
\begin{equation} \label{Intro_z.eq:12}
    u_1^{(k+1)}=\frac{1}{2}\Big(u_2^{(k)}+m_1+S_dp \Big),
\end{equation}
\begin{equation} \label{Intro_z.eq:13}
     u_2^{(k+1)}=\frac{1}{2}\Big(u_1^{(k)}+m_2 \Big),
\end{equation}
which requires each firm to know both its own marginal cost and the other firm's price at the previous step of the iteration, and also requires P1 to know the total demand $S_d$ and the consumer preference parameter $p$. In essence, P1 must know nearly all the relevant modeling information. When using the extremum seeking algorithm (\ref{Intro_c.eq:00})--(\ref{Intro_c.eq:01}), the firms only need to measure the value of their own payoff functions, $J_1$ and $J_2$. Convergence
of (\ref{Intro_z.eq:12})--(\ref{Intro_z.eq:13}) is global, whereas the convergence of the Nash seeking strategy for this example can be proved to be semi-global, following \cite{TNM:06}, or locally, by applying the theory of averaging \cite{khalil1996nonlinear}.

\section{Fundamentals of Extremum Seeking (ES)}

Extremum seeking is a method for real-time non-model based optimization.
Though ES was introduced more than a century ago, in 1922, the ``turn of the 21st century'' has been its golden
age, both in terms of the development of theory and in terms of its adoption in industry
and in fields outside of control engineering. This section provides an overview of the basic gradient-based version
of ES with periodic signals for static maps (free of infinite-dimensional dynamics).

\begin{figure}
        \centering
        \includegraphics[width=0.4\textwidth]{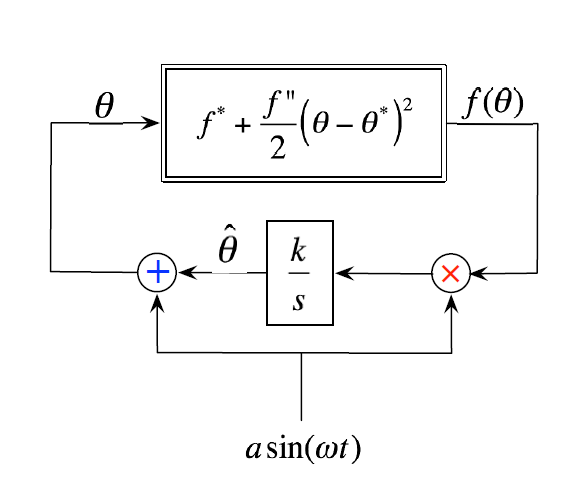}
        \caption{\hspace{0 cm} The simplest perturbation-based extremum seeking scheme for a quadratic single-input
map $f(\theta)$ in (\ref{Intro.quadratic_map}). The user has to only
know the sign of $f^{''}$, namely, whether the quadratic map has a maximum or a minimum, and
has to choose the adaptation gain $k$ such that $\sgn(k) = -\sgn\left(f^{''}\right)$. The user has to also choose
the frequency $\omega$ as relatively large compared to $a$, $k$ and $f^{''}$.}
        \label{fig:controlador_basico_es}
    \end{figure}

Many versions of ES exist, with various approaches to their study of stability 
\cite{KW:00}, \cite{LK:2012}, \cite{TNM:06}.
The most common version employs perturbation signals for the purpose of estimating
the gradient of the unknown map that is being optimized. To understand the basic idea of
extremum seeking, it is best to first consider the case of a static single-input map of the
quadratic form
\begin{equation} \label{Intro.quadratic_map}
f(\theta) = f^* + \frac{f^{''}}{2}(\theta - \theta^*)^2,
\end{equation}
where $f^*$, $f^{''}$ and $\theta^*$ are all unknown, as shown in Fig.~\ref{fig:controlador_basico_es}. 

Three different thetas appear in Fig.~\ref{fig:controlador_basico_es}:
$\theta^*$ is the unknown optimizer of the map, $\hat{\theta}(t)$
is the real-time estimate of $\theta^*$, and $\theta(t)$ is the actual input into the map. The actual input
$\theta(t)$ is based on the estimate $\hat{\theta}(t)$ but is perturbed by the signal $a \sin(\omega t)$ for the purpose
of estimating the unknown gradient $f^{''} \cdot (\theta - \theta^*)$ of the map $f(\theta)$ in (\ref{Intro.quadratic_map}).
The sinusoid is only one
choice for a perturbation signal---many other perturbations, from square waves to stochastic
noise, can be used in lieu of sinusoids, provided that they are of zero mean. The estimate $\hat{\theta} (t)$
is generated with the integrator $k/s$ with the adaptation gain $k$ controlling the speed of estimation.

The ES algorithm is successful if the error between the estimate $\hat{\theta}(t)$ and the unknown $\theta^*$, namely the signal
    \begin{equation}
        \Tilde{\theta} (t) = \hat{\theta} (t) - \theta^*,
        \label{eq:erro_de_estimacao}
    \end{equation}
converges towards zero or some small neighborhood of zero as $t \!\to\! +\infty$. Based on Fig.~\ref{fig:controlador_basico_es}, the estimate $\hat{\theta}(t)$ is governed by the differential equation $\dot{\hat{\theta}} = k \, \sin (\omega t) \, f(\theta)$, which means that the estimation error is governed by
    \begin{equation}
        \frac{d \tilde{\theta}}{dt} = 
            k \, a \, \sin (\omega t) 
            \bigg[ f^* + \frac{f^{''}}{2} 
                \big(\tilde{\theta} + a \, \sin (\omega t)
                \big)^2
            \bigg]\,.
        \label{eq:derivada_do_erro_de_estimacao}
    \end{equation}

Expanding the right-hand side, one obtains
    \begin{equation}\label{eq:slowandfastterms}
    \begin{split}
        \frac{d \tilde{\theta} (t)}{dt} =& 
            k \, a \, f^* \, \underbrace{\sin (\omega t)}_{\mbox{mean = 0}} +
            k \, a^3 \, \frac{f^{''}}{2} \underbrace{\sin^3 (\omega t)}_{\mbox{mean = 0}} \\
        &+  k \, a \, \frac{f^{''}}{2} 
                        \underbrace{ \sin (\omega t)}_{\mbox{fast, mean = 0}} \,\,\,\,\,\,\, \underbrace{\tilde{\theta} (t)^2}_{\mbox{slow}} \\
        &+  k \, a^2 \, f^{''} 
                        \underbrace{\sin^2 (\omega t)}_{\mbox{fast, mean = 1/2}} \,\,\,\,\,\,\, \underbrace{\tilde{\theta}(t)}_{\mbox{slow}}.
    \end{split}
    \end{equation}  
A theoretically rigorous time-averaging procedure \cite[Section 10.4]{khalil1996nonlinear} allows one to replace the above sinusoidal signals
by their means, yielding the ``average system''	\cite[p. 404]{khalil1996nonlinear}:
    \begin{equation}
        \frac{d \tilde{\theta}_{\rm{av}}}{dt} =
            \frac   {
                    \overbrace{k \, f^{''}}^{<0} \, 
                    a^2}
                    {2} \tilde{\theta}_{\rm{av}},
    \end{equation}
which is exponentially stable. The averaging theory guarantees that there exists a sufficiently
large $\omega$ such that, if the initial estimate $\hat{\theta}(0)$ is sufficiently close to the unknown $\theta^*$, one has
    \begin{equation}
        | \theta (t) - \theta^* | \leq 
        | \theta (0) - \theta^* | \, 
            e^{\frac{k f^{''} a^2}{2}t} + 
            \mathcal{O} \Bigg( \frac{1}{\omega} \Bigg) + \mathcal{O}(a)\,, \quad \forall t \geq 0\,.
        \label{eq:desigualdade_erro_estiamdo}
    \end{equation}
For the user, the inequality (\ref{eq:desigualdade_erro_estiamdo}) guarantees that, if $a$ is chosen small and $\omega$ is chosen
large, the input $\theta (t)$ exponentially converges to a small interval---of order $\mathcal{O}(\frac{1}{\omega} +a)$---around the unknown $\theta^*$ and,
consequently, the output $f(\theta(t))$  converges to the vicinity of the optimal output $f^*$.		

\bibliographystyle{abbrv} 
\bibliography{IEEE_Ref} 

\newpage
\begin{IEEEbiography}[{\includegraphics[width=1in,height=1.25in,clip,keepaspectratio]{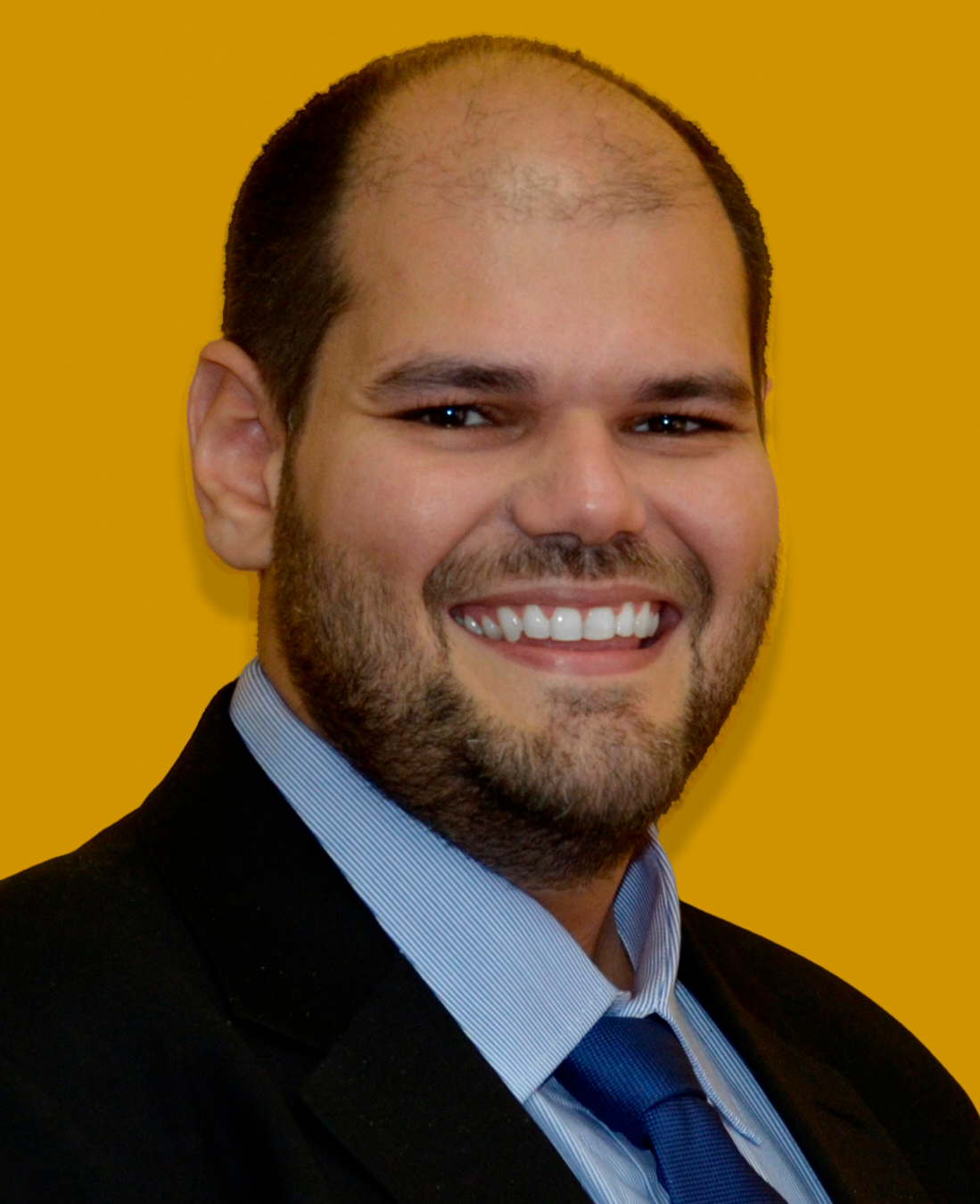}}]{Tiago Roux Oliveira} \footnotesize is a Senior Member of IEEE. He joined State University of Rio de Janeiro (UERJ), Brazil, as an Associate Professor in 2010. In 2014, he 
was a Visiting Scholar with the University of California -- San Diego (UCSD), CA, USA. 
In 2017, he was nominated as an Affiliate Member of the Brazilian Academy of Sciences (ABC). 
He has served as an Associate Editor on the Editorial Board of the Journal of the Franklin Institute, IEEE Open Journal of Control Systems, International Journal of Robust and Nonlinear Control, Systems $\&$ Control Letters, IEEE Control Systems Letters, Automatica, and other journals. In 2020, he was elected and nominated Chair of the Technical Committee 1.2 (Adaptive and Learning Systems) of the IFAC for the triennium 2020-2023 (being re-elected to the triennium 2023-2026). In 2021, he received the IEEE Transactions on Control Systems Technology Outstanding Paper Award from the IEEE CSS. He is co-author of the book entitled ``Extremum Seeking through Delays and PDEs'', published by SIAM in 2022.  
\end{IEEEbiography}

\begin{IEEEbiography}[{\includegraphics[width=1in,height=1.25in,clip,keepaspectratio]{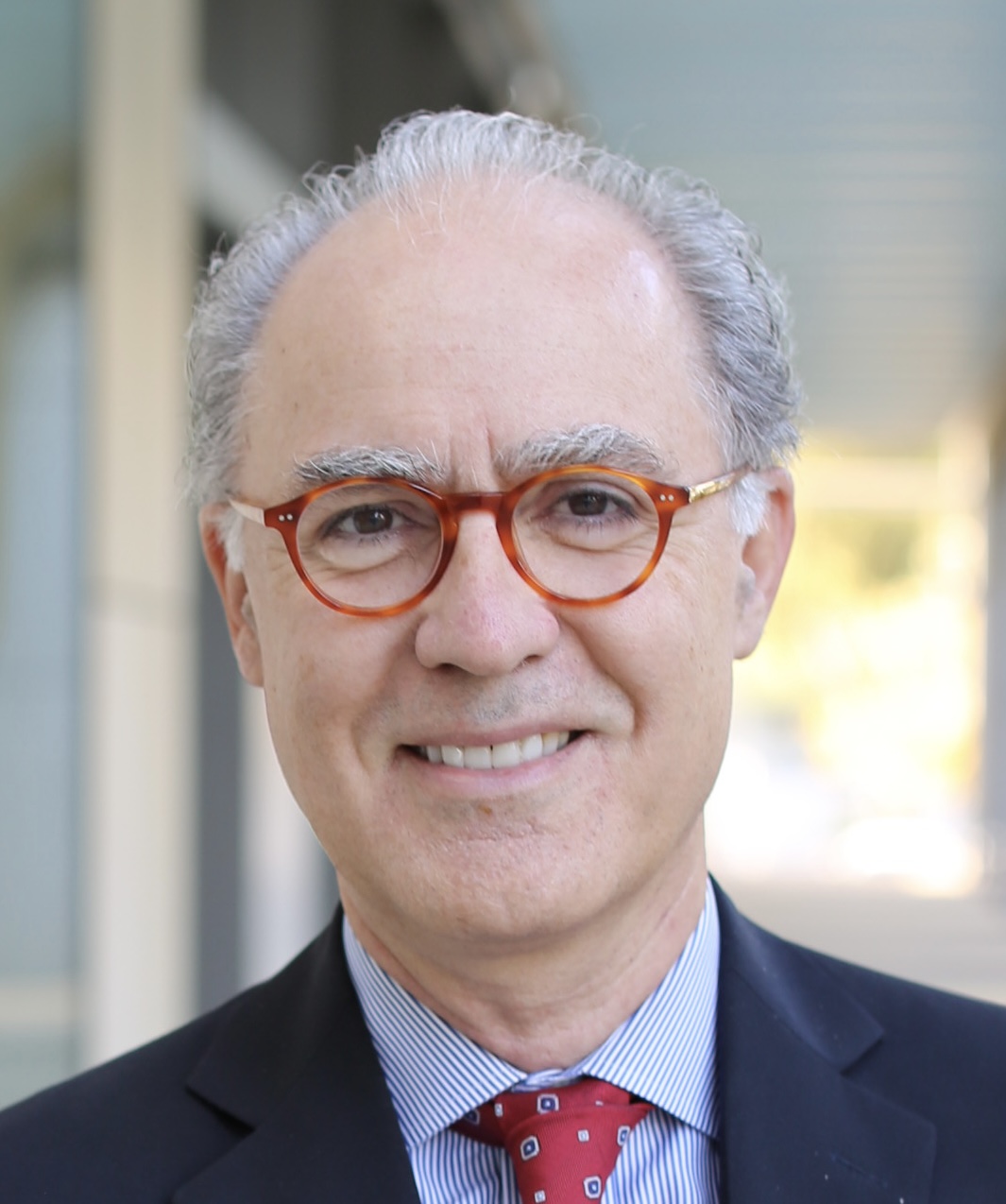}}]{Miroslav Krsti{\'c}}{\,}(krstic@ucsd.edu) is Fellow of IEEE, IFAC, SIAM, ASME, AAAS, IET, AIAA (AF), and Serbian Academy of Sciences and Arts. His awards include the Bellman, Bode Lecture, SIAM Reid, Oldenburger, Ragazzini, Chestnut, Paynter, Nyquist, IFAC Nonlinear Control, IFAC Ruth Curtain DPS, IFAC Adaptive and Learning Systems, Axelby, and Schuck (’96 and ’19). He has held chief or senior editorial positions for Systems $\&$  Control Letters, Automatica, and IEEE Transactions on Automatic Control.  Krstic has co-authored 18 books on adaptive, nonlinear, and stochastic control, extremum seeking, control of PDE systems including turbulent flows, and control of delay systems. 
\end{IEEEbiography}
\begin{IEEEbiography}[{\includegraphics[width=1in,height=1.25in,clip,keepaspectratio]{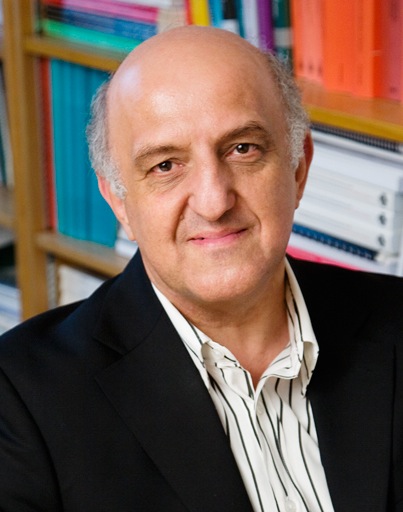}}]{Tamer Ba{\c s}ar}{\,}(basar1@illinois.edu) has been with the University
of Illinois Urbana-Champaign (UIUC), Champaign, IL, USA, since 1981, where he is currently the Swanlund Endowed Chair Emeritus; Center for Advanced Study Professor Emeritus of Electrical and Computer Engineering; a Research Professor with the Coordinated Science Laboratory; and a Research Professor with the Information
Trust Institute. He was with UIUC as the Director of the Center for Advanced Study during 2014--2020. He has more than 1000 publications in systems, control, communications, and dynamic games, including books on game theory, robust control, network security and stochastic networked control. Dr. Ba{\c s}ar was the recipient of several awards and recognitions over the years, including the highest awards of IEEE Control Systems Society (IEEE CSS), International Federation of Automatic Control (IFAC), American Automatic Control Council (AACC), and International Society of Dynamic Games (ISDG), IEEE Control Systems Award, and a number of international honorary doctorates and professorships. He is a Member of the U.S. National Academy of Engineering, a Fellow of the American Academy of Arts and Sciences, and Fellow of SIAM and IFAC, and was the President of IEEE CSS and AACC. He is a Life Fellow of IEEE. 
\end{IEEEbiography}

\end{document}